\documentclass[11pt]{amsart}
\usepackage{amsmath,amscd,amssymb,amsthm,amsfonts}

\theoremstyle{plain}
\newtheorem{thm}[equation]{Theorem}
\newtheorem{lem}[equation]{Lemma}

\newtheorem{prop}[equation]{Proposition}

\theoremstyle{definition}
\newtheorem{defi}[equation]{Definition}
\newtheorem{rmk}[equation]{Remark}

\theoremstyle{remark}

\numberwithin{equation}{section}

\newcommand{\bib}{\bibitem}
\newcommand{\ra}{\rightarrow}
\newcommand{\lra}{\longrightarrow}
\newcommand{\nin}{\noindent}
\newcommand{\bsk}{\bigskip}
\newcommand{\msk}{\medskip}

\newcommand{\contrac}{\mathop{\raise.45ex\hbox{$\underline{\hskip7pt}$}\raise.5ex\hbox{$\mkern-2mu\scriptstyle|$}}\nolimits}

\newcommand{\al}{\alpha}
\newcommand{\be}{\beta}
\newcommand{\ga}{\gamma}
\newcommand{\Ga}{\Gamma}
\newcommand{\de}{\delta}

\newcommand{\la}{\lambda}
\newcommand{\La}{\Lambda}

\newcommand{\si}{\sigma}
\newcommand{\Si}{\Sigma}
\newcommand{\om}{\omega}
\newcommand{\Om}{\Omega}
\newcommand{\ka}{\kappa}
\newcommand{\ep}{\epsilon}
\newcommand{\vep}{\varepsilon}

\newcommand{\na}{\nabla}

\newcommand{\CC}{\mathcal C}

\newcommand{\CF}{\mathcal F}

\newcommand{\CH}{\mathcal H}

\newcommand{\CL}{\mathcal L}
\newcommand{\CM}{\mathcal M}

\newcommand{\CP}{\mathcal P}

\newcommand{\CS}{\mathcal S}
\newcommand{\CT}{\mathcal T}

\newcommand{\CV}{\mathcal V}

\newcommand{\CX}{\mathcal X}

\newcommand{\CZ}{\mathcal Z}

\newcommand{\BP}{\mathbb P}
\newcommand{\BT}{\mathbb T}
\newcommand{\BR}{\mathbb R}
\newcommand{\BZ}{\mathbb Z}
\newcommand{\BC}{\mathbb C}
\newcommand{\BV}{\mathbb V}
\newcommand{\BQ}{\mathbb Q}

\newcommand{\te}{\text{e}}
\newcommand{\ti}{\text{i}}

\newcommand{\bq}{\mathbf{q}}
\newcommand{\bl}{\mathbf{l}}

\newcommand{\bw}{\mathbf{w}}
\newcommand{\bm}{\mathbf{m}}
\newcommand{\bM}{\mathbf{M}}
\newcommand{\bA}{\mathbf{A}}

\newcommand{\CLh}{\hat{\mathcal L}}
\newcommand{\CLP}{\mathcal{L}_{\mathcal{P}}}
\newcommand{\CHP}{\mathcal{H}_{\mathcal{P}}}
\newcommand{\deph}{ {\left| \text{Det} \,{\mathcal P} \right|}^{\frac{1}{2}} }
\newcommand{\dephs}{ {\left| \text{Det} \,{\mathcal P}^{*} \right|}^{\frac{1}{2}} }
\newcommand{\lag}{Lag (\mathcal V)}
\newcommand{\tlag}{\widetilde{Lag}(\mathcal V)}
\newcommand{\flag}{Lag_4 (\mathcal V)}

\begin{document}

\title{Quantization of symplectic tori in a real polarization}
\author{Mihaela Manoliu}
\address{Department of Mathematics \\
University of Texas \\
Austin, TX-78712}
\email{miha@math.utexas.edu}

\begin{abstract}
We apply the geometric quantization method with real polarizations to the quantization of a symplectic torus.
By quantizing with half-densities we canonically associate to the symplectic torus a projective Hilbert space and prove that the projective factor is expressible in terms of the Maslov-Kashiwara index.
As in the quantization of a linear symplectic space, we have two ways of resolving the projective ambiguity: (i) by introducing a metaplectic structure and using half-forms in the definition of the Hilbert space; (ii) by choosing a 4-fold cover of the Lagrangian Grassmannian of the linear symplectic space covering the torus.
We show that the Hilbert space constructed through either of these approaches realizes a unitary representation of the integer metaplectic group.
\end{abstract} 

\maketitle

\tableofcontents


\section{Introduction} \label{s:intro}

In this paper we apply the geometric quantization procedure to the quantization of a symplectic torus $(\CV/ \CZ, k \om)$, where $\CZ$ is a self-dual lattice in the symplectic vector space $(\CV, \om)$ and $k$ a positive even integer.

The basic ingredients for the geometric quantization of a symplectic manifold are a prequantum line bundle $\CL$ and a polarization $\CP$. We use invariant real polarizations of $(\CV/ \CZ, k \om)$ which correspond to 'rational' Lagrangian planes in $(\CV,\om)$.
Such a polarization $\CP$ foliates the torus by Lagrangian submanifolds and each leaf of $\CP$ has a canonically defined flat linear connection. This connection induces an operator of partial covariant differentiation along vectors in $\CP$ acting on sections of the bundle of half-densities of $\CP$. 
The Hilbert space of quantization is constructed from the space of $\CP$-parallel sections of the line bundle obtained by tensoring the prequantum line bundle $\CL$ with the half-density bundle of $\CP$.
Since the leaves of the foliation $\CP$ are compact one can only speak of distributional sections covariantly constant along these leaves. Their support defines the Bohr-Sommerfeld set in $\CV/ \CZ$.
The Hilbert space $\CHP$ constructed with the real polarization $\CP$ is a finite dimensional inner product space.

An essential problem in geometric quantization is that of comparing the Hilbert spaces of different polarizations. This requires the construction of the Blattner-Kostant-Sternberg (BKS) pairing.
In our case the pairing leads, for any two polarizations, to a unitary isomorphism between the corresponding Hilbert spaces.
Moreover, we prove that the unitary operators relating the Hilbert spaces of different polarizations satisfy a transitive composition law up to a projective factor.
The projective ambiguity is expressible in terms of the Maslov-Kashiwara index of a triple of Lagrangian subspaces in $(\CV,\om)$.
Thus quantization with half-densities canonically associates to $(\CV/ \CZ, k \om)$ a projective Hilbert space.

One can refine the above construction by introducing a metaplectic structure on $(\CV/ \CZ, k \om)$ and using half-forms instead of half-densities in the definition of the Hilbert space.
This requires an appropriate redefinition of the BKS pairing. 
As a result the projective ambiguity is resolved and we are able to canonically associate a Hilbert space to a symplectic torus $(\CV,k \om)$ with a choice of metaplectic frame bundle.

An alternative way of resolving the projective factor is by choosing a 4-fold covering space $\flag$ of the Lagrangian Grassmannian $\lag$ of $(\CV, \om)$.
The lift of the Maslov-Kashiwara index to $\flag$ can be expressed as the coboundary of the Maslov index associated to a pair of elements in $\flag$.
In this approach a Hilbert space is associated to a symplectic torus $(\CV/ \CZ, k \om)$ plus a choice of a 4-fold cover of $\lag$.

The group $Sp(\CZ)$ of symplectic transformations of $(\CV, \om)$ which preserve the lattice $\CZ$ and its double cover $Mp(\CZ)$, the integer metaplectic group, act on $(\CV/ \CZ, k \om)$ and the actions lift to the prequantum line bundle $\CL$. 
For a trivial metaplectic structure on $(\CV/ \CZ, k \om)$ the $Mp(\CZ)$-action lifts also to the metaplectic frame bundle.
We prove that the Hilbert space associated to $(\CV/ \CZ, k \om)$ and a trivial metaplectic frame bundle realizes a unitary representation of $Mp(\CZ)$.
For a choice of a nontrivial metaplectic structure only a subgroup of $Mp(\CZ)$ is represented on the respective Hilbert space.
The group $Mp(\CZ)$ acts naturally on a 4-fold cover of $\lag$.
This leads us to construct a unitary representation of $Mp(\CZ)$ also on the Hilbert space associated to $(\CV/\CZ,k \om)$ plus a choice of such 4-fold cover $\flag$.

The quantization of the torus $(\CV/ \CZ, k \om)$ is a sort of 'discrete', finite dimensional version of the quantization of the symplectic linear space $(\CV, k \om)$.
In this respect it is not surprising that in the quantization with half-densities the projective ambiguity is given in terms of the Maslov-Kashiwara triple index, similarly to the vector space case analyzed in \cite{GS,LV,Wo}.
The idea of using a covering of the Lagrangian Grassmannian of $\CV$ and the Maslov index to resolve this projective ambiguity is borrowed from \cite{GS,LV} where this procedure was shown to lead to a unitary representation of the real metaplectic group on the infinite dimensional Hilbert space of quantization of a symplectic vector space.
The use of a metaplectic structure and half-forms in the construction of the quantum Hilbert space is part of the general theory of geometric quantization \cite{Bl,GS,Wo}.

This paper arose from our work \cite{Ma} on the $U(1)$ Chern-Simons gauge theory.
There we had to address the problem of quantizing the moduli space $\CM_{\Si}$ of flat $U(1)$-connections on a Riemann surface $\Si$ using real polarizations.
We recall that $\CM_{\Si}$ is diffeomorphic to the torus $H^1(\Si;\BR)/ H^1(\Si;\BZ)$ and carries a natural symplectic structure, the push-down of the symplectic structure determined on the cohomology space $H^1(\Si;\BR)$  by the intersection pairing. 
The geometric quantization of the quotient of a symplectic vector space by a lattice using K\"{a}hler polarizations is discussed in \cite{ADW}. It is also discussed in \cite{G,R} in relation to the $U(1)$ Chern-Simons theory and the quantization of the moduli space $\CM_{\Si}$.
We have benefited from these references, as well as from \cite{JW,W} which describe the quantization of the moduli space of flat $SU(2)$-connections on a Riemann surface using a real polarization of this space.
Besides its interest in the context of the abelian Chern-Simons theory, the quantization of symplectic tori with real polarizations provides a nice example in which, like in the symplectic vector space case, the geometric quantization program (prequantization, quantization, metaplectic correction) can be fully carried through.

The organization of this paper is as follows.
In Sect.\ref{s:geomq} we give a succinct description of some aspects of geometric quantization relevant to the subsequent sections. We refer to the prequantum line bundle, real polarizations with compact leaves, definition of the Hilbert space, action of symmetries.
In Sect.\ref{s:quanttg} we apply the standard scheme of geometric quantization with half-densities to a symplectic torus $(\CV/ \CZ, k \om)$. We construct the prequantum line bundle $\CL$ and, for each choice of invariant real polarization $\CP$ of the torus, the corresponding Hilbert space $\CHP$.
In Sect.\ref{s:BKS} we describe the BKS pairing between the Hilbert spaces of any two real transverse or nontransverse polarizations. Moreover, we prove that the induced operator between the two Hilbert spaces is unitary.
In Sect.\ref{s:Heis} we describe the representation of the finite Heisenberg group on the Hilbert space $\CHP$ defined in Sect.\ref{s:quanttg}.
In Sect.\ref{s:prHsp} we prove by direct computation that the operators between the Hilbert spaces of different polarizations compose transitively up to a projective factor. 
We show that the projective ambiguity (an 8-th root of unity) can be expressed in terms of the Maslov-Kashiwara index of a triple of Lagrangian planes in the  symplectic linear space $(\CV, \om)$.
In Sect.\ref{s:metapl} we describe the metaplectic correction which enables us to construct a Hilbert space $\tilde{\CH}$.
In Sect.\ref{s:Maslov} we construct a Hilbert $\CH$ space by using a 4-fold cover of the Lagrangian Grassmannian of $(\CV,\om)$.
In Sect.\ref{s:repres} we show that the Hilbert space constructed through either approaches, of Sect.\ref{s:metapl} or of Sect.\ref{s:Maslov}, realizes a unitary representation of the integer metaplectic group.

\bsk


\section{Geometric quantization in a real polarization} \label{s:geomq}

We begin by briefly reviewing some aspects of the general scheme of geometric quantization using real polarizations with compact leaves.

The basic object is a symplectic manifold $(M , \om)$, that is, a smooth real manifold $M$ of dimension $2g$ together with a closed nondegenerate $2$-form $\om$ on $M$.
Each smooth function $f \in \CC^{\infty} (M)$ generates a vector field $ X_{f} $ on $M$, called the {\em Hamiltonian vector field\/} of $f$.
The vector field $X_{f}$ is uniquely determined by the relation $ X_{f} \, \contrac \, \om + df = 0 $, where $\contrac$ denotes the interior product.
The local flow of $X_f$ leaves the symplectic form $\om$ invariant.

The symplectic manifold $(M,\om)$ is said to be {\em 
prequantizable\/} if the cohomology class $[ \om ] \in H^{2}(M; \BR)$ defined by the symplectic form $\om$ is an integral class.
If this condition is satisfied then one can construct a {\em prequantum line bundle} for $(M, \om)$, that is a complex line bundle $\CL$ over $M$ with a hermitian metric $( \cdot , \cdot )$ on its fibres and a unitary connection $\na$ such that the curvature of $ \na$ is $-2 \pi \ti \, \om$. 
The set of isomorphism classes of such line bundles with
connections is a principal homogeneous space for the cohomology group 
$H^{1}(M;\BT)$. Hence, unless $H^{1}(M;\BT) = 0$, the choice of a prequantum line bundle $(\CL , \na)$ is not unique.

Consider an open subset $U \subset M$ such that $\CL \bigr|_{U}$ is trivializable and let $s \in \Ga(U ;\CL)$ be a nonzero section. 
Then, there exists a real $1$-form $\theta$ on $U$ such that
\begin{equation*} 
  \na s = - 2 \pi \ti \, \theta \, s
\end{equation*}
and $\om \bigr|_{U} = d \theta$. We call the $1$-form $\theta$ a local {\em symplectic potential} for $\om$ on $U$.

\msk

The quantization procedure requires the choice of a polarization, real or complex, of $(M,\om)$. We restrict ourselves to real polarizations.
A {\em real polarization} of the symplectic manifold $(M,\om)$ is a subbundle $\CP$ of the tangent  bundle $TM$, which satisfies the conditions:
(i) for each $x \in M$ , $\CP_{x}$ is a Lagrangian subspace of $T_{x} M$, that is, $\om \bigr|_{\CP_{x}} \equiv 0$ and $\dim \CP_{x} = \frac{1}{2} \, \dim T_{x} M$ ;
(ii) $\CP$ defines an integrable distribution.
Let us denote by $\CX_{\CP}(M)$ the set of all vector fields on $M$ tangent to $\CP$. Then the condition (ii) is satisfied if and only if for any $X , Y \in \CX_{\CP}(M)$ their commutator $ [X,Y] \in \CX_{\CP}(M)$.
The integral manifolds of the real polarization $\CP$ are the leaves of a foliation of $M$ by Lagrangian submanifolds.
On each leaf $\La$ of this foliation there is a canonically defined flat torsion-free linear connection \cite{We,Sn1,Wo}.
It is determined by the operator $\na^{\CP}$ of covariant differentiation  \cite{Wo}
\begin{alignat}{3} \label{e:covd}
  \na^{\CP} :  \CX_{\CP}(M)  &\times \CX_{\CP}(M) & \lra &  \;   \CX_{\CP}(M) \\
     ( X  & , Y ) & \longmapsto & \; \na^{\CP}_X \, Y  \, ,\notag
\end{alignat}
defined by setting 
\begin{equation*}
 (\na^{\CP}_X \, Y) \, \contrac \, \om = X \, \contrac \, d \, (Y \, \contrac \, \om) \, .
\end{equation*}
Let $\CC^{\infty}_{\CP}(M)$ denote the space of smooth functions on $M$ constant on the leaves of $\CP$. The Hamiltonian vector fields of functions in $\CC^{\infty}_{\CP}(M)$ belong to $\CX_{\CP}(M)$ and  are covariantly constant with respect to the connection $\na^{\CP}$.

\begin{rmk}  \label{r:angle}
Assume that the symplectic manifold $(M , \om)$ has a real
polarization $\CP$ such that the space of leaves $M/\CP$ has a manifold structure with the canonical projection map $\Pi_{\CP} : M \ra M/ \CP$ defining a smooth locally trivial fibration with compact connected fibres. 
Then each leaf $\La$ of $\CP$ is diffeomorphic to a $g$-dimensional torus $T^g$ \cite{Ar,Du,Sn1,We}.
The polarization $\CP$ is locally spanned by Hamiltonian vector fields.
If $U \subset M/\CP$ is a local coordinate neighborhood with coordinate functions $q_1, \dots ,q_g$, then the Hamiltonian vector fields $X_{q_1}, \dots ,X_{q_g}$ of the functions $q_1 \circ \Pi_{\CP}, \dots, q_g \circ \Pi_{\CP}$ commute and span $\CP \bigr|_{\Pi^{-1}_{\CP}(U)}$.
Let us recall the definition of a particular set of local coordinates on $M/\CP$ called the {\em action coordinates} \cite{Du,Ar}.
For every point in $M/ \CP$ we can choose a neighborhood $U \subset M/ \CP$ such that the symplectic form $\om$ is exact on $\Pi^{-1}_{\CP}(U)$ and such that the fibration $\Pi_{\CP}$ is trivial on $U$. Then, for each leaf $\La$ in $\Pi^{-1}_{\CP}(U)$, we can pick a basis $\ga_1(\La) , \dots , \ga_g(\La)$ for the homology group $H_1(\La;\BZ)$ such that the $\ga_i(\La)$'s vary continuously with $\La$ in $\Pi^{-1}_{\CP}(U)$.
Let $\theta$ be a symplectic potential for $\om$ on $\Pi^{-1}_{\CP}(U)$. 
The functions $j_{1} , \dots , j_{g}$ defined by
\begin{equation*} 
j_{i}(y) = \int\limits_{\ga_{i}(\La)} \theta \quad ,  \quad \text{ where } y = \Pi_{\CP}(\La) \, ,
\end{equation*}
form a local coordinate system on $U$.
The Hamiltonian vector fields $X_{j_{1}} , \dots , X_{j_{g}}$ corresponding to the action coordinates $j_{1} , \dots , j_{g}$ span $\CP \bigr|_{\Pi^{-1}_{\CP}(U)}$. 
The vector field $X_{j_{i}}$ is tangent to the curve $\ga_{i}(\La)$ and its flow is periodic with \text{period $1$}.
The distribution $\CP$ has a canonically defined density $\ka$, invariant under the Hamiltonian vector fields in $\CP$, and which assigns to each integral manifold $\La$ of $\CP$ the \text{volume $1$}.
It can be defined by setting, for each leaf $\La$ in $\Pi^{-1}_{\CP}(U)$,
\begin{equation}  \label{e:dens}
\ka \, (X_{j_{1}} \bigr|_{\La} , \dots , X_{j_{g}} \bigr|_{\La}) = 1 \, .
\end{equation}
\end{rmk}

\msk

Quantization associates to a symplectic manifold $(M ,\om)$ with real polarization $\CP$ and prequantum line bundle $(\CL , \na)$ a Hilbert space $\CHP$. The construction of $\CHP$ makes use of half-densities relative to the polarization $\CP$.
Let us first recall the definition of a density of order $p$.
Let $\BV$ be a real $n$-dimensional vector space and let $p > 0$. A {\em density of order\/} $p$ on $\BV$  is a map $ \; \nu : \overset{n}{\times} \BV \ra \BR$ such that for any $a \in GL(\BV) $ we have $\: \nu(v_1 \cdot a , \dots , v_n \cdot a ) = \left| \det \, a \right|^{p} \, \nu(v_1 , \dots , v_n) $.
We let $\text{Det} \, \BV$ stand for the highest exterior power of $\BV$, that is $\text{Det} \, \BV = \overset{n}{\wedge} \BV$.
If $\BV^{*}$ is the vector space dual to $\BV$, we let $ \left| \text{Det} \, \BV^{*} \right|^{p}$ denote the one-dimensional space of $p$-densities on $\BV$.
The above definitions and notations generalize naturally to vector bundles.

Given a polarization $\CP$, consider the line bundles $\deph$ and $\dephs$.
The canonical flat linear connection (\ref{e:covd}) defined on each leaf of $\CP$ induces operators $\na^{\CP}$ of partial covariant differentiation along vectors tangent to $\CP$ acting on sections of $\deph$, respectively $\dephs$. They are defined as follows \cite{Sn3}.
Let $X_{1} , \dots , X_{g}$ be a set of Hamiltonian vector fields spanning $\CP \bigr|_{\Pi^{-1}_{\CP}(U)}$, for some open subset $U \subset M/ \CP$. Let $X^{*}_{1} , \dots , X^{*}_{g}$ be the basis of $\CP^{*} \bigr|_{\Pi^{-1}_{\CP}(U)}$ dual to $X_{1} , \dots , X_{g}$. 
For a section $\nu$ of $\deph$, the covariant derivative $\na^{\CP}_{W} \, \nu$ is defined on $\Pi^{-1}_{\CP}(U)$ by
\begin{equation}
(\na^{\CP}_W \, \nu) (X_1^* , \dots , X_g^*)  \, = \, W( \, \nu(X_1^* , \dots , X_g^*) \, )  \, ,
\end{equation}
for any vector field $W \in \CP \bigr|_{\Pi^{-1}_{\CP}(U)} $.
Similarly, the covariant derivative $\na^{\CP}_W \mu$ of a section $\mu$ of $\dephs$ is given on $\Pi^{-1}_{\CP}(U)$ by
\begin{equation} \label{e:cdhd}
( \na^{\CP}_W \, \mu)(X_{1} , \dots , X_{g})  \, = \, W(
\, ( \mu(X_{1} , \dots , X_{g}) \, ) \, ,
\end{equation}
for any vector field $W \in \CP \bigr|_{\Pi^{-1}_{\CP}(U)}  $.

Now consider the line bundle $\CLP = \CL \otimes \deph$.
The connection $\na$ in $\CL$ and the operator $\na^{\CP}$ of covariant differentiation of sections of $\deph$ in the $\CP$-direction determine an operator $\na^{\CP}$ of covariant differentiation along vectors in $\CP$ acting on sections of $\CLP$ . It is defined by setting
\begin{equation} \label{e:cdlp}
\na^{\CP}_{X} ( s \otimes \nu) \, = \, \na_{X} s \otimes \nu + s \otimes \na_{X}^{\CP} \nu \, ,
\end{equation}
for any $X \in \CX_{\CP}(M)$ and section $ s \otimes \nu \in \Ga(M;\CLP)$.
If $\La$ is a leaf of $\CP$ then, since $\om \bigr|_{\La} \equiv 0$, the covariant derivative operator defined in (\ref{e:cdlp}) induces a flat connection on $\CLP \bigr|_{\La}$.

When the polarization $\CP$ has non-simply connected leaves one takes as the quantization of $(M, \om, \CP, \CL, \na)$ the cohomological vector space $H^{\bullet}(M;\CP,\CLP) $  of the complex \cite{Sn2,Wo}
\begin{equation*}
0 \lra  \Om^0_{\CP}(\CLP) \overset{\na^{\CP}}{\lra}  \Om^1_{\CP}(\CLP) \overset{\na^{\CP}}{\lra} \cdots \overset{\na^{\CP}}{\lra} \Om^g_{\CP}(\CLP) \lra 0
\end{equation*}
The space $\Om^k_{\CP}(\CLP)$ is the space of smooth sections of the line bundle $\overset{k}{\wedge} \CP^* \otimes \CLP$ on $M$.
The partial connection $\na^{\CP}$ on $\CLP$ extends naturally to an operator 
\begin{equation*}
\na^{\CP} \; : \; \Om^k_{\CP}(\CLP) \lra \Om^{k+1}_{\CP}(\CLP)
\end{equation*}
by
\begin{equation*}
\begin{split}
(\na^{\CP} \al) (X_1, \dots,X_{k+1}) &= \sum_{i=1}^{k+1} \, (-1)^{i+1} \, \na^{\CP}_{X_i} [ \al(X_1, \dots, \hat{X}_i, \dots, X_{k+1})] \\ 
&+ \, \sum_{i < j} \, (-1)^{i+j} \, \al ([X_i, X_j], X_1, \dots, \hat{X}_i, \dots, \hat{X}_j, \dots, X_{k+1}) \; , \\
\end{split}
\end{equation*} 
where $X_1, \dots, X_{k+1}$ are local vector fields tangent to $\CP$.
One can check that $(\na^{\CP})^2 \al =0$, for any $\al$ in $\Om^{\bullet}_{\CP}(\CLP)$.
If the polarization $\CP$ has compact leaves then \cite{Sn2} proves that the only nontrivial cohomology space is  $H^g (M;\CP,\CLP) $, where $g= \frac{1}{2} \text{dim} M$ .
Assuming the distribution $\CP$ is orientable, the space  $H^g (M;\CP,\CLP) $ can be endowed with a natural Hilbert space structure \cite{Sn2, Wo}.
In order to describe it, let us introduce the line bundle $\CL \otimes \dephs$ over $M$ obtained by tensoring the prequantum line bundle $\CL$ with the bundle of $\frac{1}{2}$-densities on $\CP$. 
The connection $\na$ in $\CL$ and the operator $\na^{\CP}$ of covariant differentiation along $\CP$ of sections in $\dephs$ define the covariant differentiation of sections of $\CL \otimes \dephs$ in the direction of vectors tangent to $\CP$:
\begin{equation} \label{e:cdlps}
\na^{\CP}_{W} ( s \otimes \mu ) \, = \, \na_{W} s \otimes \mu + s \otimes \na^{\CP}_{W} \mu \, ,
\end{equation}
for any $s \otimes \mu \in \Ga(M;\CL \otimes \dephs)$  and vector $W$ in $\CP$.
The union of the leaves $\La$ for which the holonomy group of the flat connection $\na^{\CP}$ in $( \CL \otimes \dephs) \bigr|_{\La}$ is trivial define the {\em Bohr-Sommerfeld set\/} $\mathcal{BS}_{\CP}$.
For each $\La \subset \mathcal{BS}_{\CP}$ we let $S_{\La}$ denote the one-dimensional space of covariantly constant sections of the line bundle $( \CL \otimes \dephs) \bigr|_{\La}$. 
Now, recall from (\ref{r:angle}), there is a canonical density $\ka$ on $\CP$ which gives to each leaf $\La$ of $\CP$ total \text{volume 1}. 
Given an orientation on $\CP$, the density $\ka$ extends uniquely to a $g$-form $\hat{\ka}$ on $\CP$. 
Let $\de$ denote the $\frac{1}{2}$-density on $\CP$ obtained as the square root of $\ka$ and let $\de^{-1}$ denote the dual of $\de$, that is, $\de^{-1}$ is a section of $\deph$.
For every $\La \subset \mathcal{BS}_{\CP}$ we can choose a covariantly constant section $s_{\La}$ of $\CL \bigr|_{\La}$. 
Then $s_{\La} \otimes \de^{-1}$ is a covariantly constant section of $\CLP \bigr|_{\La}$.
Let $[\al]$ be a cohomology class in $H^g (M;\CP,\CLP) $ and $\al \in \Om^g_{\CP}(\CLP)$ a representative of $[\al]$. Then we can express $\al |_{\La} = f_{\La} \, \hat{\ka} \otimes s_{\La} \otimes \de^{-1}$ for some function $f_{\La} \in \CC^{\infty}(\La)$.
The cohomology class $[\al]$ is uniquely determined by the value of the integrals $c_{\La} = \int_{\La} f_{\La} \hat{\ka}$, when $\La$ runs over the set of leaves contained in the Bohr-Sommerfeld set (\cite{Wo}, \S 10.6).
Thus, for each $\La \subset \mathcal{BS}_{\CP}$, we can associate to $[\al]$ the element $c_{\La} s_{\La} \otimes \de$ in $S_{\La}$.
This establishes a natural isomorphism
\begin{equation*}
 H^g (M;\CP,\CLP)\, \cong \, \underset{\La \subset
\mathcal{BS}_{\CP}}{\oplus} \, S_{\La}
\end{equation*}
Thus we can define the Hilbert space of quantization $\CHP$ as the vector space
\begin{equation} \label{e:hsp}
\CHP = \underset{\La \subset \mathcal{BS}_{\CP}}{\oplus} \, S_{\La}
\end{equation}
with inner product determined by
\begin{equation} \label{e:inprod}
\langle s \otimes \mu , s' \otimes \mu' \rangle  = 
  \begin{cases}
0 , &\text{ if $s \otimes \mu \in S_{\La} \, , \, s' \otimes \mu' \in S_{\La'} $ and $ \La \neq \La'$}  \\
\int\limits_{\La} (s , s') \, \mu \ast \mu' ,  &\text{ if $s \otimes \mu , \, s' \otimes \mu' \in S_{\La}$ }
  \end{cases}
\end{equation}
Here $\mu \ast \mu'$ is the density on $\La$ defined by setting, for any $x \in \La$ and basis $(W_{1} , \dots , W_{g})$ of $\CP_{x} = T_{x} \La$, 
\begin{equation*}
(\mu \ast \mu') (W_{1} , \dots , W_{g}) \, = \, \mu(W_{1} , \dots , W_{g}) \: \mu'(W_{1} , \dots , W_{g})
\end{equation*}
and $(s, s')$ is the function on $\La$ obtained by taking the hermitian inner product in the fibre of $\CL$.
If the symplectic manifold $M$ is compact, then the Hilbert space $\CHP$ is finite dimensional.
The same definition (\ref{e:hsp}) of the Hilbert space of quantization in a real polarization with compact leaves is used in \cite{JW} and it is this reference which inspired us to take this approach to quantization.

\msk

In order to have a meaningful quantization one would like to find a way of identifying the Hilbert spaces of different polarizations, at least projectively.
Let $\CP_{1}$ and $\CP_{2}$ be two polarizations of $(M , \om)$. Then, for a fixed prequantum line bundle $(\CL , \na)$, let $\CH_{\CP_1}$ and $\CH_{\CP_2}$ be the corresponding Hilbert spaces.
Under certain conditions on $\CP_1$ and $\CP_2$ there exists a sesquilinear pairing 
\begin{equation*}
\langle \langle \cdot , \cdot \rangle \rangle : \CH_{\CP_2} \times \CH_{\CP_1} \ra \BC 
\end{equation*}
called the {\em Blattner-Kostant-Sternberg} (BKS) {\em pairing}. 
The polarizations $\CP_1$ and $\CP_2$ are said to be {\em unitarily related} if there is an unitary isomorphism
\begin{equation} \label{e:fp}
 F_{\CP_2 \CP_1} : \CH_{\CP_1} \ra \CH_{\CP_2}
\end{equation} 
with
\begin{equation*}
\langle \psi_2 , F_{\CP_2 \CP_1} \, \psi_1 \rangle \, = \, \langle \langle \psi_2 , \psi_1 \rangle \rangle
, \quad \text{ for any } \, \psi_1 \in \CH_{\CP_1} , \, \psi_2 \in \CH_{\CP_2} \, .
\end{equation*}
One also requires that if $\CP_1 = \CP_2$ then $F_{\CP_2 \CP_1} = I$.  
The operator  $F_{\CP_2 \CP_1}$ is called the {\em intertwining isomorphism} for $\CH_{\CP_1}$ and $\CH_{\CP_2}$.
For arbitrary $\psi_1 , \, \psi'_1 \in \CH_{\CP_1} $ we have the sequence of equalities
\begin{align*}
\langle F_{\CP_2 \CP_1} \psi'_1 , F_{\CP_2 \CP_1} \psi_1 \rangle \,
& = \, \langle \langle F_{\CP_2 \CP_1} \psi'_1 , \psi_1 \rangle \rangle \, 
= \, \overline{ \langle \langle \psi_1 , F_{\CP_2 \CP_1} \psi'_1 \rangle \rangle} \\
& = \, \overline{ \langle \psi_1 , F_{\CP_1 \CP_2} \circ F_{\CP_2 \CP_1} \psi'_1 \rangle} \,
= \, \langle F_{\CP_1 \CP_2} \circ F_{\CP_2 \CP_1} \psi'_1 , \psi_1 \rangle \notag
\end{align*}
Thus, unitarity of $F_{\CP_2 \CP_1}$ is equivalent to the condition $ F_{\CP_1 \CP_2} \circ F_{\CP_2 \CP_1} = I$.
Two polarizations $\CP_1$ and $\CP_2$ are called {\em transverse}
if $\CP_1 \cap \CP_2 = 0$.
The formal construction of the BKS pairing for a transverse pair of real polarizations can be found in \cite{Bl,Sn3,GS,Wo}. 
The pairing construction can also be generalized to nontransverse real polarizations $\CP_1$ and $\CP_2$ for which the space $M/ \CP_{12}$ of integral manifolds of the distribution $\CP_{12} \, = \, \CP_1 \cap \CP_2$ is   a quotient manifold of $M \:$ \cite{Sn3,Wo}. 
The unitarity of the intertwining operator $F_{\CP_2 \CP_1}$ has been proved only for some particular examples.
One of these is the quantization of the flat space $(\BR^{2g} , \om)$ in a constant polarization \cite{GS,LV,Wo}. The symplectic form is $\om \, = \, \sum_{i = 1}^{g} d p^{i} \wedge d q^{i}$ in the standard coordinates $p^{1} , \dots , p^{g} , q^{1} , \dots , q^{g}$ of $\BR^{2g}$. 
The intertwining operator $F_{\CP_2 \CP_1}$ for a pair $\CP_1 , \, \CP_2$ of transverse constant real polarizations of $\BR^{2g}$ is a generalization of the Fourier transform. It turns into the ordinary Fourier transform when $\CP_1$ and $\CP_2$ are the polarizations defined by \text{$p^{i} =$ constant} and \text{$q^{i} =$ constant}, respectively. 
For three real polarizations $\CP_1 , \CP_2 , \CP_3$ the corresponding intertwining operators satisfy 
$F_{\CP_1 \CP_3} \circ F_{\CP_3 \CP_2} \circ F_{\CP_2 \CP_1} = c( \CP_1 , \CP_2 , \CP_3) \; I$,
where $c(\CP_1 , \CP_2 , \CP_3)$ is a complex number of modulus one \cite{LV,GS}.

\msk

Now let us consider a symplectic diffeomorphism $f : M \ra M$ of $(M , \om)$. Assume that $f$ can be lifted to an automorphism $\tilde{f}$ of the prequantum line bundle $(\CL , \na)$.
Let $\CP$ be a polarization of $(M , \om)$ and assume that $f$ sends $\CP$ into a polarization $f_{*} \CP$. 
Then $f$ lifts to a map of line bundles $\tilde{f} : \CL \otimes {\left| \text{Det} \,{\CP}^* \right|}^{\frac{1}{2}}  \ra \CL \otimes {\left| \text{Det} \,(f_* \CP)^* \right|}^{\frac{1}{2}}$ which gives rise to a unitary map of Hilbert spaces $\tilde{f}  : \CHP \ra \CH_{f_* \CP} $. 
If the polarizations $\CP$ and $f_* \CP$ are unitarily related, we have an intertwining isomorphism $F_{\CP , f_* \CP} : \CH_{f_* \CP} \ra \CHP$. Then we are able to associate to $f$ an unitary operator $U_{\CP}(f) \, = \, F_{\CP , f_{*} \CP} \circ \tilde{f}$ on the Hilbert space $\CHP$.
This scheme can be used sometimes to produce a unitary (projective) representation of a group $G$ acting on $(M , \om)$ by symplectic diffeomorphisms.

\bsk


\section{Quantization of tori with half-densities} \label{s:quanttg}

In this section we apply the geometric quantization scheme with half-densities, outlined in the previous section, to a symplectic torus foliated by real invariant polarizations.

Let $\CV$ be a $2g$-dimensional  real vector space with affine symplectic form $\om$ and $\CZ$ a self-dual lattice in $(\CV,\om)$. 
We recall that $\CZ$ is a self-dual lattice if $\CZ = \CZ^*$, where $\CZ^* = \{ X \in \CV \mid \om(X,W) \in \BZ, \text{ for all }  W \in \CZ \}$.
Under the quotient map $p: \CV \ra \CV/ \CZ$  the symplectic form $\om$ pushes down to a symplectic form on the torus $\CT = \CV/ \CZ$ which we continue to denote by the same letter. Notice that $\frac{\om^g}{g!}$ gives $\CT$ total \text{volume 1}.
The symplectic form $\om$  on $\CV$ is exact. We let $\theta_0$ denote the symplectic potential invariant under the group $Sp(\CV)$ of linear transformations of $\CV$ which preserve $\om$. Thus, $\om = d \theta_0$ and $\theta_0$ is defined by
\begin{equation*}
(X' \, \contrac \, \theta_{0})_{X} = \frac{1}{2} \, \om(X,X') \,  ,  \quad \text{ for any } X,  X' \in \CV .
\end{equation*}

The problem we are considering is that of quantizing  the symplectic torus $(\CT, k \, \om)$, where the normalizing factor $k$ is assumed to be a positive integer.
Later we will restrict to $k$ even.
Let us consider the trivial hermitian line bundle $\CLh = \CV \times \BC \overset{\pi}{\lra} \CV $ over $\CV$.
For each positive $k \in \BZ$ we define the connection $\hat{\na}$ on $\CLh$ by setting
\begin{equation*}
   \hat{\na} \, \hat{s}_{0} \, = \, - 2 \pi \ti k \, \theta_{0} \,   \hat{s}_{0} \, , \notag
\end{equation*}
where $\hat{s}_{0} \in \Ga(\CV;\CLh)$ is the unit section of $\CLh$, i.e. $\hat{s}_{0} (X) = (X,1)$ for all $X \in \CV$.
The curvature of the connection $\hat{\na}$ on $\CLh$ is $- 2 \pi  \ti \, k \, \om$.
The integer lattice $\CZ$ acts on $\CV$ by translations and we consider lifting this action to the line bundle $\CLh$ such that it preserves the connection and hermitian structure.
Thus we define the $\CZ$-action on $\CLh$ by \cite{Mu1}
\begin{equation}  \label{e:ZonL}
  W \cdot (X, \la) \, = \, (X+W, \ep(W) \: \te^{\pi i  k \om(W,X)}  \, \la ) , 
\end{equation}
for all  $W \in \CZ  \text{ and } (X,\la) \in \CLh$, where the map $\ep  : \CZ \ra \BT$ satisfies
\begin{equation*}
\ep(W_{1} + W_{2}) \, = \, \ep(W_{1}) \, \ep(W_{2}) \, \te^{\pi \ti k \om(W_{1}, W_{2})} , 
\end{equation*}
for any $W_{1}, \, W_{2} \in \CZ$.
Let us consider now the standard left-action of the group of symplectic transformations  $Sp(\CV)$ on the space $\CV$.
This action can be trivially lifted to $\CLh$, preserving the hermitian metric and the connection $\hat{\na}$, by setting
\begin{equation}  \label{e:SonL}
b \cdot (X, \la) \, = \, (b X, \la), 
\end{equation}
for all $ b \in Sp(\CV) \text{ and }  (X,\la) \in \CLh$.
The subgroup $Sp(\CZ)$ of $Sp(\CV)$ contains those elements which leave invariant the integer lattice $\CZ$. 
The semi-direct product $Sp(\CZ) \ltimes \CZ$ acts on $\CV$ and we consider lifting this action to $\CLh$ such that it satisfies (\ref{e:ZonL}) and (\ref{e:SonL}).
A straightforward computation \cite{G} shows that for such a lift of the $Sp(\CZ) \ltimes \CZ$-action to $\CLh$ to exist one must take $k \in 2 \BZ$ and the multiplier $\ep(W) = 1$, for all $W \in \CZ$.
With these assumptions the $Sp(\CZ) \ltimes \CZ$-action on $\CLh$ is defined by
\begin{align}  \label{e:SZonL}  
W \cdot (X, \la) & = (X + W, \te^{\pi \ti k \om(W,X)} \, \la )  \\
b \cdot (X, \la) & = ( b X , \la) ,  \notag
\end{align}
for any $ W \in \CZ$, $\, b \in Sp(\CZ)$ and $ \, (X, \la) \in \CLh$.
From now on we assume $k$ to be a positive even integer and the lift of the $\CZ$-action to $\CLh$ to be given by (\ref{e:SZonL}).

We define the prequantum line bundle $\CL$ over the symplectic torus $\CT$ to be the quotient of $\CLh$ under the $\CZ$-action as described by the commutative diagram:
\begin{equation*}
\begin{CD}
\CLh @>>> \CL = \CLh / \CZ \\
@V{\pi}VV             @VV{\pi}V \\
\CV @>>> \CT = \CV / \CZ
\end{CD}
\end{equation*}
The hermitian metric and connection $\hat{\na}$ in $\CLh$ induce naturally a hermitian metric and a connection $\na$ in $\CL$. Moreover, the $Sp(\CZ)$-action on $\CLh \ra \CV$ determines a $Sp(\CZ)$-action on $\CL \ra \CT$ preserving the hermitian structure and the connection.
We have the identification
\begin{equation*}
\Ga(\CT;\CL) \, \cong \, \Ga(\CV;\CLh)^{\CZ}
\end{equation*}
between the space $\Ga(\CT;\CL)$ of sections of $\CL \ra \CT$ and the space $\Ga(\CV;\CLh)^{\CZ}$ of $\CZ$-invariant sections of $\CLh \ra \CV$.

\msk

In order to carry out the quantization of the symplectic torus $(\CT,k \om)$ one needs to choose a polarization. 
We will consider only invariant real polarizations of $\CT$. 
Before proceeding further let us recall some definitions.
A subspace $L \subset \CV$ is called an {\em isotropic\/} subspace of the $2g$-dimensional symplectic vector space $(\CV,\om)$ if $\om(X, X') = 0$, for all $ X, X' \in L$. 
A {\em Lagrangian\/} subspace $L$ of $\CV$ is just a maximal isotropic subspace, that is $\dim L = g$. 
A subspace $\tilde{\CV}$ of $\CV$ is a {\em symplectic\/} subspace if $(\tilde{\CV}, \om \bigr|_{\tilde{\CV}})$ is a symplectic vector space.
We have the following lemma (\cite{LV}, \S 2.1):
\begin{lem} \label{l:Zbasis}
Let $\CZ$ be a self-dual lattice in $(\CV, \om)$ and $L \subset \CV$ an isotropic subspace of dimension $l \, (l \leq g)$, such that $\CZ \cap L$ generates $L$ as a vector space. 
Then there exists a symplectic basis $( W_{1} , \dots , W_{g} ; W_{1}^{\perp} , \dots , W_{g}^{\perp} )$ of $(\CV, \om)$ such that
\begin{align*}
&\CZ \,  = \, \BZ W_{1} \oplus \cdots \oplus \BZ W_{g} \oplus \BZ W_{1}^{\perp} \oplus \cdots \oplus \BZ W_{g}^{\perp} \\
& L \,  = \, \BR W_{1} \oplus \cdots \oplus \BR W_{l} \\
& \om(W_{i}, W_{j}^{\perp}) \,  = \, \de_{ij}, \; \om(W_{i}, W_{j}) \, = \, \om(W_{i}^{\perp}, W_{j}^{\perp}) \, = \, 0
\end{align*} 
\end{lem}
\nin For the proof we refer to the above mentioned reference.
Now let us consider a constant real polarization $\CP$ of the $2g$-dimensional linear space $(\CV, \om)$, defined by a foliation of $\CV$ by {\em rational} Lagrangian planes. 
The term {\em rational\/} refers to the fact that, under the identification of the tangent space $T_{x} \CV$ with $\CV$ at any point $x \in \CV$, $\CP_{x}$ corresponds to a Lagrangian subspace $L_{\CP} \subset \CV$ with the property that $\CZ \cap L_{\CP}$ generates $L_{\CP}$ as a vector space. 
The polarization $\CP$ of $(\CV , \om)$ maps under the quotient $\CV \ra \CT$ to an invariant polarization of $(\CT,\om)$ which we continue to denote by the same letter. $\CP$ is the tangent bundle along the leaves of a foliation of $\CT$ by Lagrangian submanifolds diffeomorphic to $g$-tori. 
From Lemma (\ref{l:Zbasis}) it follows that we can find a symplectic basis $( W_{1} , \dots , W_{g} ; W_{1}^{\perp} , \dots , W_{g}^{\perp} )$ of $(\CV, \om)$ such that
\begin{align}  \label{e:Zbasis}
& \CZ \,  = \, \BZ W_{1} \oplus \cdots \oplus \BZ W_{g} \oplus \BZ    W_{1}^{\perp} \oplus \cdots \oplus \BZ W_{g}^{\perp}  \\
& L_\CP \,  = \, \text{span}_{\BR} \{ W_{1} , \dots ,W_{g} \} \notag
\end{align}
We are going to use the same notation $W_i , W_i^{\perp}$ for the corresponding invariant vector fields on $\CT$.

We can choose a symplectic potential $\theta_{\CP}$ for the form $\om$ on $\CV$ adapted to the polarization $\CP$, that is, such that $W \, \contrac \, \theta_{\CP} = 0$ for any $W \in \CP$. We define $\theta_{\CP}$ by the equations
\begin{alignat}{2}
W_{i} \, \contrac \, \theta_{\CP} \, &= \, 0 &  i=1, \dots , g \\
(W_{i}^{\perp} \, \contrac \, \theta_{\CP} )_{X} \, &= \, \om(X, W_{i}^{\perp})  \qquad  & i = 1, \dots ,g  \notag
\end{alignat}
Since $d \, \theta_{\CP} =  d \, \theta_{0} = \om$, there is a real function $K_{\CP} \in \CC^{\infty}(\CV)$ such that
\begin{equation*}
\theta_{\CP} \, = \, \theta_{0} + d K_{\CP}. 
\end{equation*}
The above relation determines $K_{\CP}$ up to a constant. We fix this arbitrariness by requiring that $K_{\CP} (0) =0$. Then it is easy to check that $K_{\CP}$ satisfies the following relations:
\begin{align}  \label{e:Kp}
K_{\CP}(X+ W) - K_{\CP}(X) \, & = \quad \frac{1}{2} \om(W ,X)  \\
K_{\CP} (X+ W^{\perp}) - K_{\CP}(X) \, & = \, - \frac{1}{2} \om(W^{\perp} ,X) , \notag
\end{align}
for any $X \in \CV, \,W \in \CP$ and $W^{\perp} \in  \, \text{span}_{\BR} \{ W_1^{\perp}, \dots , W_g^{\perp} \}$. 
Let us denote by $\hat{s}_{\CP}$ the section of $\CLh \ra \CV$ defined by
\begin{equation} \label{e:sigp}
\hat{s}_{\CP}(X) \, = \, \te^{- 2 \pi \ti k \, K_{\CP}(X)} \hat{s}_{0}(X), \quad \text{ for all } \, X \in \CV.
\end{equation}
Then we have
\begin{equation} \label{e:cdsp}
\hat{\na} \hat{s}_{\CP} \, = \, - 2 \pi \ti  k \, \theta_{\CP} \, \hat{s}_{\CP}.
\end{equation}

Having chosen a polarization $\CP$ of $\CT$ we construct the Hilbert space $\CHP$ of quantization as outlined in Sect.\ref{s:geomq}.
In order to do that let us first determine the Bohr-Sommerfeld set on $\CT$. 
The Bohr-Sommerfeld set is the union of all the leaves $\La$ of the polarization $\CP$ for which the line bundle $(\CL \otimes \dephs) \bigr|_{\La}$ is trivial as a line bundle with flat connection.
From the definition of the invariant vector fields $W_{1}, \dots ,W_{g}$ spanning $\CP$ it is obvious that the canonical density $\ka$ on $\CP$ defined by (\ref{e:dens}) satisfies $\ka(W_{1}, \dots ,W_{g}) = 1$.
Then also the $\frac{1}{2}$-density $\de$ on $\CP$ defined as the square-root of $\ka$ satisfies
\begin{equation} \label{e:chdens}
\de (W_{1}, \dots ,W_{g}) \, = \, 1 .
\end{equation}
For any leaf $\La$ of $\CP$, the $\frac{1}{2}$-density $\de \bigr|_{\La}$ is a covariantly constant section of the line bundle $(\dephs) \bigr|_{\La}$ with respect to the canonical flat connection $\na^{\CP}$ defined in (\ref{e:cdhd}).
Therefore the Bohr-Sommerfeld set is determined only by the prequantum connection $\na$ in $\CL$ and it is thus given by the union of the leaves $\La$ of $\CP$ for which the holonomy of $\na$ restricted to $\CL \bigr|_{\La}$ is trivial. 
A straightforward computation shows that the Bohr-Sommerfeld orbits are determined by the condition
\begin{equation} \label{e:BScond}
\te^{2 \pi \ti  k  \om(W,X)} \, = \, 1 , \qquad \text{ for any } \,  W \in \CP \cap \CZ .
\end{equation}
With respect to the basis vectors $W_{1}, \dots ,W_{g}$ spanning $\CP \cap \CZ$ the above condition reads:
\begin{equation*}
k \, \om(W_{i} , X) \in \BZ , \qquad i = 1, \dots , g \, ,
\end{equation*}
and shows that there are $k^{g}$ distinct Bohr-Sommerfeld orbits on the torus $\CT$.
We label these Bohr-Sommerfeld orbits by $\La_{\bq}$, with $\bq \in (\BZ/ k \BZ)^g$. 
The Bohr-Sommerfeld set is then given by the union $\mathcal{BS}_{\CP} = \underset{\bq \in (\BZ/ k \BZ)^g}{\cup} \, \La_{\bq}$. 
An orbit $\La_{\bq}$ is described by the equations
\begin{equation*}
\La_{\bq = (q_{1}, \dots ,q_{g})} : \qquad k \, \om(W_{i}, X) \, = \, q_{i} \pmod{k} , \quad i=1, \dots ,g .
\end{equation*}
For each $\bq \in (\BZ/ k \BZ)^g$, the orbit $\La_{\bq}$ on $\CT$ is covered by the family $\{ \hat{\La}_{\bq , \bl} \}_{\bl = (l_{1}, \dots ,l_{g}) \in \BZ^g}$ of Lagrangian planes in $\CV$, where
\begin{equation} \label{e:Laql}
\hat{\La}_{\bq , \bl} \, = \, \{ X \in \CV \mid \, k \, \om(W_{i}, X) = q_{i} + k l_{i} , \; i=1, \dots ,g \} .
\end{equation}
We let $s_{\bq}$ denote the unitary covariantly constant section of $\CL \bigr|_{\La_{\bq}}$ represented by the $\CZ$-invariant distributional section $\hat{s}_{\bq}$ of the line bundle $\CLh \ra \CV$:
\begin{equation} \label{e:sigq}
\hat{s}_{\bq} \, = \, \underset{\bl = (l_{1}, \dots , l_{g}) \in \BZ^g}{\sum} \; \hat{s}_{\bq , \bl} \, .
\end{equation}
$\hat{s}_{\bq , \bl}$ is the section supported on the Lagrangian plane $\hat{\La}_{\bq , \bl}$, obtained by restricting to $\hat{\La}_{\bq , \bl}$ the section $\hat{s}_{\CP}$ from (\ref{e:sigp}), that is,
\begin{equation} \label{e:sigql}
\hat{s}_{\bq , \bl} (X) \, = \, \te^{ - 2 \pi \ti k K_{\CP}(X) } \hat{s}_{0}(X) , \quad \text{ for } X \in \hat{\La}_{\bq , \bl} .
\end{equation}
For each $\La_{\bq} \subset \mathcal{BS}_{\CP}$ let $\de_{\bq}$ denote $\de \bigr|_{\La_{\bq}}$. 
We can take $\si_{\bq} = s_{\bq} \otimes \de_{\bq}$ as a basis vector for the one-dimensional space $S_{\La_{\bq}}$ of covariantly constant sections of $( \CL \otimes \dephs) \bigr|_{\La_{\bq}}$. 
Then, according to (\ref{e:hsp}), the Hilbert space $\CHP$ is the $k^{g}$-dimensional vector space
\begin{equation*}
\CHP \, = \, \underset{ \bq \in (\BZ/ k \BZ)^g}{\oplus} \; S_{\La_{\bq}} \, = \, \underset{ \bq \in (\BZ/ k \BZ)^g}{\oplus} \; \BC \, \si_{\bq} \, ,
\end{equation*}
with inner product defined as in (\ref{e:inprod}).
Since $\de_{\bq} \ast \de_{\bq} = \ka \bigr|_{\La_{\bq}}$ and since $\ka \bigr|_{\La_{\bq}}$ gives $\La_{\bq}$ total volume $1$, the inner product  in $\CHP$ gives
\begin{equation*}
\langle \si_{\bq} , \si_{\bq'} \rangle \, = \, \de_{\bq \bq'}
\end{equation*}
We call $\{ \si_{\bq} = s_{\bq} \otimes \de_{\bq} \}_{\bq \in  (\BZ/ k \BZ)^g}$ a {\em standard unitary basis} of $\CHP$.
We emphasize that the choice of such a basis is uniquely determined by the choice of a symplectic basis $( W_i ; W_i^{\perp} )$ satisfying (\ref{e:Zbasis}).

\bsk


\section{The Blattner-Kostant-Sternberg pairing} \label{s:BKS}

In this section we consider two arbitrary real polarizations $\CP_1$ and $\CP_2$ of $\CT$ and construct the Blattner-Kostant-Sternberg (BKS) pairing $\langle \langle \cdot , \cdot \rangle \rangle : \CH_{\CP_2} \times \CH_{\CP_1} \lra \BC$ between the Hilbert spaces $\CH_{\CP_1}$ and $\CH_{\CP_2}$ associated to the two polarizations.
Then we prove that the operator $F_{\CP_2 \CP_1} : \CH_{\CP_1} \lra \CH_{\CP_2}$ induced from the BKS pairing is unitary.
The main reference that we use for the construction of the BKS pairing is \cite{Sn3}.

Recall that $\CH_{\CP_1} = \underset{\La_{1} \subset \mathcal{BS}_{\CP_1}}{\oplus} \; S_{\La_{1}}$ and $\CH_{\CP_2} = \underset{\La_{2} \subset \mathcal{BS}_{\CP_2}}{\oplus} \; S_{\La_2}$, where $S_{\La_i}$ is the one-dimensional vector space of covariantly constant sections of the line bundle $(\CL \otimes \left| \text{Det} \, \CP_i^* \right|^{\frac{1}{2}}) \bigr|_{\La_{i}}$ over the Bohr-Sommerfeld orbit $\La_{i} \subset \mathcal{BS}_{\CP_i}, \, i=1,2$.
The BKS pairing is defined by setting
\begin{equation*}
\langle \langle s_2 \otimes \mu_2 , s_1 \otimes \mu_1  \rangle \rangle \, = \,  \int\limits_{\La_2 \cap \La_1} \: (s_2 , s_1) \; \mu_2 \ast \mu_1 \, ,
\end{equation*}
for any $\La_1 \subset \mathcal{BS}_{\CP_1}, \, \La_2 \subset \mathcal{BS}_{\CP_2}$ and $ s_1 \otimes \mu_1 \in S_{\La_1}, \, s_2 \otimes \mu_2 \in S_{\La_2}$.
The density $\mu_2 \ast \mu_1$ on $\La_1 \cap \La_2$ is defined as follows \cite{Sn3}.
First note that $T(\La_{1} \cap \La_{2}) = T \La_{1} \cap T \La_{2}$. 
Then for any point $x \in \La_1 \cap \La_2$, we choose a basis for $T_{x}(\CT)$ of the form $(\underline{V_2} , \underline{W} ; \underline{V_1} , \underline{T})$, where
\begin{align} \label{e:pbas}
\underline{W} \, &= \, ( W_{h+1}, \dots , W_{g} ) \qquad \: \text{ is a basis of } \; T_{x} \La_{1} \cap T_{x} \La_{2} \notag \\
\underline{U_1} = (\underline{V_1} , \underline{W}) \, &= \, ( V_{11}, \dots , V_{1h}, W_{h+1}, \dots , W_g ) \; \text{ is a basis of } \, T_{x} \La_1  \\
\underline{U_2} = (\underline{V_2} , \underline{W}) \, &= \, ( V_{21}, \dots , V_{2h}, W_{h+1}, \dots , W_g ) \; \text{ is a basis of } \, T_{x} \La_2 \notag \\
  & \qquad \qquad \qquad \qquad  \qquad \, h=g - \dim(\La_{1} \cap \La_{2}) \,  , \notag
\end{align}
and such that the following relations hold:
\begin{align} \label{e:psym}
k  \om( W_i, T_j) \, &= \, \de_{ij} \qquad \qquad i,j=h+1, \dots ,g \notag \\
k  \om( V_{2i}, V_{1j}) \, &= \, \de_{ij} \qquad \qquad i,j=1, \dots ,h \\
k \om( V_{1i} ,  T_j) \, &= \, k \om( V_{2i} , T_j) \, = \, 0 .\notag
\end{align}
The density $\mu_{2} \ast \mu_{1}$ on $\La_1 \cap \La_2$ is defined by setting
\begin{equation} \label{e:pdens}
(\mu_{2} \ast \mu_{1}) (\underline{W}) \, = \, \mu_2 (\underline{U_2}) \, \mu_1 (\underline{U_1}) .
\end{equation}

\msk

The BKS pairing induces a linear operator $F_{\CP_2 \CP_1} : \CH_{\CP_1} \lra \CH_{\CP_2}$ as defined in (\ref{e:fp}). 
Then, with respect to standard unitary bases $\{ \si_{1 \bq_1} = s_{1 \bq_1} \otimes \de_{1 \bq_1} \}$ of $\CH_{\CP_1}$ and $\{ \si_{2 \bq_{2}} = s_{2 \bq_2} \otimes \de_{2 \bq_2} \}$ of $\CH_{\CP_2}$, we can represent $F_{\CP_2 \CP_1}$ by the matrix
\begin{equation}
F_{\CP_2 \CP_1} (\si_{1 \bq_1} ) \, = \, \underset{\bq_{2} \in (\BZ/ k \BZ)^g}{\sum} \:  \bM(2,1)_{\bq_{2} \bq_{1}} \; \si_{2 \bq_{2}} \, ,
\end{equation}
where
\begin{equation} \label{e:imat}
\bM(2,1)_{\bq_2 \bq_1}  \, = \, \int\limits_{\La_{2 \bq_2} \cap \La_{1 \bq_1} } \: (s_{2 \bq_2} , s_{1 \bq_1} ) \; \de_{2 \bq_2} \ast \de_{1 \bq_1} . 
\end{equation}
Having in view the definition (\ref{e:sigq}-\ref{e:sigql}) of the unitary sections $s_{2 \bq_2}$ of $\CL \bigr|_{\La_{2 \bq_2}}$ and $s_{1 \bq_1}$ of $\CL \bigr|_{\La_{1 \bq_1}}$, the function $(s_{2 \bq_2} , s_{1 \bq_1} )$ is given, at each point $x \in \La_{2 \bq_2} \cap \La_{1 \bq_1}$, by the expression
\begin{equation} \label{e:sig}
(s_{2 \bq_2} , s_{1 \bq_1} )(x) \, = \, \te^{2 \pi \ti k K_{\CP_2}(X) - 2 \pi \ti k K_{\CP_1}(X) } \, ,
\end{equation}
where $x = p(X)$ under the quotient map $p : \CV \ra \CT$. 
Using (\ref{e:Kp}) and (\ref{e:BScond}) we can show that, for any two Bohr-Sommerfeld orbits $\La_{1} \subset \mathcal{BS}_{\CP_1}$ and $\La_{2} \subset \mathcal{BS}_{\CP_2}$, the function $\te^{ 2 \pi \ti k (K_{\CP_2} - K_{\CP_1})(X) }$ on $p^{-1}(\La_1 \cap \La_2)$ is invariant under $\CZ$-translations, that is,
$ \te^{2 \pi \ti k ( K_{\CP_2} - K_{\CP_1})(X+W)} \, = \, \te^{2 \pi \ti k ( K_{\CP_2} - K_{\CP_1})(X)}$,
for any $W \in \CZ$ and $X \in p^{-1}(\La_1 \cap \La_2)$.
Thus the left-hand side of (\ref{e:sig}) is well defined.

In order to prove the unitarity of the intertwining operator $F_{\CP_2 \CP_1}$, which comes down to proving the unitarity of the matrix $\bM(2,1)$, we have to discuss separately the case of transverse and that of nontransverse polarizations.

\msk

\nin {\bf  Transverse polarizations.}
Let $\CP_1$ and $\CP_2$ be two transverse real polarizations of $(\CT , \om)$ and, as in the previous section, let $L_{\CP_1}, \, L_{\CP_2} \subset \CV$ denote the corresponding rational Lagrangian subspaces.
Then, following (\ref{l:Zbasis}), we can find symplectic bases $( W_{11} , \dots , W_{1g} ; W_{11}^{\perp} , \dots , W_{1g}^{\perp} )$ and $( W_{21} , \dots , W_{2g} ; W_{21}^{\perp} , \dots , W_{2g}^{\perp} )$ of $(\CV, \om)$, which are at the same time bases for the integer lattice $\CZ \subset \CV$ and such that
$L_{\CP_1} \, = \, \text{span}_{\BR} \{ W_{11} , \dots , W_{1g} \}$ and
$L_{\CP_2} \, = \, \text{span}_{\BR} \{ W_{21} , \dots , W_{2g} \}$.
In order to simplify the future notation let us define with respect to the  above bases the following $(g \times g)$-matrices with $\BZ$-coefficients:
\begin{align} \label{e:ommat}
\om(2,1)_{ij} &= \om(W_{2i}, W_{1j}) & \quad \om(1,2)_{ij} &= \om(W_{1i}, W_{2j})  \\
\om(2,1^{\perp})_{ij} &= \om(W_{2i}, W_{1j}^{\perp}) & \quad \om(1^{\perp},2)_{ij} &= \om(W_{1i}^{\perp}, W_{2j})   \notag \\
\om(2^{\perp},1)_{ij} &= \om(W_{2i}^{\perp}, W_{1j}) & \quad \om(1,2^{\perp})_{ij} &= \om(W_{1i}, W_{2j}^{\perp})  \notag \\
\om(2^{\perp},1^{\perp})_{ij} &= \om(W_{2i}^{\perp}, W_{1j}^{\perp}) & \quad \om(1^{\perp},2^{\perp})_{ij} &= \om(W_{1i}^{\perp}, W_{2j}^{\perp})  \notag
\end{align}
We note that
\begin{align*}
{ ^t \om(2,1)}  &=  - \om(1,2) & \quad
{ ^t \om(2,1^{\perp})}  &=  - \om(1^{\perp},2) \\
{ ^t \om(2^{\perp},1)}  &=  - \om(1,2^{\perp}) & \quad
{ ^t \om(2^{\perp},1^{\perp})}  &=  - \om(1^{\perp},2^{\perp}) 
\end{align*}
Since the polarizations $\CP_1$ and $\CP_2$ are assumed transverse, i.e. $L_{\CP_1} \cap L_{\CP_2}=0$, the vectors $( W_{11} , \dots , W_{1g}, W_{21} , \dots , W_{2g} )$ form a basis for $\CV$ and the matrix $\om(2,1)$ is nonsingular.
We will need later the following:
\begin{lem} \label{l:symm}
The matrix $\om(2^{\perp},1) \om(2,1)^{-1}$ is symmetric.
\end{lem}
\begin{proof}
If we express the vectors $W_{1i} $ with respect to the symplectic basis
$( W_{2i} ; W_{2i}^{\perp} )$ then, since $ W_{11} , \dots , W_{1g}$ span the Lagrangian subspace $L_{\CP_1}$, we find that
$ 0 = \om(W_{1i}, W_{1k}) = \om(2,1)_{ji} \om(2^{\perp},1)_{jk} -\om(2^{\perp},1)_{ji} \om(2,1)_{jk}$ and the lemma follows.
\end{proof}
Before proceeding further let us state the following lemma whose proof is straightforward:
\begin{lem} \label{l:ind}
Given a positive integer $n$ and a nonsingular $(n \times n)$-matrix $A$ with integer coefficients consider the set $\BZ^{n} / A \BZ^{n}$. Then the cardinality of the set $\BZ^{n} / A \BZ^{n}$ is equal to $\left| \det(A) \right|$.
\end{lem}
\nin The equivalence class of $\bl = (l_1, \dots ,l_g) \in \BZ^n$ in $\BZ^n / A \BZ^n$ will be denoted by $[\bl]$.

The Bohr-Sommerfeld orbits $\La_{1 \bq_{1}}$ and $\La_{2 \bq_{2}}$, $ \, \bq_{1}, \bq_{2} \in (\BZ / k \BZ)^g$, determined on the $2g$-torus $\CT$ by the two polarizations $\CP_1$ and $\CP_2$ are described by the equations:
\begin{align*}
\La_{\bq_{1} = (q_{11}, \dots ,q_{1g})} : &\qquad k \, \om(W_{1i}, X) \, = \, q_{1i} \pmod{k} , \quad i=1, \dots ,g \\
\La_{\bq_{2} = (q_{21}, \dots ,q_{2g})} : &\qquad k \, \om(W_{2i}, X) \, = \, q_{2i} \pmod{k} , \quad i=1, \dots ,g   \notag
\end{align*}
$\La_{1 \bq_{1}}$ and $\La_{2 \bq_{2}}$ are transverse submanifolds of $\CT$. 
Hence $\La_{1 \bq_{1}} \cap \La_{2 \bq_{2}}$ is a finite set and we have:

\begin{prop} \label{p:intpts}
For each $\bq_{1}, \bq_{2} \in (\BZ / k \BZ)^g$, the number of intersection points of the Bohr-Sommerfeld orbits $\La_{1 \bq_{1}}$ and $\La_{2 \bq_{2}}$ is equal to $\left| \det \om(2,1) \right|$.
\end{prop}

\begin{proof}
Let $\BZ^{g} / \om(2,1) \BZ^{g}$ be the set defined as in (\ref{l:ind}) with respect to the nonsingular integer $(g \times g)$-matrix $\om(2,1)$. 
We can use the set $\BZ^{g} / \om(2,1) \BZ^{g}$ to index the intersection points in  $\La_{1 \bq_{1}} \cap \La_{2 \bq_{2}}$. 
For a given $[\bl ] = [(l_{1}, \dots , l_{g})] \in \BZ^{g} / \om(2,1) \BZ^{g}$, the point $X_{\bl \bq_{2} \bq_{1}} \in \La_{1 \bq_{1}} \cap \La_{2 \bq_{2}}$ is determined as the solution (modulo $\CZ$) of the system of linear equations:
\begin{alignat*}{2}
\om(W_{1i}, X) \, &= \, \frac{q_{1i}}{k} ,& & \qquad i= 1, \dots , g \\
\om(W_{2i}, X) \, &= \, \frac{q_{2i}}{k} + l_{i} ,& & \qquad i= 1, \dots , g \notag
\end{alignat*}
Then, with respect to the basis $(  W_{11} , \dots , W_{1g}, W_{21} , \dots , W_{2g} )$ of $\CV$, the solution $X_{\bl \bq_{2} \bq_{1}}$ of the above system reads:
\begin{equation} \label{e:Xlqq}
X_{\bl \bq_{2} \bq_{1}} \, = \, \om(2,1)_{ij}^{-1} \, (\frac{q_{2j}}{k} + l_{j}) \, W_{1i} - \,{ ^t \om(2,1)^{-1}_{ij} } \, \frac{q_{1j}}{k} \, W_{2i} .
\end{equation}
Obviously, another $\bl' = (l'_{1}, \dots , l'_{g}) \in \BZ^{2g}$ such that $[\bl] = [\bl']$ gives (modulo $\CZ$) the same solution as $\bl = (l_{1}, \dots , l_{g})$. 
The expression (\ref{e:Xlqq}), together with the fact that the cardinality of $\BZ^{g} / \om(2,1) \BZ^{g}$ is $\left| \det \om(2,1) \right|$, proves the proposition.
\end{proof}

\nin Recall that as mentioned in (\ref{e:chdens}) the $\frac{1}{2}$-densities $\de_1 \in \left| \text{Det} \CP_1^* \right|^{\frac{1}{2}}$ and $\de_2 \in \left| \text{Det} \CP_2^* \right|^{\frac{1}{2}}$ satisfy the relations
\begin{align*}
\de_1 (  W_{11} , \dots , W_{1g}) \, &= \, 1 \\
\de_2 (  W_{21} , \dots , W_{2g}) \, &= \, 1  , \notag
\end{align*}
with respect to the bases $\underline{W_{1}} =(  W_{11} , \dots , W_{1g}) $ of $\CP_1$ and $\underline{W_{2}} =(  W_{21} , \dots , W_{2g})$ of $\CP_2$. 
Now define the frames  $\underline{U_2} = \underline{W_2}$  for $\CP_2$ and $\underline{U_1} = \underline{W_1} \cdot [k \, \om(2,1)]^{-1}$ for $\CP_1$. 
For any $x \in \La_{1 \bq_{1}} \cap \La_{2 \bq_{2}}$, we have $(\underline{U_2};\underline{U_1})$ as basis for $T_{x}(\CT)$ satisfying the conditions (\ref{e:pbas})-(\ref{e:psym}). 
Then, according to the definition (\ref{e:pdens}), we find that the density $\de_{2 \bq_2} \ast \de_{1 \bq_1}$ is the number
\begin{equation} \label{e:trdens}
\de_{2 \bq_2} \ast \de_{1 \bq_1} \, = \, \de_{2}(\underline{U_2})  \, \de_{1}(\underline{U_1}) \, = \, \left| k^{g} \, \det \om(2,1) \right|^{- \frac{1}{2}} \, .
\end{equation}
It follows then, from (\ref{e:sig}), (\ref{p:intpts}) and (\ref{e:trdens}), that the matrix (\ref{e:imat}) representing $F_{\CP_2 \CP_1}$ for transverse polarizations is given by the expression
\begin{equation*}
\bM(2,1)_{\bq_{2} \bq_{1}} \, = \, \left| k^{g} \, \det \om(2,1) \right|^{-\frac{1}{2}} \: \sum_{[\bl] \in \BZ^{g} / \om(2,1) \BZ^{g}} \, \te^{ 2 \pi \ti k ( K_{\CP_2} - K_{\CP_1})(X_{\bl \bq_{2} \bq_{1}}) } \,.
\end{equation*}
We proceed to explicitly compute $\bM(2,1)$. 
Expressing $X_{\bl \bq_{2} \bq_{1}}$, in turn, with respect to the symplectic bases $( W_{1i}\, ; W_{1i}^{\perp} )$ and $( W_{2i}\, ; W_{2i}^{\perp} )$ and making then use of (\ref{e:Kp}), we find that
\begin{align*}
K_{\CP_1} (X_{\bl \bq_2 \bq_1}) &=  \frac{1}{2} \frac{\,{ ^t \bq_1}}{k} \om(2,1)^{-1} (\frac{\bq_{2}}{k} + \bl ) - \frac{1}{2} \frac{\,{ ^t \bq_1}}{k} \om(2,1)^{-1} \om(2,1^{\perp}) \frac{\bq_{1}}{k} \notag \\
K_{\CP_2} (X_{\bl \bq_2 \bq_1})  &=  - \frac{1}{2} \frac{\,{ ^t \bq_1}}{k} \om(2,1)^{-1} (\frac{\bq_{2}}{k} + \bl ) - \frac{1}{2} \,{ ^t ( \frac{\bq_{2}}{k} + \bl)}    \om(2^{\perp},1) \om(2,1)^{-1}  (\frac{\bq_{2}}{k} + \bl) \, .
\end{align*}
This leads to the expression
\begin{equation} \label{e:trmat}
\bM(2,1)_{\bq_2 \bq_1} \, = \, \left| k^g \det \om(2,1) \right|^{- \frac{1}{2}} \; \sum_{[\bl ] \in \BZ^g / \om(2,1) \BZ^g} \, \te^{\frac{\pi \ti}{k} \bA(2,1)_{ \bl \bq_2 \bq_1} } \; ,
\end{equation}
with
\begin{equation*}
\begin{split}
\bA(2,1)_{ \bl \bq_{2} \bq_{1}} = 
 & \,{ ^t \bq_1}  \om(2,1)^{-1} \om(2,1^{\perp}) \bq_1 - 2 \,{ ^t \bq_1} \om(2,1)^{-1} ( \bq_2 + k \bl) \\
 & - \,{ ^t (\bq_2 + k \bl)} \om(2^{\perp},1)  \om(2,1)^{-1} (\bq_{2} + k \bl) \; .
\end{split}
\end{equation*}
In order to prove the unitarity of the matrix $\bM(2,1)$ we have to show that
\begin{equation} \label{e:unitmatr}
\sum_{\bq_2 \in (\BZ / k \BZ)^g} \overline{\bM(2,1)}_{\bq_2 \bq'_1} \, \bM(2,1)_{\bq_2 \bq_1} \, 
= \, \de_{\bq'_1 \bq_1} \; .
\end{equation}
It is not difficult to see from the expression (\ref{e:trmat}) that we can write
\begin{gather}  \label{e:usum}
\sum_{\bq_2 \in (\BZ / k \BZ)^g} \overline{\bM(2,1)}_{\bq_2 \bq'_1} \, \bM(2,1)_{\bq_2 \bq_1} = \\
= \frac{1}{\left| \det \om(2,1) \right|^g } \: 
\sum_{ \substack{\bq_2 \\ 0 \leq q_{2i} \leq k \left| \det \om(2,1) \right|-1 }} \:  \sum_{[\bl],[\bl']} \:  [ \te^{- \frac{\pi \ti}{k} \bA(2,1)_{ \bl' \bq_2 \bq'_1} } ] \: [ \te^{\frac{\pi \ti}{k} \bA(2,1)_{ \bl \bq_2 \bq_1} }  ]. \notag
\end{gather}
We isolate the terms in (\ref{e:usum}) containing the $\bq_2$ variable and, making use of (\ref{l:symm}), we find that their sum reads:
\begin{equation*} 
\sum_{\substack{\bq_2 \\ 0 \leq q_{2i} \leq k \left| \det \om(2,1)  \right|-1}}
 \: \te^{\frac{2 \pi \ti}{k} [ \,{ ^t (\bq'_1 - \bq_1)} + k \,{ ^t ( \bl' - \bl)} \om(2^{\perp},1) ] \om(2,1)^{-1} \bq_{2} } \, 
= \, k^g  \left| \det \om(2,1) \right|^g \, \de_{\bq'_1 \bq_1} \de_{\bl' \bl} 
\end{equation*}
Inserting the above in (\ref{e:usum}), we obtain (\ref{e:unitmatr}).

\msk

\nin {\bf Nontransverse polarizations.} 
We consider now two real polarizations $\CP_1$ and $\CP_2$ of the torus $(\CT, \om)$ which do not intersect transversally, that is $L_{12} = L_{\CP_1} \cap L_{\CP_2} \neq 0$. 
It is easy to see that the isotropic subspace $L_{12} \subset \CV$ is generated by $\CZ \cap L_{12}$ as a vector space. 
Then, using Lemma (\ref{l:Zbasis}), it follows that we can choose symplectic bases $( W_{1i} ; W_{1i}^{\perp})$ and $( W_{2i} ; W_{2i}^{\perp})$ of $(\CV, \om)$ with
\begin{align} \label{e:intbas}
W_{1i} \, = \, W_{2i} \, ,\quad W_{1i}^{\perp} \, = \, W_{2i}^{\perp} \, , \quad & \text{ for } i=h+1, \dots ,g  \\
& \quad (h= g - \dim L_{12})  \notag
\end{align}
and such that
\begin{align} \label{e:intbases}
&L_{\CP_1} \, = \, \oplus_{i=1}^{g} \BR W_{1i} 
\qquad \quad L_{\CP_2} \, = \, \oplus_{i=1}^{g} \BR W_{2i} \notag \\
&L_{12} \, = \, \oplus_{i=h+1}^{g} \BR W_{1i}  \, = \, \oplus_{i=h+1}^{g} \BR W_{2i}  \\
&\CZ \, =  \oplus_{i=1}^{g} \BZ W_{1i} \oplus \oplus_{i=1}^{g} \BZ W_{1i}^{\perp} = \oplus_{i=1}^{g} \BZ W_{2i} \oplus \oplus_{i=1}^{g} \BZ W_{2i}^{\perp} \notag 
\end{align}
With respect to the above chosen bases, the $(g \times g)$-matrices with $\BZ$-coefficients introduced in (\ref{e:ommat}) have now the particular form
\begin{alignat}{2}  \label{e:ntrom}
\om(2,1) \, &= \, 
      \begin{pmatrix}
       \check{\om}(2,1) & 0 \\
        0 & 0
       \end{pmatrix}  \qquad
&\om(2^{\perp},1^{\perp}) \, &= \,
       \begin{pmatrix}
        \check{\om}(2^{\perp},1^{\perp}) & 0 \\
         0 & 0
        \end{pmatrix} \\
\om(2^{\perp},1) \, &= \, 
         \begin{pmatrix}
        \check{\om}(2^{\perp},1) & 0 \\
         0 & -\text{I}
        \end{pmatrix}   \qquad
&\om(2,1^{\perp}) \, &= \,
          \begin{pmatrix}
        \check{\om}(2,1^{\perp}) & 0 \\
         0 & \text{I}
        \end{pmatrix} \notag
\end{alignat}
where we introduced the notation $ \check{\om}(2,1)_{ij} = \om(2,1)_{ij}$, for $ i,j=1, \dots ,h$ and similarly for the other matrices.

The Bohr-Sommerfeld orbits $\La_{1 \bq_{1}}$ and $\La_{2 \bq_{2}}$, $ \, \bq_{1}, \bq_{2} \in (\BZ/ k \BZ)^g$, corresponding to the two polarizations $\CP_1$ and $\CP_2$ of $(\CT, \om)$ are described by the equations:
\begin{align*}
\La_{\bq_{1} = (q_{11}, \dots ,q_{1g})} : &\qquad k \, \om(W_{1i}, X) \, = \, q_{1i} \pmod{k} , \quad i=1, \dots ,g \\
\La_{\bq_{2} = (q_{21}, \dots ,q_{2g})} : &\qquad k \, \om(W_{2i}, X) \, = \, q_{2i} \pmod{k} , \quad i=1, \dots ,g  \notag
\end{align*}
The intersection set  $\La_{1 \bq_{1}} \cap \La_{2 \bq_{2}}$ is nonempty  only if $q_{1i} = q_{2i}$, for \text{$i=h+1, \dots ,g$}.
If this is the case, then  $\La_{1 \bq_{1}} \cap \La_{2 \bq_{2}}$ is a submanifold of $\CT$. 
Moreover, each connected component of $\La_{1 \bq_1} \cap \La_{2 \bq_2}$ is diffeomorphic to a torus of dimension $g-h$ and we have
\begin{prop}
For each $\bq_{1}, \bq_{2} \in (\BZ/ k \BZ)^g$ such that $q_{1i} = q_{2i}$, \text{$i=h+1, \dots ,g$}, the number of connected components of $\La_{1 \bq_{1}} \cap \La_{2 \bq_{2}}$ is equal to $\left| \det \check{\om}(2,1) \right|$.
\end{prop}

\begin{proof}
Consider the set $\BZ^{h} / \check{\om}(2,1) \BZ^{h}$ of (\ref{l:ind}) with respect to the nonsingular integer $(h \times h)$-matrix $\check{\om}(2,1)$. 
We use the set $\BZ^{h} / \check{\om}(2,1) \BZ^{h}$ to index the components of the submanifold  $\La_{1 \bq_{1}} \cap \La_{2 \bq_{2}}$ of $\CT$. 
Given an $\bl = (l_{1}, \dots , l_{h}) \in \BZ^{h}$ representing an element $[\bl ] \in \BZ^{h} / \check{\om}(2,1) \BZ^{h}$, the space of solutions (modulo $\CZ$) of the system of linear equations
\begin{align} \label{e:ntrsyst}
\om(W_{1i}, X) \, &= \, \frac{q_{1i}}{k} ,  & i= 1, \dots , h \notag \\
\om(W_{1i}, X) \, &= \, \om(W_{2i}, X) \, = \, \frac{q_{1i}}{k} = \frac{q_{2i}}{k}, & i=h+1, \dots , g  \\ 
\om(W_{2i}, X) \, &= \, \frac{q_{2i}}{k} + l_{i} ,  & i= 1, \dots , h \notag
\end{align}
gives all the points belonging to one component of $\La_{1 \bq_{1}} \cap \La_{2 \bq_{2}}$. 
It is obvious that an $\bl' = (l'_{1}, \dots , l'_{h}) \in \BZ^{h}$ such that $[\bl ] = [\bl' ] $ gives the same space of solutions (modulo $\CZ$) as $\bl = (l_{1}, \dots ,l_{h})$. 
Since, according to Lemma (\ref{l:ind}), $\BZ^{h} / \check{\om}(2,1) \BZ^{h}$ has $\left| \det \check{\om}(2,1) \right|$ elements, this proves the proposition.
\end{proof}
We let $(\La_{1 \bq_{1}} \cap \La_{2 \bq_{2}})_{[\bl]}$ denote the component of $\La_{1 \bq_{1}} \cap \La_{2 \bq_{2}}$ corresponding to the element $[\bl] \in \BZ^{h} / \check{\om}(2,1) \BZ^{h}$.
According to the expression (\ref{e:imat}) and the previous observations, we have
\begin{equation} \label{e:ntrmat}
\bM(2,1)_{\bq_2 \bq_1}  \, = \, \prod_{i=h+1}^q \, \de_{q_{2i} q_{1i}}  \int\limits_{\La_{2 \bq_2} \cap \La_{1 \bq_1} } \: (s_{2 \bq_2} , s_{1 \bq_1} ) \: \de_{2 \bq_2} \ast \de_{1 \bq_1} 
\end{equation}
The function $(s_{2 \bq_2} , s_{1 \bq_1} )$ is constant on each connected component of $\La_{2 \bq_{2}} \cap \La_{1 \bq_{1}}$ and has the expression
\begin{equation} \label{e:ntrsig}
( s_{2 \bq_{2}} , s_{1 \bq_{1}} ) \bigr|_{ (\La_{2 \bq_{2}} \cap \La_{1 \bq_{1}} )_{[\bl]} } \, = \, \te^{ 2 \pi \ti (K_{\CP_2} - K_{\CP_1})(X_{\bl \bq_{2} \bq_{1}}) } \; ,
\end{equation}
where $X_{\bl \bq_{2} \bq_{1}}$ is a solution of the system (\ref{e:ntrsyst}). 
The left-hand side of (\ref{e:ntrsig}) is well-defined since, as previously remarked, $\; \te^{ 2 \pi \ti (K_{\CP_2} - K_{\CP_1})(X_{\bl \bq_{2} \bq_{1}}) }$ is invariant under a $\CZ$-translation of $ X_{\bl \bq_{2} \bq_{1}}$.

Now let us consider for any $x \in \CT$ the basis $(\underline{V_2} , \underline{W} ; \underline{V_1} , \underline{T})$ of $T_x(\CT)$, where
\begin{align} \label{e:ntrbasis}
V_{2i} &= W_{2i} & i=1, \dots , h \notag \\
W_j &= W_{1j} = W_{2j} & j=h+1, \dots , g  \\
V_{1i} &= W_{1l} \; [k \check{\om}(2,1)]^{-1}_{li} \qquad & i=1, \dots ,h \notag \\
 T_{j} &= \frac{1}{k} \,  W_{1j}^{\perp} & j=h+1, \dots ,g  \notag
\end{align}
This basis satisfies the conditions (\ref{e:pbas})-(\ref{e:psym}). 
From the definition (\ref{e:pdens}), it follows that, for any pair of Bohr-Sommerfeld orbits $\La_{1 \bq_{1}}$ and $\La_{2 \bq_{2}}$ with $\La_{1 \bq_{1}} \cap \La_{2 \bq_{2}} \neq \emptyset$, the density $\de_{2 \bq_2} \ast \de_{1 \bq_1}$ satisfies
\begin{equation*}
( \de_{2 \bq_2} \ast \de_{1 \bq_1} ) (\underline{W}) \, = \,
\de_{2}(\underline{V_2} , \underline{W}) \: \de_{1}(\underline{V_1} , \underline{W}) \, = \, \left| k^{h} \det \check{\om}(2,1) \right|^{- \frac{1}{2}} 
\end{equation*}
Then, since $\underline{W}=( W_{h+1},\dots,W_g)$ is a basis for $\CZ \cap L_{12}$, we have
\begin{equation*}
\int\limits_{(\La_{2 \bq_2} \cap \La_{1 \bq_1})_{[\bl]} } \: \de_{2 \bq_2} \ast \de_{1 \bq_1} \, = \, 
\left| k^{h} \det \check{\om}(2,1) \right|^{- \frac{1}{2}}
\end{equation*}
Therefore the matrix $\bM(2,1)$ can be written as
\begin{equation*}
\bM(2,1)_{\bq_{2} \bq_{1}} \, = \, \left| k^{h} \, \det \check{\om}(2,1) \right|^{-\frac{1}{2}} \: \prod_{i=h+1}^{g} \, \de_{q_{2i} q_{1i}} \: \sum_{[\bl] \in \BZ^{h} / \check{\om}(2,1) \BZ^{h}} \, \te^{ 2 \pi \ti k ( K_{\CP_2} - K_{\CP_1})(X_{\bl \bq_2 \bq_1}) } 
\end{equation*}
Expressing the solution $X_{\bl \bq_{2} \bq_{1}}$ of the system (\ref{e:ntrsyst}) in terms of the bases $( W_{1i}\, ; W_{1i}^{\perp} )$ and $( W_{2i}\, ; W_{2i}^{\perp} )$ of $\CV$ and using the properties (\ref{e:Kp}), we arrive at
\begin{equation} \label{e:ntrmatrix}
\bM(2,1)_{\bq_{2} \bq_{1}} \, = \, \left| k^{h} \, \det \check{\om}(2,1) \right|^{-\frac{1}{2}} \: \prod_{i=h+1}^{g} \, \de_{q_{2i} q_{1i}} \: \sum_{[\bl] \in \BZ^{h} / \check{\om}(2,1) \BZ^{h}} \, \te^{ \frac{\pi \ti}{k} \check{\bA}_{\bl \check{\bq}_{2} \check{\bq}_{1}}  } \;  ,
\end{equation}
with
\begin{equation*}
\begin{split}
\check{\bA}(2,1)_{ \bl \check{\bq}_{2} \check{\bq}_{1}} = 
 & \,{ ^t \check{\bq}_1} \check{\om}(2,1)^{-1} \check{\om}(2,1^{\perp}) \check{\bq}_{1} - 2 \,{ ^t \check{\bq}_1} \check{\om}(2,1)^{-1} ( \check{\bq}_{2} + k \bl) \\
 & - \,{ ^t (\check{\bq}_2 + k \bl)} \check{\om}(2^{\perp},1)  \check{\om}(2,1)^{-1} (\check{\bq}_2 + k \bl) \, .
\end{split}
\end{equation*}
Here we used the notation $\check{\bq}_{1} = (q_{11}, \dots, q_{1h}) , \, \check{\bq}_{2} = (q_{21}, \dots , q_{2h})$.
From this point on the proof of the fact that $\sum_{\bq_{2}}  \overline{\bM(2,1)}_{\bq'_{1} \bq_{2}} \, \bM(2,1)_{\bq_{2} \bq_{1}} \, 
= \, \de_{\bq'_{1} \bq_{1}} $ goes as in the previous case of transverse polarizations.
Thus we can state

\begin{thm}
For any two polarizations $\CP_1$ and $\CP_2$ of the symplectic torus $(\CT,k \om)$, the linear operator $F_{\CP_2 \CP_1} : \CH_{\CP_1} \ra \CH_{\CP_2}$ is unitary.
\end{thm}

\bsk


\section{Representation of the Heisenberg group} \label{s:Heis}

Let us consider again the $2g$-dimensional symplectic vector space $(\CV, \om)$ with self-dual integer lattice $\CZ$ in $\CV$.
Let $Z_{k}$ be the finite abelian group $\frac{1}{k} \CZ / \CZ$. 
We define the Heisenberg group $H_k$ as a central extension of $Z_k$ by the circle group $\BT$:
\begin{equation*}
1 \lra \BT \lra H_k \lra Z_k \lra 0 \, .
\end{equation*}
As a set $H_k  = \BT \times Z_{k}$ and the composition law in $H_{k}$ is given by
\begin{equation*}
(\la,v) \cdot (\la',v') = (\la \la' \, c(v,v') , v+v'), \quad \text{ for any } (\la,v), (\la',v') \in \BT  \times Z_k \, .
\end{equation*}
The $2$-cocycle $c : Z_{k} \times Z_{k} \ra \BT$ is defined by
\begin{equation*}
c(v,v') \, = \, \te^{- \pi \ti k \om(V,V')} \, ,
\end{equation*}
where $V,V'$ are the preimages of $v,v'$ under the quotient map $\frac{1}{k} \CZ \ra  \frac{1}{k} \CZ / \CZ$. 
As before, we work under the assumption that $k \in 2 \BZ_{+}$, so that $c$ is well-defined by the above expression.

We will show that the vector space $\CHP$ obtained by quantizing the $2g$-torus $(\CT, k \om)$ in a real polarization $\CP$ realizes a unitary representation of the Heisenberg group $H_{k}$. 

Let us choose as in (\ref{e:Zbasis}) an integer symplectic basis $( W_{i}\, ; W_{i}^{\perp} )$ for the space $(\CV, \om)$, with $\CP = \text{span}_{\BR} \{ W_{1}, \dots , W_{g} \}$, and let us introduce as in Sect.\ref{s:quanttg} the corresponding standard unitary basis $\{ s_{\bq} \otimes \de_{\bq} \}_{\bq \in (\BZ / k \BZ)^g}$ for $\CHP$. The sections $s_{\bq}$ of the prequantum line bundle $\CL$ over the Bohr-Sommerfeld orbits $\La_{\bq}$ are represented by the $\CZ$-invariant distributional sections $\hat{s}_{\bq} \in \Ga(\CV; \CLh)$ defined by (\ref{e:sigq})-(\ref{e:sigql}). 
The abelian translation group $\CV$ acts on the space $\Ga(\CV;\CLh)$ of sections of the line bundle $\CLh \ra \CV$ by
\begin{equation} \label{e:act}
(V \cdot s)(X) \, = \, V \cdot s(X-V), \quad \text{ for any } V \in \CV, \, s \in \Ga(\CV;\CLh) \, .
\end{equation}
We compute the action of the subgroup $\frac{1}{k} \CZ \subset \CV$ on the sections $\hat{s}_{\bq}$ according to the definition (\ref{e:act}). 
Using the relations (\ref{e:Kp}), we find that with respect to the basis $( \frac{1}{k} W_{i}\, ; \frac{1}{k} W_{i}^{\perp} )$ of $\frac{1}{k} \CZ$ the action of $\frac{1}{k} \CZ$  on $\{ \hat{s}_{\bq} \}$ reads 
\begin{align*}
\frac{1}{k} W_{i} \cdot \hat{s}_{\bq} \, &= \, \te^{ 2 \pi \ti \frac{q_{i}}{k} } \: \hat{s}_{\bq} \\
\frac{1}{k} W_{i}^{\perp} \cdot \hat{s}_{\bq} \, &= \, \hat{s}_{\bq_{(i+1)}} \, , \notag
\end{align*}
where $\bq_{(i+1)}$ stands for $\bq_{(i+1)} = (q_{1}, \dots ,q_{i-1}, q_{i} +1, q_{i+1}, \dots ,q_{g})$. 
Since the sections $\{ \hat{s}_{\bq} \}$ are $\CZ$-invariant, we get a well defined action of the group $Z_{k} =  \frac{1}{k} \CZ / \CZ$. 

The section $\de $ of $\dephs$ was obtained as the square root of the invariant density $\ka$ on $\CP$. 
Hence the translation group $\CV$ leaves $\de$ invariant. 
For any $X \in \hat{\La}_{\bq,\bl}$, we find from (\ref{e:Laql}) that $X + \frac{1}{k} W_{i} \in \hat{\La}_{\bq,\bl}$ and $X + \frac{1}{k} W_{i}^{\perp} \in \hat{\La}_{\bq_{(i+1)},\bl}$. Thus we can write:
\begin{align*}
\frac{1}{k} W_{i} \cdot \de_{\bq} \, &= \,  \de_{\bq} \\
\frac{1}{k} W_{i}^{\perp} \cdot \de_{\bq} \, &= \, \de_{\bq_{(i+1)}} .\notag
\end{align*}
Therefore, we define the action of the Heisenberg group $H_{k}$ on the basis $\{ s_{\bq} \otimes \de_{\bq} \}$ of the Hilbert space $\CHP$ by
\begin{align} \label{e:Hact}
(\la, w_{i}) \cdot ( s_{\bq} \otimes \de_{\bq} ) \, &= \, \la \, \te^{2 \pi \ti \frac{q_{i}}{k} } \: s_{\bq} \otimes \de_{\bq} \\
(\la, w_{i}^{\perp}) \cdot ( s_{\bq} \otimes \de_{\bq} ) \, &= \, \la \, s_{\bq_{(i+1)}} \otimes \de_{\bq_{(i+1)}} \; , \notag
\end{align}
where $\la \in \BT$ and $w_{i}, w_{i}^{\perp}$ are the images in $Z_{k}$ of the basis vectors $W_{i}, W_{i}^{\perp}$. 
For any $\la \in \BT$ the element $(\la,0) \in H_{k}$ is represented on $\CHP$ by the operator $\la \, I$. 
Hence, by the Stone-Von Neumann theorem (see \cite{Mu2}) the above representation of $H_{k}$ is, up to isomorphism, the unique irreducible unitary representation of $H_{k}$ on $\CHP$.

\begin{thm}
For any  two invariant real polarizations $\CP_1, \CP_2$ of $(\CT, k \om)$, the unitary operator $F_{\CP_2 \CP_1} : \CH_{\CP_1} \ra \CH_{\CP_2}$ intertwines the representations of the Heisenberg group $H_{k}$ on the vector spaces $\CH_{\CP_1}$ and $\CH_{\CP_2}$.
\end{thm}

\begin{proof}
To prove this statement it is enough to show that
\begin{align*}
F_{\CP_2 \CP_1} [(\la, w_{1i}) \cdot ( s_{1 \bq_{1}} \otimes \de_{1 \bq_{1}} ) ] \, &= \, (\la, w_{1i}) \cdot F_{\CP_2 \CP_1} ( s_{1 \bq_{1}} \otimes \de_{1 \bq_{1}} ) \\
F_{\CP_2 \CP_1} [(\la, w_{1i}^{\perp}) \cdot ( s_{1 \bq_{1}} \otimes \de_{1 \bq_{1}} ) ] \, &= \, (\la,w_{1i}^{\perp}) \cdot F_{\CP_2 \CP_1} ( s_{1 \bq_{1}} \otimes \de_{1 \bq_{1}} ) \, . \notag
\end{align*}
We make use of the fact that the bases $( W_{1i}; W_{1i}^{\perp} )$ and  $( W_{2i}; W_{2i}^{\perp} )$ of $(\CV, \om)$, chosen as in (\ref{e:Zbasis}), are related by
\begin{align*}
W_{1i} \, &= \, - \om(2^{\perp},1)_{ji} W_{2j} + \om(2,1)_{ji} W_{2j}^{\perp} \\
W_{1i}^{\perp} \, &= \, - \om(2^{\perp},1^{\perp})_{ji} W_{2j} + \om(2,1^{\perp})_{ji} W_{2j}^{\perp} . \notag
\end{align*}
We also use the relations:
\begin{align*}
\om(2^{\perp},1^{\perp}) &= \,{ ^t \om(2,1)^{-1} } + \om(2^{\perp},1) \om(2,1)^{-1} \om(2,1^{\perp}) , \\
& \qquad \qquad \text{if $\CP_1$ and $\CP_2$ are transverse}, \\
\check{\om}(2^{\perp},1^{\perp}) &= \,{ ^t \check{\om}(2,1)^{-1}} + \check{\om}(2^{\perp},1) \check{\om}(2,1)^{-1} \check{\om}(2,1^{\perp}) , \\
& \qquad \qquad \text{if $\CP_1$ and $\CP_2$ are not transverse}.
\end{align*}
Thus we can write:
\begin{align*}
(\la, w_{1i}) &= \bigl( \la \te^{ \pi \ti k [\,{ ^t \om(2,1)} \om(2^{\perp},1) ]_{ii} } , \sum_{j} \om(2,1)_{ji} w_{2j}^{\perp} \bigr) \cdot \bigl(1, - \sum_{j} \om(2^{\perp},1)_{ji} w_{2j} \bigr) \\
(\la, w_{1i}^{\perp}) &= \bigl(\la \te^{ \pi \ti k [\,{ ^t \om(2,1^{\perp})} \om(2^{\perp},1^{\perp}) ]_{ii} } , \sum_{j} \om(2,1^{\perp})_{ji} w_{2j}^{\perp} \bigr) \cdot \bigl(1, - \sum_{j} \om(2^{\perp},1^{\perp})_{ji} w_{2j} \bigr) .
\end{align*}
Then, using (\ref{e:Hact}), we have for $\CP_1$ and $\CP_2$ transverse:
\begin{equation*}
\begin{split}
& (\la, w_{1i}) \cdot F_{\CP_2 \CP_1} ( s_{1 \bq_1} \otimes \de_{1 \bq_1} ) \, = 
\, \sum_{\bq_2} \bM(2,1)_{\bq_2 \bq_1} (\la, w_{1i}) \cdot ( \si_{2 \bq_2} \otimes \de_{2 \bq_2} ) \\
&= \la \te^{ \pi \ti k [\,{ ^t \om(2,1)} \om(2^{\perp},1) ]_{ii} } \, \sum_{\bq_2} \bM(2,1)_{\bq_2 \bq_1} \te^{- 2 \pi \ti \sum_j \om(2^{\perp},1)_{ji} \frac{q_{2j}}{k} } ( s_{2 \bq'_2} \otimes \de_{2 \bq'_2}),  \\
& \qquad \qquad \qquad \qquad \text{ where} \quad q'_{2j} = q_{2j}+ \om(2,1)_{ji} \\
&= \la \te^{ \pi \ti k [\,{ ^t \om(2,1)} \om(2^{\perp},1) ]_{ii} } \, \sum_{\bq''_2} \bM(2,1)_{\bq''_2 \bq_1} \te^{- 2 \pi \ti \sum_j \om(2^{\perp},1)_{ji} \frac{q_{2j}}{k} } ( s_{2 \bq_2} \otimes \de_{2 \bq_2} ), \\
& \qquad \qquad \qquad \qquad \text{ where} \quad q''_{2j} = q_{2j}- \om(2,1)_{ji} \\
&= \la \te^{ 2 \pi \ti \frac{q_{1i}}{k} } \, \sum_{\bq_2} \bM(2,1)_{\bq_2 \bq_1} \, ( s_{2 \bq_2} \otimes \de_{2 \bq_2} ) \, = \,
F_{\CP_2 \CP_1} [(\la, w_{1i}) \cdot ( s_{1 \bq_1} \otimes \de_{1 \bq_1} ) ] 
\end{split}
\end{equation*}
and
\begin{equation*}
\begin{split}
& (\la, w_{1i}^{\perp}) \cdot F_{\CP_2 \CP_1} ( s_{1 \bq_1} \otimes \de_{1 \bq_1} ) \, = \,
\sum_{\bq_2} \bM(2,1)_{\bq_2 \bq_1} (\la, w_{1i}^{\perp}) \cdot ( s_{2 \bq_2} \otimes \de_{2 \bq_2} ) \\
&= \la \te^{ \pi \ti k [\,{ ^t \om(2,1^{\perp})} \om(2^{\perp},1^{\perp}) ]_{ii} } \, \sum_{\bq_2} \bM(2,1)_{\bq_2 \bq_1} \te^{- 2 \pi \ti \sum_j \om(2^{\perp},1^{\perp})_{ji} \frac{q_{2j}}{k} } ( s_{2 \bq'_2} \otimes \de_{2 \bq'_2} ), \\
&\qquad \qquad \qquad \qquad \text{ where} \quad q'_{2j} = q_{2j}+ \om(2,1^{\perp})_{ji} \\
&= \la \te^{ \pi \ti k [\,{ ^t \om(2,1^{\perp})} \om(2^{\perp},1^{\perp}) ]_{ii} } \, \sum_{\bq_2} \bM(2,1)_{\bq''_2 \bq_1} \te^{- 2 \pi \ti \sum_j \om(2^{\perp},1^{\perp})_{ji} \frac{q_{2j}}{k} } ( s_{2 \bq_2} \otimes \de_{2 \bq_2} ), \\
&\qquad \qquad \qquad \qquad \text{ where} \quad q''_{2j} = q_{2j}- \om(2,1^{\perp})_{ji} \\
&= \la \sum_{\bq_2} \, \te^{\frac{2 \pi \ti}{k} [\om(2,1)^{-1} \om(2,1^{\perp}) \bq_{1}]_{i} + \frac{\pi \ti}{k} [\om(2,1)^{-1} \om(2,1^{\perp})]_{ii} - \frac{2 \pi \ti}{k} [\om(2,1)^{-1}(\bq_2 + k \bl)]_{i} } \,  \times \\
& \qquad \qquad \qquad \qquad \times \bM(2,1)_{\bq_2 \bq_1} \, ( s_{2 \bq_2} \otimes \de_{2 \bq_2} ) \\
&= F_{\CP_2 \CP_1}[ \la \, ( s_{1 {\bq_1}_{(i+1)}} \otimes \de_{1 {\bq_1}_{(i+1)}} ) ] \, = \, F_{\CP_2 \CP_1} [(\la, w_{1i}^{\perp}) \cdot ( s_{1 \bq_1} \otimes \de_{1 \bq_1} ) ] \, .
\end{split}
\end{equation*}
A similar proof goes for $\CP_1$ and $\CP_2$ nontransverse polarizations and we omit it.
\end{proof}

In view of the Stone-von Neumann theorem, the fact that the unitary operators between the Hilbert spaces of different polarizations intertwine the representations of the Heisenberg group $H_k$ on these spaces, leads naturally to the conclusion that these operators compose transitively up to a phase factor. 
In the next section we will prove this fact through a direct computation which will give us at the same time the explicit expression of the projective factor involved.

\bsk


\section{The projective Hilbert space} \label{s:prHsp}

We begin by recalling the definition of the Maslov-Kashiwara index (or signature index) and some of its properties \cite{LV}.
Let $L_{1}, L_{2}, L_{3}$ be three Lagrangian subspaces of the 
$2g$-dimensional symplectic vector space $(\CV, \om)$. Then
\begin{defi}
The Maslov-Kashiwara index $\tau( L_{1}, L_{2}, L_{3})$ is defined as the signature of the quadratic form $G$ on the $3g$-dimensional vector space $L_{1} \oplus L_{2} \oplus L_{3}$, given by
\[ G(X_1 \oplus X_2 \oplus X_3) \, = \, \om(X_1, X_2) + \om(X_2, X_3) + \om(X_3, X_1) \]
\end{defi}
Let us  assume that $L_1$ and $L_3$ are transverse, i.e. $L_1 \cap L_3 = 0$. 
Then $\CV = L_1 \oplus L_3$ and we let $p_{31} : \CV \ra L_3$ denote the projection onto $L_3$ along $L_1$. 
The Lagrangian subspaces $L_1$ and $L_3$ determine a symmetric bilinear form on $L_2$ given by
\begin{equation} \label{e:Hmat}
H(X_{2}, X'_{2}) \, = \, \om( X_{2} , p_{31} X'_{2}) , \quad \text{ for any } X_{2}, X'_{2} \in L_{2} \, .
\end{equation}
Then we have the lemma (\cite{LV}, \S 1.5):

\begin{lem} \label{l:Milemma}
If $L_1$ and $L_3$ are transverse, then $\tau( L_1, L_2, L_3)$ is equal to the signature of the quadratic form $H(X_2, X_2)$ on $L_2$:
\[ \tau( L_1, L_2, L_3) \, = \, \mathrm{sgn } \, H \]
\end{lem}

The properties of the Maslov-Kashiwara index $\tau$ are summarized in the following proposition \cite{LV,Go}:
\begin{prop} \label{p:propMKind}
The function $\tau : \lag \times \lag \times \lag \ra \BZ$ has the following properties: \\
(1) $\tau(a L_1, a L_2, a L_3) = \tau(L_1, L_2, L_3)$, for any $a \in Sp(\CV)$. \\
(2) $\tau(L_1, L_2, L_3)$ is unchanged (changes sign) under an even (odd) permutation of the triple $(L_1, L_2, L_3)$. \\
(3) $\tau(L_1, L_2, L_3) - \tau(L_1, L_2, L_4) + \tau(L_1, L_3, L_4) - \tau(L_2, L_3, L_4) =0$, for any \\
$L_1, L_2, L_3, L_4 \in \lag$. \\
(4) For any triple $L_1, L_2, L_3$ of Lagrangian subspaces of $\CV \: (\dim \CV = 2g)$ we have
$ \tau(L_1, L_2, L_3) \equiv g + \dim(L_1 \cap L_2) + \dim(L_2 \cap L_3) + \dim(L_3 \cap L_1) \, \pmod{2} $. 
\end{prop}
 
Now let us consider three real invariant polarizations $\CP_1, \CP_2, \CP_3$ of the symplectic $2g$-torus $(\CT, \om)$. 
As before  we let $L_{\CP_1}, L_{\CP_2}, L_{\CP_3}$ denote the corresponding Lagrangian subspaces of $(\CV, \om)$.
Quantization associates to each polarization $\CP_i$ the Hilbert space $\CH_{\CP_i}$. 
These Hilbert spaces are related by the unitary intertwining operators $F_{\CP_1 \CP_3}$, $ \, F_{\CP_3 \CP_2}$,
$ \, F_{\CP_2 \CP_1}$ introduced in Sect.\ref{s:BKS}.
We will prove that the operator $F_{\CP_1 \CP_3} \circ F_{\CP_3 \CP_2} \circ F_{\CP_2 \CP_1}$ on $\CH_{\CP_1}$ is proportional to  the identity operator. 
More precisely:
\begin{thm} \label{t:trans-proj}
For any three invariant real polarizations $\CP_1, \CP_2, \CP_3$ of the symplectic torus $(\CT, k \om)$ we have
\begin{equation*}
 F_{\CP_1 \CP_3} \circ F_{\CP_3 \CP_2} \circ F_{\CP_2 \CP_1}  \, = \, \te^{- \frac{\pi \ti}{4} \tau(L_{\CP_1}, L_{\CP_2}, L_{\CP_3}) } \;  I .
\end{equation*}
\end{thm}
\nin For the proof of the above theorem we will need the following proposition.
\begin{prop}[Reciprocity formula for Gauss sums] \label{p:Gauss}
Let $Q$ be a nonsingular $(g \times g)$-matrix with integer coefficients, $a \in 2 \BZ$ a positive even integer and $\bw = (w_{1}, \dots, w_{g}) \in \BQ^g$ such that $a \bw \in \BZ^g$. 
Then we have:
\[ \sum_{ \bq \in (\BZ/ a \BZ)^{g} } \, \te^{ \frac{\pi \ti}{a} \,{ ^t \bq} Q \bq + 2 \pi \ti \,{ ^t \bw} \cdot \bq }
\, = \, \left| \frac{a^{g}}{\det Q} \right|^{\frac{1}{2}} \: \te^{ \frac{\pi \ti}{4} \mathrm{sgn} Q} \; \sum_{ \bm \in \BZ^{g} / Q \BZ^{g} } \, \te^{ - \pi \ti a \,{ ^t (\bm + \bw)} Q^{-1} (\bm + \bw) } \]
\end{prop}
\nin The proof of the above relation is a straightforward generalization of the argument used in (\cite{La}, ch.4) to prove the one-dimensional Gauss sum formula. We prove now Theorem (\ref{t:trans-proj}).
\begin{proof}
{\em Case 1}. Let us consider first the case of three mutually transverse polarizations $\CP_1, \CP_2, \CP_3$ of $(\CT, k\om)$. 
As in (\ref{e:Zbasis}), 
we choose integer symplectic bases $( W_{1i}\, ;W_{1i}^{\perp} )$, $( W_{2i}\, ;W_{2i}^{\perp} )$, $( W_{3i}\, ;W_{3i}^{\perp} )$ of the symplectic vector space $(\CV, \om)$, with
$L_{\CP_1} \, = \, \text{span}_{\BR} \{ W_{1i} \}_{i=1, \dots ,g} , \,
L_{\CP_2} \, = \, \text{span}_{\BR} \{ W_{2i} \}_{i=1, \dots ,g}$, 
$L_{\CP_3} \, = \, \text{span}_{\BR} \{ W_{3i} \}_{i=1, \dots ,g}$.
With respect to the above bases we introduce as in (\ref{e:ommat}) the matrices $\om(2,1)$, $\om(2^{\perp},1)$,$ \om(3,2)$, $\om(3,1)$, etc. 
Since $L_{\CP_1}, L_{\CP_2}, L_{\CP_3}$ are mutually transverse, the matrices $\om(2,1)$, $\om(3,2)$ and $\om(3,1)$ are nonsingular.
The operators $F_{\CP_1 \CP_3}, F_{\CP_3 \CP_2} $ and $ F_{\CP_2 \CP_1}$ are represented by the unitary matrices $\bM(1,3), \bM(3,2)$ and $ \bM(2,1)$, as expressed by (\ref{e:trmat}). 
Then we have to compute the sum
\begin{gather} \label{e:trsum}
\sum_{\bq_2, \bq_3 \in (\BZ / k \BZ)^g } \: \bM(1,3)_{\bq'_1 \bq_3} \, \bM(3,2)_{\bq_3 \bq_2} \, \bM(2,1)_{\bq_2 \bq_1} \\
= \, \sum_{\bq_2, \bq_3 \in (\BZ / k \BZ)^g } \: \overline{\bM(3,1)}_{ \bq_3 \bq'_1 } \, \bM(3,2)_{\bq_3 \bq_2} \, \bM(2,1)_{\bq_2 \bq_1} \quad \notag \\
= \, \left| k^{3g} \, \det \om(2,1) \, \det \om(3,2) \, \det \om(3,1) \right|^{- \frac{1}{2}} \times \quad \notag \\
\begin{split} \notag
\times \sum_{\bq_2, \bq_3 \in (\BZ / k \BZ)^g } \; 
& [ \; \sum_{[\bl_{31} ] \in \BZ^g / \om(3,1) \BZ^g} \, \te^{- \frac{\pi \ti}{k} \bA(3,1)_{ \bl_{31} \bq_3 \bq'_1} } \; ] \\
& [ \; \sum_{[\bl_{32} ] \in \BZ^g / \om(3,2) \BZ^g} \, \te^{\frac{\pi \ti}{k} \bA(3,2)_{ \bl_{32} \bq_3 \bq_2} } \; ] \\
& [ \; \sum_{[\bl_{21} ] \in \BZ^g / \om(2,1) \BZ^g} \, \te^{\frac{\pi \ti}{k} \bA(2,1)_{ \bl_{21} \bq_2 \bq_1} } \; ]  \; . \\
\end{split}
\end{gather}
Isolating the terms with the $\bq_2$ variable, we get for the corresponding sum the expression
\begin{equation*}
\CS_{\bq_2} \, = \, \sum_{\bq_2 \in (\BZ / k \BZ)^g } \: 
\te^{ \frac{\pi \ti}{a} \,{ ^t \bq_2} Q \bq_2 + 2 \pi \ti \,{ ^t \bw} \cdot \bq_2 }
\end{equation*}
with
\begin{align*}
\bw \, &= \, \,{ ^t \om(2,1)^{-1} } \frac{\bq_{1}}{k} - \om(3,2)^{-1} (\frac{\bq_{3}}{k} + \bl_{32}) - \om(2^{\perp},1) \om(2,1)^{-1} \bl_{21} \\
a \, &= \, k \left| \det \om(2,1) \, \det \om(3,2) \right| \, \in 2 \BZ
\end{align*}
and $Q$ the $(g \times g)$-matrix with $\BZ$-coefficients
\begin{equation*} 
Q \, = \, \left| \det \om(2,1) \, \det \om(3,2) \right| \: [ \om(2^{\perp},3) \om(2,3)^{-1} - \om(2^{\perp}, 1) \om(2,1)^{-1} ] \, .
\end{equation*}
From Lemma (\ref{l:symm}) it follows that $Q$ is a symmetric matrix. 
A simple computation shows that $Q$ can be expressed as
\begin{equation} \label{e:Qmat} 
Q \, = \, \left| \det \om(2,1) \, \det \om(3,2) \right| \: \om(3,2)^{-1} \om(3,1) \om(2,1)^{-1}
\end{equation}
Thus $\det Q \neq 0$.
It is not difficult to see that we can extend the summation range of the $\bq_{2}$-variable in (\ref{e:trsum}), so that we can write
\begin{equation*} 
\CS_{\bq_{2}} \, = \, \frac{1}{(a/k)^{g}} \: 
\sum_{\bq_{2} \in (\BZ / a \BZ)^{g} } \: 
\te^{ \frac{\pi \ti}{a} \,{ ^t \bq_2} Q \bq_{2} + 2 \pi \ti \,{ ^t \bw} \cdot \bq_{2} }
\end{equation*}
By applying the Gauss sum reciprocity formula (\ref{p:Gauss}) this becomes
\begin{align*}
\begin{split}
\CS_{\bq_{2}} \, &= \, \frac{1}{(a/k)^{g}} \, \left| \frac{a^{g}}{\det Q} \right|^{\frac{1}{2}} \: \te^{ \frac{\pi \ti}{4} \text{sgn} Q} \; \sum_{ \bl \in \BZ^{g} / Q \BZ^{g} } \, \te^{ - \pi \ti a \,{ ^t (\bl + \bw)} Q^{-1} (\bl + \bw) } \\
&= \, \frac{1}{\left| \det \om(2,1) \, \det \om(3,2) \right|^{g}} \, \left| \frac{ k^{g} \det \om(2,1) \, \det \om(3,2)}{\det \om(3,1)} \right|^{\frac{1}{2}} \, 
\te^{ \frac{\pi \ti}{4} \text{sgn} Q} \times \\
& \times \sum_{ \bl \in \BZ^{g} / Q \BZ^{g} } \, \te^{ - \pi \ti k \,{ ^t (\bl + \bw)} \om(2,1) \om(3,1)^{-1} \om(3,2) (\bl + \bw)}  \; .
\end{split}
\end{align*}
We introduce this back into (\ref{e:trsum}) and proceed to compute the sum over the $\bq_{3}$-variable. 
Retaining only the terms which contain $\bq_{3}$, we find that this sum can be expressed as
\begin{equation*} 
\CS_{\bq_{3}}  = 
\sum_{\bq_{3} \in (\BZ / k \BZ)^{g} } \, 
\te^{ 2 \pi \ti \,{ ^t \bq_{3}} \,{ ^t \om(3,1)^{-1} }  [ (\frac{\bq'_{1}}{k} - \frac{\bq_{1}}{k}) + \,{ ^t \om(3^{\perp},1)} (\bl_{31} - \bl_{32} - \om(3,2^{\perp}) \bl_{21} + \om(3,2) \bl) ] } \; .
\end{equation*}
Letting $b= k | \det \om(3,1) |$, an inspection of (\ref{e:trsum}) shows that we can write this sum as
\begin{equation*}
\begin{split}
\CS_{\bq_{3}}  &= \frac{1}{(b/k)^{g}} 
\sum_{\bq_{3} \in (\BZ / b \BZ)^{g} } \, 
\te^{ 2 \pi \ti \,{ ^t \bq_{3}} \,{ ^t \om(3,1)^{-1} } [ (\frac{\bq'_{1}}{k} - \frac{\bq_{1}}{k}) + \,{ ^t \om(3^{\perp},1)} (\bl_{31} - \bl_{32} - \om(3,2^{\perp}) \bl_{21} + \om(3,2) \bl) ] } \\
&= \, k^{g} \, \de_{\bq'_{1} \bq_{1}} \: \de_{\bl_{31}, \bl_{32} + \om(3,2^{\perp}) \bl_{21} - \om(3,2) \bl}
\end{split}
\end{equation*}
Finally, plugging the above result back into (\ref{e:trsum}) and performing the remaining sums over $\bl_{21}, \bl_{32}, \bl_{31}$ and $\bl$, we obtain
\begin{equation*}
\sum_{\bq_{2}, \bq_{3} \in (\BZ / k \BZ)^{g} } \: \bM(1,3)_{\bq'_{1} \bq_{3}} \, \bM(3,2)_{\bq_{3} \bq_{2}} \, \bM(2,1)_{\bq_{2} \bq_{1}}  \, = \, \te^{ \frac{\pi \ti}{4} \text{sgn} Q} \: \de_{\bq'_{1} \bq_{1}} 
\end{equation*}

\nin{\em Case 2}. We consider now the case of three polarizations $\CP_1, \CP_2, \CP_3$ such that $\CP_3$ is transverse to $\CP_1$ and $\CP_2$, but $\CP_1$ and $\CP_2$ are not transverse to each other.
We choose for $(\CV, \om)$ integer symplectic bases $( W_{1i}\, ;W_{1i}^{\perp} )$ and $( W_{2i}\, ;W_{2i}^{\perp} )$, such that $L_{\CP_1} \, = \, \text{span}_{\BR} \{ W_{1i} \}_{i=1, \dots ,g}$, $\: L_{\CP_2} \, = \, \text{span}_{\BR} \{ W_{2i} \}_{i=1, \dots ,g}$ and satisfying the conditions (\ref{e:intbas})-(\ref{e:intbases}), and $( W_{3i}\, ;W_{3i}^{\perp} )$ with $L_{\CP_3} \, = \, \text{span}_{\BR} \{ W_{3i} \}_{i=1, \dots ,g}$.
The matrices $\om(3,1)$ and $\om(3,2)$ defined with respect to these bases are nonsingular, while $\om(2,1)$ has the form shown in (\ref{e:ntrom}) and only the reduced matrix $\check{\om}(2,1)$ is nonsingular.
Using the expressions (\ref{e:trmat}) and (\ref{e:ntrmatrix}), we proceed to compute the sum
\begin{gather} \label{e:ntrsum}
 \sum_{\bq_2, \bq_3 \in (\BZ / k \BZ)^g } \: \bM(1,3)_{\bq'_1 \bq_3} \, \bM(3,2)_{\bq_3 \bq_2} \, \bM(2,1)_{\bq_2 \bq_1} \\
= \, \left| k^{2g} k^{h} \, \det \check{\om}(2,1) \, \det \om(3,2) \, \det \om(3,1) \right|^{- \frac{1}{2}} \: \prod_{i=h+1}^{g} \de_{q_{2i} q_{1i}} \times \notag \\ 
\begin{split} \notag
\times \sum_{\bq_2, \bq_3 \in (\BZ / k \BZ)^g } \; 
& [ \; \sum_{[\bl_{31} ] \in \BZ^g / \om(3,1) \BZ^g} \, \te^{- \frac{\pi \ti}{k} \bA(3,1)_{ \bl_{31} \bq_3 \bq'_1} } \; ] \\
& [ \; \sum_{[\bl_{32} ] \in \BZ^g / \om(3,2) \BZ^g} \, \te^{\frac{\pi \ti}{k} \bA(3,2)_{ \bl_{32} \bq_3 \bq_2} } \; ] \\
& [ \; \sum_{[\bl_{21} ] \in \BZ^h / \check{\om}(2,1) \BZ^h} \, \te^{\frac{\pi \ti}{k} \check{\bA}(2,1)_{ \bl_{21} \check{\bq}_2 \check{\bq}_1} } \; ] \\
\end{split}
\end{gather}
Let us deal first with the sum over $\bq_2$ and make use of the fact that $q_{2i} = q_{1i} $ for $ i=h+1, \dots ,g$. 
As before we use the notation \text{$\check{\bq}_2 = (q_{21}, \dots, q_{2h})$}. We find that the sum over the terms which contain $\check{\bq}_2$ reads:
\begin{equation*}
\CS_{\check{\bq}_{2}} \, = \, \sum_{\check{\bq}_{2} \in (\BZ / k \BZ)^{h} } \: 
\te^{ \frac{\pi \ti}{k} \,{ ^t \check{\bq}_2} B \check{\om}(2,1)^{-1} \check{\bq}_{2} + 2 \pi \ti \,{ ^t \bw} \cdot \check{\bq}_{2} } \;  ,
\end{equation*}
where we introduced the notations
\begin{align*}
B_{ij} \, &= \, ( \om(3,2)^{-1} \om(3,1) )_{ij} , \qquad i,j=1, \dots, h  \\
w_{i} \, &= \, \sum_{j=h+1}^{g} \, (\om(2^{\perp},3) \om(2,3)^{-1})_{ij} \frac{q_{1j}}{k} - \sum_{j=1}^{h} \, \check{\om}(2,1)^{-1}_{ji} \frac{q_{1j}}{k}  \\
 & - \sum_{j=1}^{g} \om(3,2)^{-1}_{ij} (\frac{\bq_{3}}{k} + \bl_{32})_{j} - \sum_{j=1}^{h} \, ( \check{\om}(2^{\perp},1) \check{\om}(2,1)^{-1})_{ij} (\bl_{21})_{j} 
\end{align*}
Now let us define
\begin{align*}
a \, &= \, k \left| \det \check{\om}(2,1) \, \det \om(3,2) \right|  \, \in 2 \BZ \\
\check{Q} \, &= \, \left| \det \check{\om}(2,1) \, \det \om(3,2) \right| \: B \check{\om}(2,1)^{-1} 
\end{align*}
The $(h \times h)$-matrix $\check{Q}$ is symmetric and has integer coefficients. 
We also have 
$ \om(3,2)^{-1} \om(3,1) = \left( \begin{smallmatrix} B & 0 \\ 0 & I \end{smallmatrix} \right)$
and thus $(B^{-1})_{ij} = (\om(3,1)^{-1} \om(3,2))_{ij}, \, i,j=1, \dots , h$. 
Then, since (\ref{e:ntrsum}) is unchanged if we appropriately extend the summation range of the $\check{\bq}_{2}$-variable, we can write
\begin{align*}
\begin{split}
\CS_{\check{\bq}_{2}} \, &= \, \frac{1}{(a/k)^{h}} \: 
\sum_{\check{\bq}_{2} \in (\BZ / a \BZ)^{h} } \: 
\te^{ \frac{\pi \ti}{a} \,{ ^t \check{\bq}_{2}} \check{Q} \check{\bq}_{2} + 2 \pi \ti \,{ ^t \bw} \cdot \check{\bq}_{2} } \\
&= \, \frac{1}{\left| \det \check{\om}(2,1) \, \det \om(3,2) \right|^{h}} \, \left| \frac{ k^{h} \det \check{\om}(2,1) \, \det \om(3,2)}{\det \om(3,1)} \right|^{\frac{1}{2}} \, 
\te^{ \frac{\pi \ti}{4} \text{sgn} \check{Q}} \times \\
& \times \sum_{ \bl \in \BZ^{h} / \check{Q} \BZ^{h} } \, \te^{ - \pi \ti k \,{ ^t (\bl + \bw)} \check{\om}(2,1) B^{-1} (\bl + \bw)}  \,  .
\end{split}
\end{align*}
For the last equality in the above expression we used the Gauss sum formula of Proposition (\ref{p:Gauss}).
Now, introducing the notation $\bl'_{21} = (\bl_{21},0, \dots ,0)$ and $\bl' = (\bl , 0, \dots ,0)$ for the extensions of $\bl_{21}$ and $\bl$ to $g$-dimensional vectors, we find that the resulting sum over the terms containing $\bq_{3}$ is
\begin{equation*}
\begin{split} 
\CS_{\bq_{3}}  &= 
\sum_{\bq_{3} \in (\BZ / k \BZ)^{g} } \, 
\te^{ 2 \pi \ti \,{ ^t \bq_{3}} \,{ ^t \om(3,1)^{-1} }  [ (\frac{\bq'_{1}}{k} - \frac{\bq_{1}}{k}) + \,{ ^t \om(3^{\perp},1)} (\bl_{31} - \bl_{32} - \om(3,2^{\perp}) \bl'_{21} + \om(3,2) \bl') ] } \\
&= \, k^{g} \de_{\bq'_{1} \bq_{1}} \: \de_{\bl_{31}, \bl_{32} + \om(3,2^{\perp}) \bl'_{21} - \om(3,2) \bl' } \, .
\end{split}
\end{equation*}
Finally, inserting the above result back into (\ref{e:ntrsum}), we obtain
\begin{equation*}
\sum_{\bq_{2}, \bq_{3} \in (\BZ / k \BZ)^{g} } \: \bM(1,3)_{\bq'_{1} \bq_{3}} \, \bM(3,2)_{\bq_{3} \bq_{2}} \, \bM(2,1)_{\bq_{2} \bq_{1}}  \, = \, \te^{ \frac{\pi \ti}{4} \text{sgn} \check{Q}} \: \de_{\bq'_{1} \bq_{1}} 
\end{equation*}

Now, let us examine more closely the symmetric matrices $Q$ and $\check{Q}$ and their relation to the symmetric bilinear form $H$ from (\ref{e:Hmat}) introduced in relation to the Maslov-Kashiwara index.
Since $L_{\CP_1} \cap L_{\CP_3} = 0 $, we can take $( W_{3i}\, ,W_{1i})$ as basis for $\CV$. With respect to this basis $W_{2i}$ has the decomposition
\begin{align*}
W_{2i} \, &= \, (\om(2,1) \om(3,1)^{-1})_{ij} W_{3j} + ( \om(2,3) \om(1,3)^{-1})_{ij} W_{1j}  \,  .
\end{align*}
Using this we find for the symmetric bilinear form $H$ on $L_{\CP_2}$ the expression
\begin{align*}
H_{ij} \, &= \, H(W_{2i}, W_{2j}) \, = \, \om(W_{2i},  p_{31} \, W_{2j}) \, \\ 
&= \, - [\om(2,1) \om(3,1)^{-1} \om(3,2) ]_{ij} \; , \qquad \qquad i,j=1, \dots ,g \, . 
\end{align*}
If $L_{\CP_1}$ and $ L_{\CP_2}$ are transverse, then comparing the above formula to (\ref{e:Qmat}) we find that
\begin{equation*}
H \, = \, - \left| \det \om(2,1) \det \om(3,2) \right| \: Q^{-1} \, .
\end{equation*}
Thus $\text {sgn } H = - \text{ sgn } Q$ in this case.
If $L_{\CP_1}$ and $ L_{\CP_2}$ are not transverse, then we have
\begin{equation*}
\begin{split}
H \, &= \, - \om(2,1) \om(3,1)^{-1} \om(3,2) \, = \, -
\begin{pmatrix} \check{\om}(2,1) & 0 \\ 0 & 0 \end{pmatrix} 
\begin{pmatrix} B^{-1} & 0 \\ 0 & \text{I} \end{pmatrix} \\
&= \, - \begin{pmatrix} \check{\om}(2,1) B^{-1} & 0 \\ 0 & 0 \end{pmatrix}  \, = \, - \left| \det \check{\om}(2,1) \det \om(3,2) \right| \:  \begin{pmatrix} \check{Q}^{-1} & 0 \\ 0 & 0 \end{pmatrix} 
\end{split}
\end{equation*}
and  $\text {sgn } H = - \text{ sgn } \check{Q}$.
Therefore, having in view Lemma (\ref{l:Milemma}), we obtain in either case that
\begin{equation} \label{e:finsum}
\sum_{\bq_{2}, \bq_{3} \in (\BZ / k \BZ)^{g} } \: \bM(1,3)_{\bq'_{1} \bq_{3}} \, \bM(3,2)_{\bq_{3} \bq_{2}} \, \bM(2,1)_{\bq_{2} \bq_{1}}  \, = \, \te^{ - \frac{\pi \ti}{4} \tau(L_{\CP_{1}}, L_{\CP_{2}}, L_{\CP_{3}})} \: \de_{\bq'_{1} \bq_{1}}  \; .
\end{equation}
Given three arbitrary polarizations $\CP_{1}, \CP_{2}, \CP_{3}$ of $\CT$ we can always find a polarization $\CP_{4}$ transversal to all of them. 
Together with the cocycle property (\ref{p:propMKind})-(3)) of $\tau$ this implies that (\ref{e:finsum}) holds in general.
\end{proof}
Thus geometric quantization assigns canonically to the symplectic torus $(\CT, k \om)$ a projective Hilbert space $\BP \CH$.
A similar property to the one expressed by Theorem (\ref{t:trans-proj}) for the quantization of symplectic tori holds for the quantization of symplectic vector spaces using constant real polarizations. 
In the vector space case too, the unitary operators relating the Hilbert spaces of different polarizations satisfy a transitive composition law up to a projective factor given in terms of the Maslov-Kashiwara index, exactly as in (\ref{t:trans-proj}).
For the treatment of the quantization of a symplectic vector space we refer the reader to \cite{GS, LV, Wo}.

\bsk


\section{The metaplectic correction} \label{s:metapl}

The quantization of the $2g$-dimensional torus $\CT$ described in Sect.\ref{s:quanttg} associated to the data $(\CT, k \om, \CP)$, with $\CP$ an invariant real polarization, a Hilbert space $\CHP$.
In constructing $\CHP$ we made use of the bundle of half-densities on $\CP$ and we proved in Sect.\ref{s:BKS} and Sect.\ref{s:prHsp} that the Hilbert spaces of different polarizations are projectively identified.
In this section we show that by using half-forms instead of half-densities we can associate to $(\CT, k \om)$ a true Hilbert space $\tilde{\CH}$, not merely a projective one $\BP \CH$.

\msk

\subsection{Metaplectic and metalinear structures} \label{ss:met}

A consistent choice of the bundle of half-forms on $\CP$, for each polarization $\CP$, requires the choice of a metaplectic frame bundle on $(\CT, k \om)$.

Let $SP(\CT)$ denote the symplectic frame bundle of $(\CT, k \om)$. The bundle $SP(\CT)$ is the principal $Sp(2g,\BR)$-bundle over $\CT$ whose fibre at $x$ is the set of all ordered bases 
$(U_1, \dots ,U_g ; V_1, \dots ,V_g)$ of $T_x(\CT)$ which satisfy $k \om(U_i,V_j) = \de_{ij}$ and $k \om(U_i,U_j) = k \om(V_i,V_j) = 0$.
Let $Mp(2g,\BR)$ denote the metaplectic group, that is, the double cover of the symplectic group $Sp(2g,\BR)$, and let $\rho : Mp(2g,\BR) \ra Sp(2g,\BR)$ be the covering homomorphism.
A {\em metaplectic frame bundle} on $(\CT, k \om)$ is a principal $Mp(2g,\BR)$-bundle $MP(\CT)$ over $\CT$ together with a $2:1$ covering map $p : MP(\CT) \ra SP(\CT)$, such that the following diagram commutes:
\begin{equation*}
\begin{CD}
MP(\CT) \times Mp(2g,\BR) @>>> MP(\CT) \\
@V{p \times \rho}VV    @VV{p}V \\
SP(\CT) \times Sp(2g,\BR) @>>> SP(\CT)
\end{CD}
\end{equation*}
The horizontal arrows denote the right group actions.
The set of isomorphism classes of metaplectic frame bundles on $\CT$ is a principal homogeneous space for the cohomology group $H^1(\CT;\BZ_2) \;$ (\cite{GS}, ch.V).

Let us consider now the complex linear group $GL(g,\BC)$ and its double cover, the complex metalinear group $ML(g,\BC)$.
Let $\rho : ML(g,\BC) \ra GL(g,\BC)$ be the covering projection and $\chi : ML(g,\BC) \ra \BC$ the character of $ML(g,\BC)$ defined as the holomorphic square root of the complex character $\, \det \circ \rho$ such that $\chi(\tilde{e})=1$, where $\tilde{e}$ is the identity element of $ML(g,\BC)$.
The real linear group $GL(g,\BR)$ is a subgroup of $GL(g,\BC)$ and its preimage $ML(g,\BR) = \rho^{-1} GL(g,\BR) \subset ML(g,\BC)$ is the real metalinear group.
The linear group $GL(g,\BR)$ is naturally embedded in $Sp(2g,\BR)$ and the metalinear group $ML(g,\BR)$ can be identified with a subgroup of the metaplectic group $Mp(2g,\BR)$ \cite{GS,Bl}.
The element $\tilde{\vep}$ in $Mp(2g,\BR)$, $\tilde{\vep} \neq \tilde{e}$, and which covers the identity $e$ in $Sp(2g,\BR)$,  belongs to $Mp(2g,\BR) \cap ML(g,\BR)$ and $\chi(\tilde{\vep}) = \te^{\pi \ti} = -1$.

We outline the standard construction \cite{Bl,GS} which, given a metaplectic frame bundle $MP(\CT)$, defines the bundle of metalinear Lagrangian frames and, for each polarization $\CP$, the bundle of half-forms on $\CP$.
Let $N$ denote the subgroup of $Sp(2g,\BR)$ given by
$N \, = \, \{ \left( \begin{smallmatrix} I & S \\ 0 & I \end{smallmatrix} \right) \mid \, S \in GL(g,\BR), \, ^t S = S \}$.
The group $N$ is simply connected and its preimage $\rho^{-1} N$ under the projection $\rho : Mp(2g,\BR) \ra Sp(2g,\BR)$ has two components, each diffeomorphic to $N$.
We identify the identity component with $N$ and regard $N$ as a subgroup of $Mp(2g,\BR)$.

The bundle $SP(\CT)/N$ on $\CT$ obtained by taking the quotient of the symplectic frame bundle under the right-action of the subgroup $N \subset Sp(2g,\BR)$ is naturally identified with the bundle $\CF l(\CT)$ of Lagrangian frames of $T(\CT)$.
We recall that a Lagrangian frame $\underline{V}$ of $T_x(\CT)$ is an ordered $g$-tuple $\underline{V}=(V_1, \dots ,V_g)$ of vectors $V_i \in T_x(\CT)$ such that $k \om(V_i,V_j) = 0$. 
The group $GL(g,\BR)$ acts freely on $\CF l_x(\CT)$ on the right.
There is a surjective map $\la : SP(\CT) \ra \CF l(\CT)$ which sends a symplectic frame $(\underline{V};\underline{W}) = (V_1, \dots ,V_g ;W_1, \dots ,W_g)$ of $T_x(\CT)$ to the Lagrangian frame $\underline{V} = (V_1, \dots ,V_g)$.
The map $\la : SP(\CT) \ra \CF l(\CT)$ is a principal $N$-bundle.
The group $GL(g,\BR) \subset Sp(2g,\BR)$ acts on $SP(\CT)$ on the right and $\la$ is equivariant with respect to the $GL(g,\BR)$-actions on $SP(\CT)$ and $\CF l(\CT)$.
Since $GL(g,\BR)$ normalizes $N$ in $Sp(2g,\BR)$, it also acts  on $SP(\CT)/N$ on the right and we have a $GL(g,\BR)$-equivariant isomorphism of $SP(\CT)/N$ with $\CF l(\CT)$. 

Consider now the bundle $\tilde{\CF l}(\CT) = MP(\CT)/N$ on $\CT$. 
The metalinear group $ML(g,\BR)$ normalizes $N$ in $Mp(2g,\BR)$.
Therefore, the quotient map $\tilde{\la} : MP(\CT) \ra\tilde{\CF l}(\CT)$ is $ML(g,\BR)$-equivariant and we have the commutative diagram:
\begin{equation} \label{e:cdmetLag}
\begin{CD}
MP(\CT)  @>{\tilde{\la}}>> \tilde{\CF l}(\CT) \\
@VpVV     @VVpV \\
SP(\CT)  @>{\la}>> \CF l(\CT)
\end{CD}
\end{equation}
$\tilde{\CF l}(\CT) $ is a double cover of $\CF l(\CT)$ called the bundle of {\em metalinear Lagrangian frames} on $\CT$.
One can easily see that $ML(g,\BR)$-orbits in $\tilde{\CF l}(\CT)$ project to $GL(g,\BR)$-orbits in $\CF l(\CT)$.

Now, for any real polarization $\CP$ of $(\CT, k \om)$, let $GL(\CP)$ denote the bundle of linear frames of $\CP$. It is the principal $GL(g,\BR)$-bundle on $\CT$ whose fibre over $x$ is the set of all ordered bases $\underline{W} = (W_1, \dots ,W_g)$ of $\CP_x \subset T_x(\CT)$.
A {\em metalinear frame bundle} for $\CP$ is a principal $ML(g,\BR)$-bundle $ML(\CP)$ on $\CT$ together with a $2:1$ covering map $p : ML(\CP) \ra GL(\CP)$ such that the following diagram commutes:
\begin{equation*}
\begin{CD}
ML(\CP) \times ML(g,\BR) @>>> ML(\CP) \\
@V{p \times \rho}VV    @VV{p}V \\
GL(\CP) \times GL(g,\BR) @>>> GL(\CP)
\end{CD}
\end{equation*}
Since $GL(\CP) \subset \CF l(\CT)$, we define the  bundle of metalinear frames of $\CP$ as the preimage $ML(\CP) = p^{-1}(GL(\CP))$ under the covering projection $p: \tilde{\CF l}(\CT)  \ra \CF l(\CT)$.
Thus the choice of a metaplectic structure on $(\CT, k \om)$ uniquely determines, for each real polarization $\CP$, a metalinear structure on $\CP$. 

Given a metalinear frame bundle for each polarization $\CP$, we can now  define the bundle $(\text{Det} \CP^*)^{\frac{1}{2}}$ of half-forms on $\CP$. 
It is the complex line bundle on $\CT$ associated to $ML(\CP)$ through the left-action of $ML(g,\BR)$ on $\BC$, defined for any $\tilde{a} \in ML(g,\BR)$ by multiplication with $\chi(\tilde{a})^{-1}$.
That is,  
\begin{equation}
(\text{Det} \CP^*)^{\frac{1}{2}} \, = \, ML(\CP) \times_{ML(g,\BR)} \BC\, ,    \end{equation}
where a point $[ \underline{\tilde{W}}, \la ]$ in the line bundle $(\text{Det} \CP^*)^{\frac{1}{2}}$ is defined as the equivalence class $[ \underline{\tilde{W}}, \la ] = \{ (\underline{\tilde{W}} \cdot \tilde{a}, \chi(\tilde{a}) \la) \mid \tilde{a} \in ML(g,\BR) \}$.
The space of sections $\mu$ of $(\text{Det} \CP^*)^{\frac{1}{2}}$ is isomorphic to the space of $ML(g,\BR)$-equivariant functions $\mu^{\sharp} : ML(\CP) \ra \BC$, that is, functions satisfying $\, \mu^{\sharp} (\underline{\tilde{W}} \cdot \tilde{a}) = \chi(\tilde{a}) \mu^{\sharp}(\underline{\tilde{W}})$.

We describe now the construction of a metaplectic frame bundle over $\CT$.
Let us consider again the symplectic linear manifold $(\CV, k \om)$.
We choose a symplectic frame $(U_1, \dots ,U_g ; V_1, \dots ,V_g)$ with respect to which we identify the symplectic frame bundle $SP(\CV)$ with the product bundle $\CV \times Sp(2g,\BR)$.
There is only one metaplectic structure on $\CV$ and it corresponds to the trivial metaplectic frame bundle $MP(\CV) = \CV \times Mp(2g,\BR)$ with covering map $p : MP(\CV) \ra SP(\CV)$ defined by $p(X, \tilde{c}) = (X, \rho(\tilde{c}))$, for any $X \in \CV$ and $\tilde{c} \in Mp(2g,\BR)$.
The bundle of symplectic frames $SP(\CT)$ over $\CT$ is also trivializable and we identify it with the trivial bundle $\CT \times Sp(2g,\BR)$ through the trivialization defined by the global invariant symplectic frame $(U_1, \dots ,U_g ; V_1, \dots ,V_g)$.
To define the metaplectic frame bundle $MP(\CT)$ over $\CT$ we lift the action of the integer lattice $\CZ$ on $\CV$ to the bundle $MP(\CV)$ and let $MP(\CT) = MP(\CV)/ \CZ$. 
The $\CZ$-action  on $MP(\CV)$ is determined by a homomorphism $\ep : \CZ \ra \BZ/ 2\BZ$ and is defined by
\begin{equation*}
W \cdot (X, \tilde{c}) \, = \, (X+W\, , \, \tilde{\vep}^{\ep(W)} \tilde{c}) \, ,
\end{equation*}
for any $W \in \CZ$ and $(X,\tilde{c}) \in \CV \times Mp(2g,\BR)$.
Thus we have the commutative diagram:
\begin{equation} \label{e:metf}
\begin{CD}
MP(\CV) @>>> MP^{\ep}(\CT) = MP(\CV)/ \CZ \\
@VVV  @VVV \\
\CV @>>> \CT = \CV/ \CZ
\end{CD}
\end{equation}
where $MP^{\ep}(\CT)$ stands for the metaplectic frame bundle determined by the homomorphism $\ep : \CZ \ra \BZ/2 \BZ$.

For the symplectic linear space $(\CV, k \om)$ we have, by a construction similar to the one outlined previously for $(\CT, k \om)$, the following bundles:
the bundle of Lagrangian frames $\CF l(\CV) = SP(\CV)/N$; 
the bundle of metalinear Lagrangian frames $\tilde{\CF l}(\CV) = MP(\CV)/N$; 
for any constant real polarization $\CP$ of $\CV$, the bundle $\widehat{GL}(\CP) \subset \CF l(\CV)$ of frames of $\CP$, the bundle $\widehat{ML}(\CP) \subset \tilde{\CF l}(\CV)$ of metalinear frames of $\CP$ and the bundle $(\widehat{\text{Det}} \CP^*)^{\frac{1}{2}} \, = \, \widehat{ML}(\CP) \times_{ML(g,\BR)} \BC$ of half-forms on $\CP$.
The $\CZ$-action on $MP(\CV)$ commutes with the right $Mp(2g,\BR)$-action and, therefore, also with the action of the subgroup $N \subset Mp(2g,\BR)$. 
Hence there is an induced $\CZ$-action on $\tilde{\CF l}(\CV) = MP(\CV)/N$ and on any subbundle $\widehat{ML}(\CP) \subset \tilde{\CF l}(\CV)$. Thus we have
\begin{equation*}
\tilde{\CF l}(\CT) = MP(\CT)/N = (MP(\CV)/\CZ)/N = \tilde{\CF l}(\CV) /\CZ
\end{equation*}
That is, $\tilde{\CF l}(\CT)$ is the quotient of $\tilde{\CF l}(\CV)$ under the $\CZ$-action defined, for any $W \in \CZ$ and $ (X, \tilde{c} W) \in \tilde{\CF l}(\CV)$, by
\begin{equation*}
W \cdot (X, \tilde{c} N) = (X +W, \tilde{\vep}^{\ep(W)} \tilde{c} N) \, . 
\end{equation*}

Hence, the metalinear frame bundle $ML(\CP)$ on $\CT$ for the invariant real polarization $\CP$ of $\CT$ induced from the constant polarization $\CP$ of $\CV$ is the $\CZ$-quotient $ML(\CP) = \widehat{ML}(\CP)/\CZ$.
The $\CZ$-action on $\widehat{ML}(\CP) \subset \tilde{\CF l}(\CV)$ induces a $\CZ$-action on the associated bundle of half-forms $(\widehat{\text{Det}} \CP^*)^{\frac{1}{2}} \, = \, \widehat{ML}(\CP) \times_{ML(g,\BR)} \BC$ and
\begin{equation*}
\begin{split}
W \cdot [ (X, \tilde{c} N), \la ] \, &= \, [ W \cdot (X, \tilde{c} N)\, , \la ] \\
&= \, [ (X +W, \tilde{\vep}^{\ep(W)} \tilde{c} N)\, , \la ]  \\
&= \, [ (X +W, \tilde{c} N) \cdot \tilde{\vep}^{\ep(W)}\, , \la ] \qquad ( \tilde{\vep} \, \text{ is central}) \\
&= \, [ (X +W, \tilde{c} N)\, , \chi(\tilde{\vep}^{\ep(W)})^{-1} \la ]  \\
&= \, \te^{- \pi \ti \ep(W)} [ (X +W \, , \tilde{c} N), \la ] \, , \quad \text{ for any } W \in \CZ .
\end{split}
\end{equation*}
It follows that $(\text{Det} \CP^*)^{\frac{1}{2}} \, = \,(\widehat{\text{Det}} \CP^*)^{\frac{1}{2}} /\CZ$ and that the space of sections of the line bundle $(\text{Det} \CP^*)^{\frac{1}{2}}$ over $\CT$ can be identified with the space of $\CZ$-invariant sections of the line bundle $(\widehat{\text{Det}} \CP^*)^{\frac{1}{2}}$ over $\CV$.
A section $\hat{\mu}$ of $(\widehat{\text{Det}} \CP^*)^{\frac{1}{2}} $ corresponds to a $ML(g,\BR)$-equivariant map $\hat{\mu}^{\sharp} : \widehat{ML}(\CP) \ra \BC$  and it is $\CZ$-invariant if the function $\hat{\mu}^{\sharp}$ satisfies
\begin{equation*}
\hat{\mu}^{\sharp}(X+W, \tilde{c} N) \, = \, \te^{- \pi \ti \ep(W)} \,  \hat{\mu}^{\sharp}(X, \tilde{c} N) \, ,
\end{equation*}
for any $W \in \CZ$.

\msk


\subsection{Quantization with half-forms} \label{ss:forms}

Now, let us return to the quantization of $(\CT, k \om)$. 
We construct the Hilbert space $\tilde{\CH}_{\CP}$ from sections of the line bundle $\CL \otimes (\text{Det} \CP^*)^{\frac{1}{2}}$ covariantly constant along the polarization $\CP$.
The operator of covariant differentiation along $\CP$ of sections of the bundle of half-forms $(\text{Det} \CP^*)^{\frac{1}{2}}$ is defined similarly to (\ref{e:cdhd}).
Let $\mu$ be a section of $(\text{Det} \CP^*)^{\frac{1}{2}}$. Then, for any $V \in \CP$, $\, \na^{\CP}_V \mu$ is the section with associated function $(\na^{\CP}_V \mu)^{\sharp} : ML(\CP) \ra \BC$ determined by 
\begin{equation*}
(\na^{\CP}_V \mu)^{\sharp} (\underline{\tilde{X}}) \, = \, V(\mu^{\sharp}(\underline{\tilde{X}} )) \, .
\end{equation*}
In the above expression $\underline{\tilde{X}}$ is a local metalinear frame field covering the local frame $\underline{X} = (X_1, \dots ,X_g)$, where $X_1, \dots ,X_g$ are Hamiltonian vector fields spanning $\CP$ on some contractible open subset of $\CT$.
The covariant differentiation along $V \in \CP$ of a section $s \otimes \mu$ of $\CL \otimes (\text{Det} \CP^*)^{\frac{1}{2}}$ is then defined by
\begin{equation*}
\na^{\CP}_V(s \otimes \mu) = \na_V s \otimes \mu + s \otimes \na^{\CP}_V \, .
\end{equation*}
The space of sections of the line bundle $\CL \otimes (\text{Det} \CP^*)^{\frac{1}{2}}$ over $\CT$ is identified with the space of $\CZ$-invariant sections of the line bundle $\hat{\CL} \otimes (\widehat{\text{Det}} \CP^*)^{\frac{1}{2}}$ over $\CV$.

Let us choose as in (\ref{e:Zbasis}) constant vector fields $ W_1, \dots ,W_g $ on $\CV$ such that $\CP = \text{span}_{\BR} \{ W_1, \dots ,W_g \}$ and $W_i \in \CP \cap \CZ$.
Then the Bohr-Sommerfeld set $\mathcal{BS}_{\CP}$ on $\CT$, that is, the union of all the leaves $\La$ of the polarization $\CP$ for which the line bundle $(\CL \otimes (\text{Det} \CP^*)^{\frac{1}{2}}) \bigr|_{\La}$ has a nowhere zero covariantly constant section, is determined by the condition
\begin{equation} \label{e:BScondition}
\te^{2 \pi \ti k \om(W,X) - \pi \ti \ep(W)} \, = \, 1 \, , \qquad \text{ for any } W \in \CP \cap \CZ.
\end{equation}
The $k^g$ distinct Bohr-Sommerfeld orbits on $\CT$, labeled by the index $\bq$ running through the set $(\BZ/k \BZ)^g$, are described by the equations:
\begin{equation*}
\La_{\bq =(q_1, \dots,q_g)} \, : \quad k \om(W_i,X) = q_i + \frac{1}{2} \ep(W_i)  \pmod{k}, \quad i=1, \dots ,g
\end{equation*}
The Hilbert space $\tilde{\CH}_{\CP}$ is defined as
\begin{equation*}
\tilde{\CH}_{\CP} = \underset{\La \subset  \mathcal{BS}_{\CP}}{\oplus} \, \tilde{S}_{\La} \, ,
\end{equation*}
where $\tilde{S}_{\La}$ is the one-dimensional space of parallel sections of the line bundle $(\CL \otimes (\text{Det} \CP^*)^{\frac{1}{2}} ) \bigr|_{\La}$. The inner product in $\tilde{\CH}_{\CP}$ is defined by
\begin{equation*}
\langle \si , \si' \rangle = 
 \begin{cases}
 0 & , \text{ if } \, \si  \in \tilde{S}_{\La} , \si'  \in \tilde{S}_{\La'} , \, \La \neq \La' \\
\int\limits_{\La} \Om(\si , \si') & , \text{ if } \, \si , \si' \in \tilde{S}_{\La}
\end{cases}
\end{equation*}
The density $\Om(\si , \si')$ on the Bohr-Sommerfeld orbit $\La$ is defined as follows.
For each point $x \in \La$, there is a neighborhood $U$ of $x$ such that both line bundles $\CL$ and $(\text{Det} \CP^*)^{\frac{1}{2}} $ are trivializable.
Then the sections $\si$ and $\si'$ restricted to $U \cap \La$ can be expressed as $\si = s \otimes \mu$ and $\si' = s' \otimes \mu'$, for some covariantly constant sections $s, s'$ of $\CL$ and $\mu , \mu'$ of $(\text{Det} \CP^*)^{\frac{1}{2}} $ over $U \cap \La$.
Let $\underline{W} = (W_1, \dots ,W_g)$ be a frame for $\CP_x$. 
Then we set
\begin{equation*}
\Om(\si , \si') (\underline{W}) \, = \, (s, s')(x) \; \, \overline{\mu^{\sharp}(\underline{\tilde{W}}) } \;  {\mu'}^{\sharp}(\underline{\tilde{W}})
\end{equation*}
where $\underline{\tilde{W}} \in ML(\CP)_x$ is a metalinear frame which projects onto the linear frame $\underline{W} \in GL(\CP)_x$.

Let us also reformulate the BKS pairing in the light of the new definitions based on the choice of a metaplectic frame bundle and the use of half-forms instead of half-densities.
Let $\CP_1$ and $\CP_2$ be two real polarizations of $\CT$.
The corresponding Hilbert spaces are $\tilde{\CH}_{\CP_1} = \underset{\La_1 \subset \mathcal{BS}_{\CP_1}}{\oplus} \tilde{S}_{\La_1}$ and $\tilde{\CH}_{\CP_2} = \underset{\La_2 \subset \mathcal{BS}_{\CP_2}}{\oplus} \tilde{S}_{\La_2}$.
The BKS pairing $\langle \langle \cdot , \cdot \rangle \rangle : \tilde{\CH}_{\CP_2} \times \tilde{\CH}_{\CP_1} \ra \BC$ is defined by setting
\begin{equation} \label{e:BKSpairing}
 \langle \langle \si_2 , \si_1 \rangle \rangle = \te^{- \frac{\pi \ti}{4} h} \; \int\limits_{\La_2 \cap \La_1} \; \Om(\si_2, \si_1) \quad , \quad \text{ for any }  \: \si_2 \in \tilde{S}_{\La_2} , \si_1 \in \tilde{S}_{\La_1} \, ,
\end{equation}
where $h= g - \dim(\CP_1 \cap \CP_2)$ and $ \Om(\si_2, \si_1) $ is the density on $\La_2 \cap \La_1$ defined as follows.
For any point $x \in \La_2 \cap \La_1$, choose a basis $\underline{B} = (\underline{V_2}, \underline{W} ; \underline{V_1} , \underline{T})$ of $T_x(\CT)$ as in (\ref{e:pbas})-(\ref{e:psym}).
Then $\underline{U_1}=(\underline{V_1}, \underline{W})$ is a frame for $\CP_{1x}$ and $\underline{U_2}=(\underline{V_2}, \underline{W})$ is a frame for $\CP_{2x}$.
In a neighborhood of $x$ the sections $\si_1$ and $\si_2$ can be factorized such that $\si_1 = s_1 \otimes \mu_1$ and $\si_2 = s_2 \otimes \mu_2$, with $s_1, s_2$ sections of $\CL$ and $\mu_1, \mu_2$ sections of $(\text{Det} \CP_1^*)^{\frac{1}{2}}$ and $(\text{Det} \CP_2^*)^{\frac{1}{2}}$, respectively.
Then the value of the density  $\Om(\si_2, \si_1)$ on the basis $\underline{W}$ of $T_x(\La_2 \cap \La_1)$ is defined to be
\begin{equation} \label{e:BKSdens}
\Om(\si_2, \si_1) (\underline{W}) \, = \, (s_2,s_1)(x) \; \, \overline{\mu_2^{\sharp}(\underline{\tilde{U_2}}) } \; \mu_1^{\sharp}(\underline{\tilde{U_1}}) \, .
\end{equation}
In the above expression $\underline{\tilde{U_1}}$ is a metalinear frame of $\CP_{1x}$ which projects onto $\underline{U_1}$ and $\underline{\tilde{U_2}}$ a metalinear frame of $\CP_{2x}$ which projects onto $\underline{U_2}$. The frames $\underline{\tilde{U_1}} , \underline{\tilde{U_2}}$ are chosen in a manner 'consistent' with the metaplectic structure on $(\CT, k \om)$.
Following \cite{GS,Bl} we describe below what this 'consistent' choice means.

Assume first that the polarizations $\CP_1$ and $\CP_2$ are transverse, that is $\CP_1 \cap \CP_2 = 0$. Then $\underline{B} = (\underline{U_2};\underline{U_1})$ is a symplectic basis of $(\CT, k \om)$ at $x$.
We arbitrarily pick a metalinear frame $\underline{\tilde{U_2}} \in ML(\CP_2)_x$ projecting onto $\underline{U_2}$. 
Then there exists a unique metaplectic frame $\underline{\tilde{B}} \in MP(\CT)_x$ so that $p(\underline{\tilde{B}}) = \underline{B}$ and $\tilde{\la}(\underline{\tilde{B}}) = \underline{\tilde{U_2}}$, where the maps are those of diagram (\ref{e:cdmetLag}).
Let $j$ be the element of $Sp(2g,\BR)$ represented by the matrix $j = \left( \begin{smallmatrix} 0 & -I \\ I & 0 \end{smallmatrix} \right)$. 
Then $\underline{B} \cdot j = (\underline{U_1};-\underline{U_2})$ and we have $ \la (\underline{B}) = \underline{U_2}$ and $\la(\underline{B} \cdot j) = \underline{U_1}$.   
Now let $c(t)$, with $t \in [0, \frac{\pi}{2} ]$, be the unique one-parameter subgroup of $Mp(2g,\BR)$ projecting onto $c(t) = \left( \begin{smallmatrix} I \cos t & -I \sin t \\ I \sin t & I \cos t \end{smallmatrix} \right) \in Sp(2g,\BR)$ and set $\tilde{j} = \tilde{c}(\frac{\pi}{2})$.
Then $\underline{\tilde{U_1}} = \tilde{\la}(\underline{\tilde{B}} \cdot \tilde{j}) \in ML(\CP_1)_x$ is a metalinear frame for $\CP_1$ projecting onto $\underline{U_1}$, uniquely determined by the choices of a metaplectic frame bundle and metalinear frame $\underline{\tilde{U_2}}$.

For nontransverse polarizations $\CP_{12} = \CP_1 \cap \CP_2 \neq 0$ and $\CP_{12}$ is an isotropic integrable distribution.
The space $\CT/\CP_{12}$ of all integral manifolds of $\CP_{12}$ is a quotient manifold of $\CT$ and we let $\Pi_{12} : \CT \ra \CT/\CP_{12}$ denote the canonical projection. 
The symplectic complement of $\CP_{12}$ is the integrable coisotropic distribution $\CP_{12}^{\perp} = \CP_1 + \CP_2$. For each integral manifold $M$ of $\CP_{12}^{\perp}$, the symplectic form $\om \bigr|_M$ projects to a symplectic form $\om_{M'}$ on the submanifold $M' = \Pi_{12}(M)$ of $\CT/\CP_{12}$. 
The metaplectic structure on $(\CT, k \om)$ induces a metaplectic structure on $(M', k \om_{M'})$ (\cite{Sn3}, \S 5.4). 
Then, making use of this induced metaplectic structure, we are able, just as in the case of transverse polarizations, to uniquely pick the metalinear frame $\underline{\tilde{U_1}}$ projecting onto $\underline{U_1}$ for a given a frame $\underline{\tilde{U_2}}$ over $\underline{U_2}$.

Finally one can check that the right-hand side of the expression (\ref{e:BKSdens}) does not depend on the choice of $\underline{\tilde{U_2}}$ covering $\underline{U_2}$ and that $\Om(\si_2,\si_1)$ does indeed define a density on $T(\La_2 \cap \La_1)$.

\msk


\subsection{The Hilbert space}

We are going to quantize $(\CT, k \om)$ with the trivial metaplectic structure. The bundle of metaplectic frames on $\CT$ defined by (\ref{e:metf}) is the product bundle $MP(\CT) = \CT \times Mp(2g,\BR)$. 
Since the metaplectic frame bundle is trivial so are the metalinear frame bundle and the bundle of half-forms of any polarization of $\CT$.
The Bohr-Sommerfeld condition (\ref{e:BScondition}) reads now as in Sect.\ref{s:quanttg}.
However, as we shall prove below, the use of half-forms in the construction of the quantum Hilbert space leads to intertwining operators between Hilbert spaces of different polarizations which obey a transitive composition law. 

For each polarization $\CP$ we fix a symplectic frame $(W_i\, ;W_i^{\perp})$ of $(\CT, \om)$ as in (\ref{e:Zbasis}), with $\CP = \text{span}_{\BR} \{ W_1,\dots,W_g \}$, and a metalinear frame field $\underline{\tilde{W}} \in ML(\CP)$ projecting onto the linear frame field $\underline{W}=(W_1,\dots,W_g)$ of $\CP$. This uniquely determines a {\em standard unitary basis} $\{ \si_{\bq} = s_{\bq} \otimes \de_{\bq} \}_{\bq \in (\BZ / k \BZ)^g}$ of the Hilbert space $\tilde{\CH}_{\CP}$, where: \\
$\bullet \quad s_{\bq}$ is the unitary section of $\CL$ on the Bohr-Sommerfeld leave $\La_{\bq}$, as described by the equations (\ref{e:sigq})-(\ref{e:sigql}); \\
$\bullet \quad \de_{\bq}$ is the restriction to $\La_{\bq}$ of the section $\de$ of $(\text{Det} \CP^*)^{\frac{1}{2}}$ defined by setting $\de^{\sharp}(\underline{\tilde{W}}) = 1$.

For any two polarizations $\CP_1$ and $\CP_2$, the BKS pairing (\ref{e:BKSpairing}) induces the linear operator $\tilde{F}_{\CP_2 \CP_1} : \tilde{\CH}_{\CP_1} \ra \tilde{\CH}_{\CP_2}$ defined by $\langle \si_2, \tilde{F}_{\CP_2 \CP_1} \si_1 \rangle \, = \, \langle \langle \si_2, \si_1 \rangle \rangle$, for any $\si_1 \in \tilde{\CH}_{\CP_1}$ and $\si_2 \in \tilde{\CH}_{\CP_2}$. 
With respect to the unitary bases $\{ \si_{1 \bq_1} = s_{1 \bq_1} \otimes \de_{1 \bq_1} \}$ of $\tilde{\CH}_{\CP_1}$ and $\{ \si_{2 \bq_2} = s_{2 \bq_2} \otimes \de_{2 \bq_2} \}$ of $\tilde{\CH}_{\CP_2}$, the operator $\tilde{F}_{\CP_2 \CP_1}$ is represented by the matrix
\begin{equation*}
\widetilde{\bM}(2,1)_{\bq_2 \bq_1} \, = \, \te^{- \frac{\pi \ti}{4} h} \; \int_{\La_{2 \bq_2} \cap \La_{1 \bq_1}} \; \Om(\si_{2 \bq_2} , \si_{1 \bq_1} ) 
\end{equation*}
where $h = g - \dim(\CP_1 \cap \CP_2)$.
Let us compute the explicit form of the matrix $\widetilde{\bM}(2,1)$. 
Choose linear frames $\underline{W_1} \in GL(\CP_1)$ and $\underline{W_2} \in GL(\CP_2)$ as in (\ref{e:intbas})-(\ref{e:intbases}) and metalinear frame fields $\underline{\tilde{W_1}} \in ML(\CP_1)$ projecting onto $\underline{W_1}$ and $\underline{\tilde{W_2}} \in ML(\CP_2)$ projecting onto  $\underline{W_2}$. 
The half-forms $\de_1$ and $\de_2$ are defined by $\de_1^{\sharp}(\underline{\tilde{W_1}}) =1$ and $\de_2^{\sharp}(\underline{\tilde{W_2}}) =1$. 
Let us introduce,  as in (\ref{e:ntrbasis}), the symplectic frame $\underline{B} = (\underline{V_2}, \underline{W} ; \underline{V_1}, \underline{T})$ for $(\CT, k \om)$ satisfying the conditions  (\ref{e:pbas})-(\ref{e:psym}).
Let $a_{21}$ be the element of $GL(g,\BR)$ defined by
\begin{equation*}
a_{21} \, = \, 
\begin{cases} 
      \begin{pmatrix}
       k \check{\om}(2,1) & 0 \\
        0 & I
       \end{pmatrix} &
   , \qquad \text{if}  \quad \CP_1 \cap \CP_2 \neq 0 \\
   k \om(2,1) &  , \qquad \text{if} \quad \CP_1 \cap \CP_2 = 0 
\end{cases}
\end{equation*}
Then $\underline{U_2} = (\underline{V_2}, \underline{W}) = \underline{W_2}$ is a frame for $\CP_2$ and $\underline{U_1} = (\underline{V_1}, \underline{W}) = \underline{W_1} \cdot a_{21}^{-1}$ is a frame for $\CP_1$.
Let us take $\underline{\tilde{U_2}} = \underline{\tilde{W_2}} \in ML(\CP_2)$. Then, as mentioned previously, this uniquely picks a metalinear frame $\underline{\tilde{U_1}} \in ML(\CP_1)$ projecting onto $\underline{U_1}$.
There is a unique element $\tilde{a}_{21} \in ML(g,\BR)$ so that $\underline{\tilde{W_1}} = \underline{\tilde{U_1}} \cdot \tilde{a}_{21}$ and its projection to $GL(g,\BR)$ is $\rho(\tilde{a}_{21}) = a_{21}$.
For any  pair of intersecting Bohr-Sommerfeld orbits $\La_{1 \bq_1}$ and $\La_{2 \bq_2}$, 
the pair of frames $\underline{\tilde{U_2}}$ and $\underline{\tilde{U_1}}$ is used to define the density $\Om(\si_{2 \bq_2} , \si_{1 \bq_1})$ according to (\ref{e:BKSdens}).
Thus we obtain
\begin{equation*}
\begin{split}
\Om(\si_{2 \bq_2}  , \si_{1 \bq_1})(\underline{W}) &= (s_{2 \bq_2},s_{1 \bq_1}) \; \, \overline{\de_{2 \bq_2}^{\sharp} (\underline{\tilde{U_2}}) } \; \de_{1 \bq_1}^{\sharp}(\underline{\tilde{U_1}}) \\
&= (s_{2 \bq_2},s_{1 \bq_1}) \; \, \overline{\de_{2 \bq_2}^{\sharp} (\underline{\tilde{W_2}}) } \; \de_{1 \bq_1}^{\sharp}(\underline{\tilde{W_1}} \cdot \tilde{a}_{21}^{-1}) \\
&= \chi(\tilde{a}_{21})^{-1} \; (s_{2 \bq_2},s_{1 \bq_1}) 
\end{split}
\end{equation*}
Together with the results of Sect.\ref{s:BKS} the above expression implies that the matrix    $\widetilde{\bM}(2,1)$ differs from the expression of the matrix $\bM(2,1)$ of that section by a phase factor:
\begin{equation} \label{e:phase}
\widetilde{\bM}(2,1) \, = \,  \te^{- \frac{\pi \ti}{4} h} \frac{\chi(\tilde{a}_{21})^{-1} }{| \det a_{21} |^{- \frac{1}{2}} } \:  \bM(2,1) \, .
\end{equation}
Therefore the conclusion at the end of Sect.\ref{s:BKS} is still valid, that is, 
\begin{thm} \label{t:unitary}
For any two polarizations $\CP_1$ and $\CP_2$ of the symplectic torus $(\CT, k \om)$,
the linear operator $\tilde{F}_{\CP_2 \CP_1} : \tilde{\CH}_{\CP_1} \ra \tilde{\CH}_{\CP_2}$ is unitary.
\end{thm}
\nin However, the Theorem (\ref{t:trans-proj}) of Sect.\ref{s:prHsp} is restated as follows:

\begin{thm} \label{t:transit}
For any three invariant real polarizations $\CP_1, \CP_2, \CP_3$ of the symplectic torus $(\CT, k \om)$ we have
\[ \tilde{F}_{\CP_1 \CP_3} \circ \tilde{F}_{\CP_3 \CP_2} \circ \tilde{F}_{\CP_2 \CP_1} \, = \,  I \]
\end{thm}

\begin{proof}
The proof is the same as for Theorem (\ref{t:trans-proj}), the only modification being the appropriate factor correction indicated in (\ref{e:phase}). 
Therefore, using (\ref{e:finsum}), we find that, with the assumption that $\CP_3$ is transverse to $\CP_1$ and $\CP_2$, we have
\begin{align} \label{e:comp3}
& \sum_{\bq_{2}, \bq_{3} \in (\BZ / k \BZ)^{g} } \: \widetilde{\bM}(1,3)_{\bq'_{1} \bq_{3}} \, \widetilde{\bM}(3,2)_{\bq_{3} \bq_{2}} \, \widetilde{\bM}(2,1)_{\bq_{2} \bq_{1}}  \, = \\
& = \frac{ \te^{- \frac{\pi \ti}{4} h} \, \chi( \tilde{a}_{32}^{-1} \tilde{a}_{31} \tilde{a}_{21}^{-1} ) }{ | \det a_{32} |^{- \frac{1}{2}} 
| \det a_{31} |^{\frac{1}{2}} | \det a_{21} |^{- \frac{1}{2}} } \: \te^{-\frac{\pi \ti}{4} \tau(L_{\CP_1}, L_{\CP_2}, L_{\CP_3}) } \: \de_{\bq'_1 \bq_1} \notag
\end{align}
Under the projection $\rho : ML(g,\BR) \ra GL(g,\BR)$ we have $\rho(\tilde{a}_{31}) = a_{31}, \,  \rho(\tilde{a}_{32}) = a_{32}$ and $\rho(\tilde{a}_{21}) = a_{21}$.
Making use of the results of Sect.\ref{s:prHsp} we find the following.
If $\CP_1 \cap \CP_2 = 0$ then 
\begin{equation*}
a_{32}^{-1} a_{31} a_{21}^{-1} \, = \, k^{-1} \om(3,2)^{-1} \om(3,1) \om(2,1)^{-1} \, = \, - k^{-1} H^{-1} \, ,
\end{equation*}
where $H$ is the symmetric matrix whose signature is the Maslov-Kashiwara index $\tau(L_{\CP_1}, L_{\CP_2}, L_{\CP_3})$.
If $\CP_1 \cap \CP_2 \neq 0$ then 
\begin{equation*}
a_{32}^{-1} a_{31} a_{21}^{-1}  \, = \, \om(3,2)^{-1} \om(3,1) \begin{pmatrix} k^{-1} \check{\om}(2,1)^{-1} & 0 \\ 0 & I \end{pmatrix} \, = \,
\begin{pmatrix}  k^{-1} B \check{\om}(2,1)^{-1} & 0 \\ 0 & I \end{pmatrix}
\end{equation*}
and $H = - \left( \begin{smallmatrix} \check{\om}(2,1)B & 0 \\ 0 & 0 \end{smallmatrix} \right)$.
Thus, we obtain
\begin{equation*}
\frac{ \te^{- \frac{\pi \ti}{4} h} \, \chi( \tilde{a}_{32}^{-1} \tilde{a}_{31} \tilde{a}_{21}^{-1} ) }{ | \det a_{32} |^{- \frac{1}{2}} 
| \det a_{31} |^{\frac{1}{2}} | \det a_{21} |^{- \frac{1}{2}} } \, = \, 
\te^{- \frac{\pi \ti}{4} h} \, \te^{\frac{\pi \ti}{2} h_{+}} \, = \, \te^{\frac{\pi \ti}{4} \text{sgn} H} \, ,
\end{equation*}
where $h_{+}$ equals the number of positive eigenvalues of the symmetric matrix $H$.
For the last equality we used the relation $\text{sgn} H = 2 h_{+} -g + \dim(\CP_1 \cap \CP_2) = 2 h_{+} - h$.
Thus the phase factors on the right-hand side of equation (\ref{e:comp3}) cancel out. 
\end{proof}
Hence quantization canonically associates to the data $(\CT, k \om)$ plus a trivial metaplectic frame bundle a Hilbert space $\tilde{\CH}$.
This is similar to the quantization of a symplectic vector space \cite{GS,Wo}, where the metaplectic correction resolves the projective ambiguity from the quantization with half-densities and enables the construction of a Hilbert space.
One can carry out the constructions of the present section (Hilbert space, BKS pairing) equally well for a choice of a nontrivial metaplectic structure on $(\CT,  k \om)$. The computations necessary for proving the Theorems (\ref{t:unitary}) and (\ref{t:transit}) are more involved, but the conclusions are the same.

\bsk


\section{The Maslov index and the Hilbert space} \label{s:Maslov}

In this section we present an alternative way to the metaplectic correction of resolving the projective ambiguity which results from the quantization with half-densities described in Sect.\ref{s:prHsp}.
It exploits the relationship between the Maslov-Kashiwara index of a triple of Lagrangian planes of a symplectic vector space and the Maslov index of a pair of elements in a $4$-fold cover of the Lagrangian Grassmannian of that space.
The idea is not new; the same facts are used in the construction of the Shale-Weil representation of the metaplectic group on the Hilbert space of quantization of a symplectic space \cite{GS,LV}.

\msk

We begin by summarizing some results regarding the covering groups of the symplectic group, the covering spaces of the Lagrangian Grassmannian of a symplectic vector space and the Maslov index. 
For more details and proofs we refer the reader to \cite{Go, GS,Le}.

Consider again the symplectic $2g$-dimensional vector space $(\CV,\om)$.
Let $\lag$ denote the Lagrangian Grassmannian of $(\CV,\om)$, that is, the space of all Lagrangian subspaces of $\CV$.
The symplectic group $Sp(\CV)$ acts transitively on $\lag$.
The isotropy group at any point $L$ in $\lag$ is a closed subgroup $S$ of $Sp(\CV)$ isomorphic to $GL(g,\BR) \times \BR^{\frac{g(g+1)}{2}}$.
Let $\al$ and $\be$ denote the generators of $\pi_1(Sp(\CV)) \cong \BZ$ and $\pi_1(\lag) \cong \BZ$, respectively, whose natural images in $\BZ$ are 1.
The fibration $ Sp(\CV) \ra Sp(\CV)/S = \lag$ defines a monomorphism $\pi_1(Sp(\CV)) \ra \pi_1(\lag)$ which sends $\al$ to $\be^2$  \cite{GS,Le}.
We have (\cite{Le}, ch.I):
\begin{thm} \label{t:cover}
(1) For every $q=1, 2, \dots , \infty$ the symplectic group has a unique $q$-fold covering group $Sp_q(\CV)$; 
$Sp_{\infty}(\CV)$ is the universal covering group. 
The generator $\al$ of $\pi_1(Sp(\CV))$ acts on $Sp_q(\CV)$;
$\al^r$ does not act as the identity on $Sp_q(\CV)$ unless $r \equiv \pmod{q}$. \\
(2) For every $q=1, 2, \dots , \infty$ the Lagrangian Grassmannian $\lag$ has a unique $q$-fold covering space $Lag_{q}(\CV)$;
$Lag_{\infty}(\CV)$ is the universal covering space.
The generator $\be$ of $\pi_1(\lag)$ acts on $Lag_{q}(\CV)$;
$\be^r$ does not act as the identity on $Lag_{q}(\CV)$ unless $r \equiv \pmod{q}$. \\
(3) the group $Sp_q(\CV)$ acts transitively on $Lag_{2q}(\CV)$ and 
\[ (\al \cdot a_q) \cdot L_{2q} = a_q \cdot ( \be^2 \cdot L_{2q}) = \be^2 \cdot ( a_q \cdot L_{2q}) \, ,  \]
for any $a_q \in Sp_q(\CV) $ and $ L_{2q} \in Lag_{2q}(\CV)$.
\end{thm}
For our purpose the important fact is that the double cover $Mp(\CV) = Sp_2(\CV)$ of $Sp(\CV)$, the metaplectic group, acts on the 4-fold cover $Lag_4 (\CV)$ of $\lag$.

Let us consider the universal covering space $Lag_{\infty}(\CV)$  of $\lag$ and for any $L_{\infty} \in Lag_{\infty}(\CV)$ let $L$ denote its natural projection onto $\lag$. A pair $(L_{\infty}, L_{\infty}')$ of elements in $Lag_{\infty}(\CV) \times Lag_{\infty}(\CV)$ is called {\em transverse} if $L \cap L' =0$.
We will make use of the following results \cite{Go}:
\begin{thm} \label{t:thmMi}
There exists a unique function
\[ \mu \, : \, Lag_{\infty}(\CV) \times Lag_{\infty}(\CV) \lra \BZ \]
called the Maslov index, with the following two properties: \\
(1) $\mu(L_{\infty}, L_{\infty}') - \mu(L_{\infty}, L_{\infty}'') + \mu(L_{\infty}', L_{\infty}'') \, = \, \tau(L, L', L'')$;\\
(2) $\mu(L_{\infty}, L_{\infty}') - \tau(L, L', L'')$ is locally constant on the subset $\{ (L_{\infty}, L_{\infty}', L'') \in Lag_{\infty}(\CV) \times Lag_{\infty}(\CV) \times \lag \mid L \cap L'' = L' \cap L'' = 0 \}$.
\end{thm}
\nin where $\tau(L, L', L'')$ in (\ref{t:thmMi}-(1)) is the Maslov-Kashiwara index  of a triple $L, L', L''$ of Lagrangian planes in $(\CV, \om)$ introduced  in Sect.\ref{s:prHsp}.
\begin{prop} \label{p:propMi}
The Maslov index $\mu$ has the following properties: \\
(1) $\mu(L_{\infty}, L_{\infty}') + \mu(L_{\infty}', L_{\infty})  = 0$; \\
(2) $\mu(L_{\infty}, L_{\infty}') \equiv g - \dim(L \cap L') \, \pmod{2}$ ; \\
(3) $\mu(\be^r \, L_{\infty}, \be^{r'} \, L_{\infty}')  = \mu(L_{\infty}, L_{\infty}') + 2r -2r'$; \\
(4) $\mu(a_{\infty} \, L_{\infty}, a_{\infty} \, L_{\infty}') = \mu(L_{\infty}, L_{\infty}') $, for every $a_{\infty} \in Sp_{\infty}(\CV)$.
\end{prop}

For every positive integer $q$, the Maslov index $\mu$ on $Lag_{\infty}(\CV) \times Lag_{\infty}(\CV)$ induces a $\BZ/ 2q \BZ$-valued Maslov index $\mu_{2q}$ on $Lag_{q}(\CV) \times Lag_{q}(\CV)$ as follows.
Let $(L_{\infty}, L_{\infty}')\, , \, (K_{\infty}, K_{\infty}')$ be two pairs of elements in  $Lag_{\infty}(\CV) \times Lag_{\infty}(\CV)$ with the same projection $(L_q,L_{q}')$ to $Lag_{q}(\CV) \times Lag_{q}(\CV)$. 
From Theorem (\ref{t:cover}) it follows that there exists $(r,r') \in \BZ^2$ such that
\begin{equation*}
L_{\infty} = \be^r K_{\infty} \; , \; L_{\infty}' = \be^{r¹} K_{\infty}' \quad
\text{and} \quad r \equiv r' \equiv \pmod{q} .
\end{equation*}
Then, in view of Proposition (\ref{p:propMi}), we have
\begin{equation*}
\mu(L_{\infty}, L_{\infty}') = \mu(K_{\infty}, K_{\infty}') \, \pmod{2q} .
\end{equation*}
Thus there is a well defined Maslov index $\mu_{2q}$ on $Lag_{q}(\CV) \times Lag_{q}(\CV)$  given by
\begin{equation}
\mu_{2q}(L_q, L_q') = \mu(L_{\infty}, L_{\infty}') \, \pmod{2q},
\end{equation}
where $(L_{\infty}, L_{\infty}')$ is any element of $Lag_{\infty}(\CV) \times Lag_{\infty}(\CV)$ projecting onto $(L_q , L_q')$. Similarly to (\ref{t:thmMi}) and (\ref{p:propMi}) we have:
\begin{thm} \label{t:thmMi2q}
The Maslov index $\mu_{2q}$ is the only function
\[ \mu_{2q} : Lag_2(\CV) \times Lag_2(\CV) \lra \BZ/ 2q \BZ \]
with the following two properties: \\
(1) $ \mu_{2q}(L_q, L_q') - \mu_{2q}(L_q, L_q'') + \mu_{2q}(L_q', L_q'') = \tau_{2q}(L, L', L'') $,
where $\tau_{2q}(L, L', L'') = \tau(L, L', L'') \pmod{2q}$; \\
(2) $\mu_{2q}(L_q, L_q') - \tau_{2q}(L, L', L'') $is locally constant on the subset $\{ (L_q, L_q',L'') \in Lag_{q}(\CV) \times Lag_{q}(\CV) \times Lag(\CV) \mid L \cap L'' = L' \cap L'' = 0 \}$.
\end{thm}
\begin{prop}
The Maslov index $\mu_{2q}$ has the following properties: \\
(1) $\mu_{2q}(L_q, L_q') + \mu_{2q}(L_q', L_q)  = 0$; \\
(2) $\mu_{2q}(\be^r \, L_q, \be^{r'} \, L_q¹)  = \mu_{2q}(L_q, L_q') + 2r -2r' \, \pmod{2q}$; \\
(3) $\mu_{4q} : Lag_{2q}(\CV) \times Lag_{2q}(\CV) \ra \BZ$ is invariant under the action of $Sp_q(\CV)$ on $Lag_{2q}(\CV)$, that is $\mu_{4q}(a_q \, L_{2q}, a_q \, L_{2q}') = \mu_{4q}(L_{2q}, L_{2q}') $, for every $a_q \in Sp_q(\CV)$.
\end{prop}

We introduce the notation $\tlag = Lag_4(\CV) \,  , \,  \tilde{\mu} = \mu_8$ and let $\tilde{L}$ denote a point in $\tlag$ which projects onto $L$ in $\lag$.
As for the metaplectic group $Sp_2(\CV)$ we use the standard notation $Mp(\CV)$.

\msk

In Sect.\ref{s:prHsp} we proved that for any three polarizations $\CP_1, \CP_2 , \CP_3$ of the symplectic torus $(\CT, k \om)$ we have
\begin{equation*}
F_{\CP_1 \CP_3} \circ F_{\CP_3 \CP_2} \circ F_{\CP_2 \CP_1}  \, = \, \te^{- \frac{2 \pi \ti}{8} \tau(L_{\CP_1}, L_{\CP_2}, L_{\CP_3}) } \;  I  ,
\end{equation*}
where $L_{\CP_1}, L_{\CP_2}, L_{\CP_3}$ are the corresponding rational Lagrangian subspaces of $(\CV, \om)$.
According to (\ref{t:thmMi2q}) the lift of the Maslov-Kashiwara index $\tau$ to $\tlag$ has the decomposition:
\begin{equation} \label{e:cobMi}
\tau(L_{\CP_1}, L_{\CP_2}, L_{\CP_3}) = \tilde{\mu}(\tilde{L}_{\CP_1}, \tilde{L}_{\CP_2}) +\tilde{\mu}(\tilde{L}_{\CP_2}, \tilde{L}_{\CP_3}) +
\tilde{\mu}(\tilde{L}_{\CP_3}, \tilde{L}_{\CP_2}) \, \pmod{8} .
\end{equation}
This leads us to define, for any polarization $\CP$ and for any point $\tilde{L}_{\CP} \in \tlag$ projecting onto $L_{\CP} \in \lag$, the Hilbert space $\CH_{\tilde{L}_{\CP}} = \CHP$.
Then we let 
\begin{equation} \label{e:Fproj}
F_{\tilde{L}_{\CP_2} \tilde{L}_{\CP_1} } \, : \, \CH_{\tilde{L}_{\CP_2}} \lra \CH_{\tilde{L}_{\CP_1}}
\end{equation}
denote the unitary operator given by
\begin{equation*}
F_{\tilde{L}_{\CP_2} \tilde{L}_{\CP_1} } \, = \, \te^{- \frac{\pi \ti}{4} \tilde{\mu}(\tilde{L}_{\CP_1}, \tilde{L}_{\CP_2}) } \, F_{\CP_2 \CP_1}\, .
\end{equation*}
In view of (\ref{e:Fproj}) and (\ref{e:cobMi}) we have the transitive composition law:
\begin{equation}
F_{\tilde{L}_{\CP_1} \tilde{L}_{\CP_3}} \circ F_{\tilde{L}_{\CP_3} \tilde{L}_{\CP_2}} \circ F_{\tilde{L}_{\CP_2} \tilde{L}_{\CP_1}}  \, = \,  I  
\end{equation}
Thus the Hilbert spaces  $\CH_{\tilde{L}_{\CP}}$, for all $\tilde{L}_{\CP} \in \tlag$ and all polarizations $\CP$ of $\CT$, are canonically identified.
Geometric quantization associates to the symplectic torus $(\CT, k \om)$ together with a choice of 4-fold cover $\tlag$ of the Lagrangian Grassmannian $\lag$ a Hilbert space $\CH$.

\bsk


\section{The representation of the integer metaplectic group}
\label{s:repres}

In this section we show that the integer symplectic group has a unitary projective representation on the projective Hilbert space of Sect.\ref{s:quanttg} and that the integer metaplectic group has a unitary representation on the Hilbert spaces constructed in Sect.\ref{s:metapl} and Sect.\ref{s:Maslov}.

\subsection{Projective representation of the integer symplectic group} \label{ss:proj}

We recall from Sect.\ref{s:quanttg} that the elements of the integer symplectic group $Sp(\CZ)$ act by symplectic diffeomorphisms on the torus $(\CT, k \om)$ and that this action lifts to the prequantum line bundle $\CL$, preserving the connection and the hermitian metric.

Let $\CP$ be a polarization of $\CT$ corresponding to a rational Lagrangian plane $L_{\CP}$ in the vector space $\CV$.
An element $b$ in $Sp(\CZ)$ maps $\CP$ to a polarization $b \CP$ corresponding to the rational Lagrangian plane $b L_{\CP}$ of $\CV$.
Hence the action of $b$ on $\CT$ lifts to a map of line bundles:
\begin{equation*}
\begin{CD}
\CL \otimes \dephs @>b>> \CL \otimes {\left| \text{Det} \, (b {\mathcal P})^{*} \right|}^{\frac{1}{2}} \\
@VVV     @VVV \\
\CT @>b>> \CT
\end{CD}
\end{equation*}
This gives rise to a map of sections. If $s \otimes \mu$ is a section of $\CL \otimes \dephs$, then $b \cdot (s \otimes \mu) = (b \cdot s) \otimes (b \cdot \mu)$ is the section of $\CL \otimes {\left| \text{Det} \, (b{\mathcal P})^{*} \right|}^{\frac{1}{2}}$ defined at any $x \in \CT$ by
\begin{align*}
(b \cdot s)(x) \,  &= \, b \cdot s(b^{-1}(x)) \\
(b \cdot \mu)_{b(x)}(bW_1, \dots ,bW_g)\, &= \, \mu_x(W_1, \cdots ,W_g) \, ,
\end{align*}
where $W_1, \dots , W_g $ are vectors spanning $\CP_x$.
Thus, to each $b \in Sp(\CZ)$, there corresponds a unitary operator $b: \CHP \ra \CH_{b \CP}$.
Its composition with the intertwining isomorphism $F_{\CP, b \CP} : \CH_{b \CP} \ra \CHP$ gives a unitary operator 
\begin{equation}
U_{\CP}(b) = F_{\CP, b \CP} \circ b \, : \, \CHP \ra \CHP
\end{equation}
For any two polarizations $\CP_1, \CP_2$ we have the commutative diagram
\begin{equation*}
\begin{CD}
\CH_{\CP_1} @>{F_{\CP_2 \CP_1}}>> \CH_{\CP_2} \\
@VbVV     @VVbV \\
\CH_{b \CP_1} @>{F_{b \CP_2 , b \CP_1}}>> \CH_{b \CP_2}
\end{CD}
\end{equation*}
Using this and Theorem (\ref{t:trans-proj})  we find that
\begin{equation*}
\begin{split}
U_{\CP}(b) \, U_{\CP}(b') \, &= \, F_{\CP, b \CP} \circ b \circ F_{\CP, b' \CP} \circ b' \\
&= \, F_{\CP, b \CP} \circ F_{b \CP, b b' \CP} \circ b b' \\
&= \te^{\frac{\pi \ti}{4} \tau(L_{\CP}, b L_{\CP}, b b' L_{\CP}) } \, F_{\CP, b b' \CP} \circ b b' \\
&= \te^{\frac{\pi \ti}{4} \tau(L_{\CP}, b L_{\CP}, b b' L_{\CP}) } U_{\CP}(b b')
\end{split}
\end{equation*}
for any $b, b' \in Sp(\CZ)$.
That is, the map $b \ra U_{\CP}(b)$ from $Sp(\CZ)$ into the set of unitary operators on $\CHP$ defines a projective representation of the group $Sp(\CZ)$ with associated cocycle $c_{\CP}(b, b') = \te^{\frac{\pi \ti}{4} \tau(L_{\CP}, b L_{\CP}, b b' L_{\CP}) }$.
Hence we have:
\begin{thm}
The group $Sp(\CZ)$ has a unitary projective representation on the projective Hilbert space $\BP \CH$ associated through quantization to $(\CT, k \om)$.
\end{thm}
We determine now the explicit form of this representation. 
Let us fix a polarization $\CP$ and a symplectic frame $( W_i \, ; W_i^{\perp})$ of $(\CT, k \om)$ as in (\ref{e:Zbasis}), with $\CP = \text{span}_{\BR} \{ W_1, \dots ,W_g \}$.
With respect to the basis $( W_i \, ; W_i^{\perp} )$ we identify the group $Sp(\CZ)$ with the symplectic group $Sp(2g, \BZ)$ with $\BZ$-coefficients.
The generators of $Sp(\CZ)$ are the elements with matrix form
\begin{equation} \label{e:Spgen}
\al = \begin{pmatrix} A & 0 \\ 0 & \,{ ^t A^{-1}} \end{pmatrix}
\quad
\be = \begin{pmatrix} I & B \\ 0 & I \end{pmatrix} \quad
\ga = \begin{pmatrix} 0 & I \\ -I & 0 \end{pmatrix}
\end{equation}
with $A \in GL(g, \BZ)$ and $ B \in M(g, \BZ), \, \,{ ^t B} = B$.
As mentioned in Sect.\ref{s:quanttg}, the choice of basis $( W_i \, ; W_i^{\perp} )$ uniquely determines a standard unitary basis $\{ \si_{\bq} = s_{\bq} \otimes \de_{\bq} \}_{\bq \in (\BZ /k \BZ)^g}$ of $\CHP$ and we have
\begin{prop} \label{p:Sprepr}
The projective representation of $Sp(\CZ)$ on $\CHP$ is described by the following unitary operators representing the generators:
\begin{align*}
U_{\CP}(\al) \, \si_{\bq} \, &=  \, \si_{^t \! A^{-1} \bq} \\
U_{\CP}(\be)  \, \si_{\bq} \, &= \,  \te^{\frac{\pi \ti}{k} \,{ ^t \bq} B \bq} \, \si_{\bq} \\
U_{\CP}(\ga) \, \si_{\bq} \, &= \,   k^{- \frac{g}{2}} \sum_{\bq_1 \in (\BZ /k \BZ)^g} \, \te^{\frac{2\pi \ti}{k} \,{ ^t \bq} \bq_1} \, \si_{\bq_1}
\end{align*}
\end{prop}
\begin{proof}
{\em Case 1}. Assume that $b \in Sp(\CZ)$ maps $\CP$ to a polarization $\CP' = b \CP$ and $\CP' \neq \CP$.
The symplectic basis $( W_i' = b W_i \, ; W_i^{'\perp} = b W_i^{\perp} )$ uniquely determines a standard unitary basis $\{ \si_{\bq}' = s_{\bq}' \otimes \de_{\bq}' \}$ for $\CH_{\CP'}$.
We can easily deduce that the symplectic potential $\theta_{\CP'}$ defined by the relations
\begin{align*}
(W_i' \contrac \theta_{\CP'})_X \, &= \, 0 \\
(W_i^{'\perp} \contrac \theta_{\CP'})_X \, &= \, \om(X, W_i^{'\perp})
\end{align*}
satisfies $b^{*} \theta_{\CP'} = \theta_{\CP}$.
From the above expression together with the invariance $b^* \theta_0 = \theta_0$ and the defining relations $\theta_{\CP} = \theta_0 + d K_{\CP} \, , \, \theta_{\CP'} = \theta_0 + d K_{\CP'}$ and $K_{\CP}(0) = K_{\CP'}(0) = 0$ for the functions $K_{\CP}$ and $K_{\CP'}$, we find that
\begin{equation*}
K_{\CP'} \cdot b - K_{\CP} \, = \, 0\,.
\end{equation*}
Now, a direct application of the definition (\ref{e:sigq})-(\ref{e:sigql}) of the section $s_{\bq}$ of $\CL$ and of the definition (\ref{e:SonL}) of the $Sp(\CZ)$-action leads to the conclusion that 
\begin{equation*}
b \cdot s_{\bq} \, = \, s_{\bq}' \, ,
\end{equation*}
where $s_{\bq}'$ is the section of $\CL$ over the Bohr-Sommerfeld leave $\La_{\bq}'$ described by the equations: $ \, k \om(W_i',X) = q_i \pmod{k}$.
Also, due to the invariance of the canonical $\frac{1}{2}$-density $\de$ under $Sp(\CZ)$ we have, letting $\de_{\bq}' = \de \bigr|_{\La_{\bq}'}$,
\begin{equation*}
b \cdot \de_{\bq} \, = \, \de_{\bq}' \, .
\end{equation*}
Let us now specialize the above formulas to the element $\ga = \left( \begin{smallmatrix} 0 & I \\ -I & 0 \end{smallmatrix} \right)$ and $\CP' = \ga \CP$.
Then $( W_i' = - W_i^{\perp} \, ; W_i^{'\perp} = W_i )$ and $\CP' \cap \CP = 0$. 
Moreover, any two Bohr-Sommerfeld leaves $\La_{\bq}$ of $\CP$ and $\La_{\bq'}'$ of $\CP'$ intersect in only one point. 
With the appropriate identifications, the expression of the matrix  representing the operator $F_{\CP \CP'} : \CH_{\CP'} \ra \CHP$ with respect to the standard unitary bases $\{ \si_{\bq}' \}$ and $\{ \si_{\bq} \}$ can be read from (\ref{e:trmat}). 
We find that $U_{\CP}(\ga) \si_{\bq} \, = \, (F_{\CP \CP'} \circ \ga) (\si_{\bq}) \, = \, F_{\CP \CP'} (\si_{\bq}')$ has the form show in (\ref{p:Sprepr}).

\nin {\em Case 2}. Assume that the element $b \in Sp(\CZ)$ preserves the polarization $\CP$, that is, $b \CP = \CP$.
Then $b$ must have the following matrix form with respect to the basis $( W_i \, ; W_i^{\perp} )$:
\begin{equation*}
b \, = \, \begin{pmatrix} A & B \\ 0 & \,{ ^t A^{-1}} \end{pmatrix}
\qquad \begin{matrix} A \in GL(g, \BZ) \, , \, B \in M(g,\BZ) \\
            \,{ ^t (A^{-1} B)} = A^{-1} B \end{matrix}
\end{equation*}
Also, for any $W \in \CP$, we have
\begin{equation*}
W \contrac \, d(K_{\CP} \cdot b - K_{\CP}) \, = \, W \contrac (b^* \theta_{\CP} - \theta_{\CP}) \, = \, 0
\end{equation*}
This means that the function $K_{\CP} \cdot b - K_{\CP}$ is constant along the leaves of the polarization $\CP$. Making use of the relations (\ref{e:Kp}) we find that the value of this constant on a Bohr-Sommerfeld leave $\La_{\bq}  : \, k \om(W_i, X) = q_i \pmod{k}$ is
\begin{equation*}
K_{\CP}(b X) - K_{\CP}(X) \, = \, \frac{1}{2} \frac{^t \bq}{k} (A^{-1} B) \frac{\bq}{k}
\end{equation*}
If $x \in \La_{\bq}$ then $b(x) \in \La_{ ^t \! A^{-1} \bq }$, since from the relations $b W_i = A_{ji} W_j$ and from $k \om(b W_i, b X) = k \om(W_i, X) = q_i \pmod{k}$,  it follows that $k \om(W_j, b X) = ( ^t A^{-1})_{ij} q_i \pmod{k}$ which is the equation defining the Bohr-Sommerfeld orbit $\La_{ ^t \! A^{-1} \bq }$.
Having this in view we find that 
\begin{align*}
b \cdot s_{\bq} \, &= \, \te^{\frac{\pi \ti}{k} \, ^t \bq (A^{-1}B) \bq} \, s_{ ^t \! A^{-1} \bq} \\ 
b \cdot \de_{\bq} \, &= \, \de_{ ^t \! A^{-1} \bq} \notag
\end{align*}
Taking $b$ equal to the element $\al$ or to $\be$ leads to the expressions
given in (\ref{p:Sprepr}).
\end{proof}
The projective representation of $Sp(\CZ)$ on $\BP \CH$ constructed above coincides with the projective representation of this group defined in \cite{G} on the vector space of theta functions at level $k$ obtained through the quantization of a symplectic torus in a holomorphic (K\"{a}hler) polarization. 

\msk


\subsection{The Maslov index and the representation of the integer metaplectic group} \label{ss:rMaslov}

Let $\rho : Mp(\CV) \ra Sp(\CV)$ be the $2:1$ covering homomorphism from the metaplectic group to the symplectic group of $(\CV, \om)$.
The integer metaplectic group $Mp(\CZ)$ is defined as the preimage $\rho^{-1}(Sp(\CZ))$.

The group $Mp(\CZ)$ acts on the 4-fold cover $\tlag$ of $\lag$ and preserves the subset of elements in $\tlag$ which project to $\lag$ onto rational Lagrangian subspaces.
Having in view the definitions and notations of Sect.\ref{s:Maslov}, we associate to each element $\tilde{b} \in Mp(\CZ)$ covering  the element $b \in Sp(\CZ)$ the unitary map of Hilbert spaces $\tilde{b} : \CH_{\tilde{L}_{\CP}} \ra \CH_{\tilde{b} \tilde{L}_{\CP}}$ defined by setting $\tilde{b} = b : \CHP \ra \CH_{b \CP}$. 
The composition of $\tilde{b}$ with $F_{\tilde{L}_{\CP} , \tilde{b} \tilde{L}_{\CP}} : \CH_{\tilde{b} \tilde{L}_{\CP}} \ra \CH_{\tilde{L}_{\CP}}$ is a unitary operator 
$U_{\CP}(\tilde{b}) = F_{\tilde{L}_{\CP} , \tilde{b} \tilde{L}_{\CP}} \circ \tilde{b}$ on  $\CH_{\tilde{L}_{\CP}} = \CHP$. 
Moreover, for any two polarizations $\CP_1, \CP_2$ and for any $\tilde{L}_{\CP_1} , \, \tilde{L}_{\CP_2} \in \tlag$ covering $L_{\CP_1}$ and $L_{\CP_2}$, we have 
$F_{\tilde{b} \tilde{L}_{\CP_2}, \tilde{b} \tilde{L}_{\CP_1} } \circ \tilde{b} \, = \, \tilde{b} \circ F_{\tilde{L}_{\CP_2}, \tilde{L}_{\CP_1} }$.
Together with Theorem (\ref{t:transit}) this implies that
\begin{equation*} 
U_{\CP}(\tilde{b}) U_{\CP}(\tilde{b}') = U_{\CP}(\tilde{b} \tilde{b}'), \quad \text{for any } \,\tilde{b}, \tilde{b}' \in Mp(\CZ).
\end{equation*}
The map $\tilde{b} \mapsto U_{\CP}(\tilde{b}) $ defines a unitary representation of the group $Mp(\CZ)$ on the Hilbert space $\CHP$. 
One can easily check that $F_{\tilde{L}_{\CP_2}, \tilde{L}_{\CP_1} } \circ U_{\CP_1}(\tilde{b}) \circ F_{\tilde{L}_{\CP_1}, \tilde{L}_{\CP_2} } = U_{\CP_2}(\tilde{b})$.
Thus we can state:
\begin{thm}
The group $Mp(\CZ)$ has a unitary representation on the Hilbert space $\CH$ associated through quantization to $(\CT, k \om, \tlag)$.
\end{thm}
\nin This representation of $Mp(\CZ)$ is the analog to the present setting of the Shale-Weil representation of the real metaplectic group on the infinite dimensional Hilbert space associated through quantization to a symplectic linear space \cite{GS,LV}.

In order to give an explicit expression of the unitary operators representing the elements of $Mp(\CZ)$ on the Hilbert space we seek a description of $Mp(\CZ)$  in terms of generators.
To that end we make use of the following results from (\cite{Go}, \S 3).
The notations are those of Sect.\ref{s:Maslov}.
Fix a Lagrangian plane $L \in \lag$. Then we have: \\
$\bullet \,$ The $q$-fold cover $Lag_q(\CV)$ can be identified with the set $[ \lag \times \BZ_{2q}]_L = \{ (L', \la) \mid \la \equiv g - \dim (L \cap L') \pmod{2} \}$. 
The latter can be equipped with a topology which makes it homeomorphic to $Lag_q(\CV)$ and the natural projection $[ \lag \times \BZ_{2q}]_L \ra \lag$ continuous.
If $\hat{L}_1=(L_1, \la_1), \, \hat{L}_2=(L_2, \la_2) \in [ \lag \times \BZ_{2q}]_L$ then the Maslov index $\mu_{2q} (\hat{L}_1,\hat{L}_2)$ is given by
\begin{equation}
\mu_{2q}(\hat{L}_1,\hat{L}_2) \, = \, \la_1 - \la_2 + \tau(L, L_1, L_2)  \pmod{2q}\end{equation}
$\bullet \,$ The $q$-fold covering group $Sp_q(\CV)$ of the symplectic group $Sp(\CV)$ can be identified with the group \text{$[ Sp(\CV) \times \BZ_{4q} ]_L \, = \, \{ (b,z) \mid z \equiv \mu(b_{\infty} L_{\infty},L_{\infty}) \pmod{4} \}$}, where $b_{\infty} \in Sp_{\infty}(\CV)$ and $L_{\infty} \in Lag_{\infty}(\CV)$ are arbitrary elements projecting onto $b$ and $L$, respectively.
We note that due to (\ref{t:cover}) and (\ref{p:propMi}) we have, for every $r \in \BZ$, $\: \mu((\al^r b_{\infty}) L_{\infty}, L_{\infty}) = \mu(b_{\infty} L_{\infty},L_{\infty}) + 4r \,$ and $\, \mu(b_{\infty} (\be^r L_{\infty}), \be^r L_{\infty}) = \mu(b_{\infty} L_{\infty},L_{\infty})$, so that the expression $\mu(b_{\infty} L_{\infty},L_{\infty}) \pmod{4}$ is independent of the choices $b_{\infty}$ and $L_{\infty}$ projecting onto $b$ and $L$.
The composition law in the group $[ Sp(\CV) \times \BZ_{4q} ]_L$ is given by
\begin{equation} \label{e:mult}
(b,z) \, (b',z') \, = \, \bigl( b b' \, , z + z' + \tau(L, b L, b b' L) \pmod{4q} \bigr)
\end{equation}
The group $[ Sp(\CV) \times \BZ_{4q} ]_L$ can be equipped with a topology which makes it isomorphic to $Sp_q(\CV)$ as topological groups and for which the projection map $[ Sp(\CV) \times \BZ_{4q} ]_L \ra Sp(\CV)$ is continuous. \\  
$\bullet \,$ The group $Sp_q(\CV) \cong [ Sp(\CV) \times \BZ_{4q} ]_L$ acts transitively on the space $Lag_{2q}(\CV) \cong [Lag(\CV) \times \BZ_{4q}]_L$ by
\begin{equation} \label{e:action}
(b,z) \, (L',\la) \, = \, \bigl( b L', z + \la + \tau(L, b L, b L') \pmod{4q} \bigr)
\end{equation}

Therefore, for the integer metaplectic group $Mp(\CZ) \subset Mp(\CV) = Sp_2(\CV)$ which interests us, we get the identification
\begin{equation*}
Mp(\CZ) \cong  [ Sp(\CZ) \times \BZ_8 ]_L = \{ (b,z) \mid z \equiv \mu(b_{\infty} L_{\infty},L_{\infty}) \pmod{4} \}
\end{equation*}
This will enable us to give a description of the generators of $Mp(\CZ)$. 
Let us first make the following remarks.
\begin{rmk}
Let $(b, z_b)$ be an element of $[ Sp(\CZ) \times \BZ_8 ]_L$ such that $b L = L$ and let $L^+$ denote the Lagrangian $L$ with a choice of orientation.
Since $Sp(\CV)$ acts on the space $Lag_2(\CV)$ of oriented Lagrangian planes, we can easily deduce, using (\ref{e:action}) and the equality $\mu(b_{\infty} L_{\infty},L_{\infty}) \equiv g - \dim(bL \cap L) \pmod{2} \equiv 0 \pmod{2}$, that: \\
1) If $b L^+ = L^+$ then $z_b \equiv 0 \pmod{4}$. \\
2) If $b L^+ \neq L^+$ then $z_b \equiv 2 \pmod{4}$.
\end{rmk}
\begin{rmk}
Let us consider the elements $(\ga,z_{\ga})$ and $(\be_1, z_{\be_1})$ of $[ Sp(\CZ) \times \BZ_8 ]_L$. The elements  $\ga$ and $\be_1$ belong to the set (\ref{e:Spgen}) of generators of $Sp(\CZ)$ and are described by the matrices $\ga = \left( \begin{smallmatrix} 0 & I \\ -I & 0 \end{smallmatrix} \right)$ and $\be_1 = \left( \begin{smallmatrix} I & I \\ 0 & I \end{smallmatrix} \right)$, with respect to a symplectic basis $(W_i \, ;W_i^{\perp})$ of $\CV$ with $L = \text{span}_{\BR} \{ W_1,\dots ,W_g \}$. 
Since $\ga L \cap L = 0$ and $\mu(\ga_{\infty} L_{\infty}, L_{\infty}) \equiv g - \dim(\ga L \cap L) \pmod{2} \equiv g \pmod{2}$ we must have the equality $z_{\ga} \equiv  g \pmod{2}$.
Applying the multiplication rule (\ref{e:mult}) and using the fact that $(\ga \be_1)^3 = e$ we find that $( (\ga,z_{\ga}) (\be_1,z_{\be_1}))^3 = ( e, 3 z_{\ga} + 3 z_{\be_1} + g)$.
By the previous remark we must have $z_{\be_1} \equiv 0 \pmod{4}$ and $3 z_{\ga} + 3 z_{\be_1} + g \equiv 0 \pmod{4}$.
This leads to the conclusion that $z_{\ga} \equiv g \pmod{4}$.
\end{rmk}
In view of the above remarks, we can see that the group $Mp(\CZ) \cong [ Sp(\CZ) \times \BZ_8 ]_L$ has as generators the elements:
\begin{alignat}{2} \label{e:Mpgen}
\tilde{\vep} &= (e,4)  & & \\
\tilde{\be} &= (\be,0), &  & \quad \text{ where } \; \be = \left( \begin{smallmatrix} I & B \\ 0 & I \end{smallmatrix} \right) , \; B \in M(g,\BZ), \, ^t B =B \notag \\
\tilde{\al} &= (\al,z_{\al}), & &  \quad \text{ where } \; \al =\left( \begin{smallmatrix} A & 0 \\ 0 & ^t A^{-1} \end{smallmatrix} \right) , \; A \in GL(g,\BZ), \notag \\
& & &\quad \text{ and }\; z_{\al} = \begin{cases} 0, &\, \text{ if } \; \det A >0 \\ 2, &\, \text{ if } \; \det A <0 \end{cases} \notag \\
\tilde{\ga} &= (\ga,z_{\ga}), & & \quad \text{where} \; \ga =\left( \begin{smallmatrix} 0 & I \\ -I & 0 \end{smallmatrix} \right) , \notag \\
& & &\quad \text{ and } \; z_{\ga} =i, \text{ if } g\equiv i \pmod{8}, \; i=0,1,2,\dots,7 \notag \\
\tilde{e} &= (e,0), & & \quad \text{ the identity element } \notag
\end{alignat}
where the matrix representations of the corresponding elements of $Sp(\CZ)$ are with respect to the  symplectic basis $(W_i \, ; W_i^{\perp})$.
We note that 
\begin{equation*}
(\tilde{\ga} \tilde{\be}_1)^3 = \begin{cases} \tilde{e},  &\; \text{ if } \, g= \text{even} \\
\tilde{\vep},  &\; \text{ if } \, g= \text{odd}
\end{cases}
\end{equation*}
where $\tilde{\be}_1 =(\be_1,0)$ with $\be_1 = \left(\begin{smallmatrix} I & I \\ 0 & I \end{smallmatrix} \right)$, and that
\begin{equation*}
(\tilde{\ga})^4 = \begin{cases} \tilde{e}, &\; \text{ if } \, g= \text{even} \\
\tilde{\vep}, &\; \text{ if } \, g= \text{odd}.
\end{cases}
\end{equation*}

In particular, for $g=1$, we have the group $Sp(\CZ) = SL(2,\BZ)$ with standard generators the matrices $S=\left( \begin{smallmatrix} 0 &1 \\ -1 &0 \end{smallmatrix} \right)$ and $T=\left( \begin{smallmatrix} 1 &1 \\ 0 &1 \end{smallmatrix} \right)$ subject to the relations $(ST)^3=I,\, S^4 =I$. The integer metaplectic group $Mp(\CZ)=Mp(2,\BZ)$ has generators $\tilde{S}=(S,5),\; \tilde{T}=(T,0),\; \tilde{\vep}=(e,4)$ and the identity $\tilde{e}=(e,0)$, satisfying the relations $(\tilde{S} \tilde{T})^3 = \tilde{e}$ and $(\tilde{S})^4 = \tilde{\vep}$.

\msk

Let us now take $L$ to be a rational Lagrangian plane $L_{\CP}$ which determines a polarization $\CP$ on the torus $(\CT = \CV/\CZ, k \om)$.
We choose an element $\tilde{L}_{\CP}$ in $\tlag \cong [\lag \times \BZ_8]_{L_{\CP}}$ projecting onto $L_{\CP}$ and look for the representation of the group $Mp(\CZ)$ on the Hilbert space $\CH_{\tilde{L}_{\CP}}$.
To each element $\tilde{b} = (b,z_b)$ in $Mp(\CZ) \cong [ Sp(\CZ) \times \BZ_8]_{L_{\CP}}$ there corresponds a unitary operator on  $\CH_{\tilde{L}_{\CP}} = \CHP$
\begin{equation*}
U_{\CP}(\tilde{b})= F_{\tilde{L}_{\CP}, \tilde{b} \tilde{L}_{\CP}} \circ \tilde{b} \, = \te^{- \frac{\pi \ti}{4} \tilde{\mu}(\tilde{L}_{\CP}, \tilde{b} \tilde{L}_{\CP})} F_{\CP, b \CP} \circ b \, = \te^{\frac{\pi \ti}{4} z_b} U_{\CP}(b)
\end{equation*}
Using the results (\ref{p:Sprepr}) for the representation $U_{\CP}$ of $Sp(\CZ)$ on $\CHP$ we obtain:
\begin{prop}
The representation of $Mp(\CZ)$ on $\CHP$ is described by the following unitary operators representing the generators
\begin{align*}
U_{\CP}(\tilde{\vep}) \, \si_{\bq} \, &= \, \te^{\pi \ti} \: \si_{\bq} \\
U_{\CP}(\tilde{\al}) \, \si_{\bq} \, &= \, \te^{\frac{\pi \ti}{4} z_{\al}} \: \si_{^t \! A^{-1} \bq} \\
U_{\CP}(\tilde{\be}) \, \si_{\bq} \, &= \, \te^{\frac{\pi \ti}{k} \,{ ^t \bq} B \bq} \: \si_{\bq} \\
U_{\CP}(\tilde{\ga}) \, \si_{\bq} \, &= \, k^{- \frac{g}{2}} \, \te^{\frac{\pi \ti}{4} z_{\ga}} \,  \sum_{\bq_1 \in (\BZ /k \BZ)^g} \, \te^{\frac{2\pi \ti}{k} \,{ ^t \bq} \bq_1} \: \si_{\bq_1}
\end{align*}
\end{prop} 
In particular, for $g=1$, we obtain the representation of $Mp(2,\BZ)$ described by
\begin{align*}
U_{\CP}(\tilde{\vep})_{q q'} \, &= \, \te^{\pi \ti} \: \de_{q q'}  \\
U_{\CP}(\tilde{T})_{q q'} \, &= \, \te^{\frac{\pi \ti}{k} \, q^2} \: \de_{q q'} \\
U_{\CP}(\tilde{S})_{q q'} \, &= \,  k^{- \frac{1}{2}} \, \te^{\frac{5 \pi \ti}{4} }  \, \te^{\frac{2\pi \ti}{k} \,q q'} 
\end{align*}

\msk


\subsection{The metaplectic correction and the representation of the integer metaplectic group}

We recall that in Sect.\ref{s:metapl} we fixed an invariant symplectic frame $(U_1, \dots , U_g ;$$V_1, \dots ,V_g)$ of $(\CT,k \om)$ which determined the trivialization $SP(\CT) = \CT \times Sp(2g,\BR)$ and chose the metaplectic frame bundle to be the product bundle $MP(\CT)= \CT \times Mp(2g,\BR)$. With respect to the symplectic basis $(U_i \, ; V_i)$ of $(\CV, k\om)$ we identify the group $Sp(\CV)$ with the real symplectic group $Sp(2g,\BR)$, the subgroup $Sp(\CZ) \subset Sp(\CV)$ with the integer symplectic group $Sp(2g,\BZ)$, the group $Mp(\CV)$ with the real metaplectic group $Mp(2g,\BR)$ and the subgroup $Mp(\CZ) \subset Mp(\CV)$ with the integer metaplectic group $Mp(2g,\BZ)$. 

Any element $b$ in $Sp(\CZ)$ induces a bundle automorphism $b : SP(\CT) \ra SP(\CT)$ which, under the above mentioned identifications, is described by the left-action
\begin{equation*}
(x,c) \in \CT \times Sp(2g,\BR) \overset{b}{\lra} (b(x), bc) \in \CT \times Sp(2g,\BR)
\end{equation*}
Any $\tilde{b} \in Mp(\CZ)$ acts on $(\CT,k \om)$ as a symplectic diffeomorphism through its projection $b = \rho(\tilde{b}) \in Sp(\CZ)$.
The action lifts to an automorphism $\tilde{b} : MP(\CT) \ra MP(\CT)$ of the metaplectic frame bundle described by 
\begin{equation*}
(x,\tilde{c}) \in \CT \times Mp(2g,\BR) \overset{\tilde{b}}{\lra} (b(x), \tilde{b} \tilde{c}) \in \CT \times Mp(2g,\BR)
\end{equation*}
\begin{rmk}
Let $MP^{\ep}(\CT)$ denote a nontrivial metaplectic frame bundle determined by a homomorphism $\ep : \CZ \ra \BZ/2\BZ$ (see Sect.\ref{s:metapl}).
Then, only the elements $\tilde{b}$ in $Mp(\CZ)$ with projection $b= \rho(\tilde{b})$ to $Sp(\CZ)$ satisfying $\ep(b W) = \ep(W)$, for all $W \in \CZ$, define bundle automorphisms $\tilde{b} : MP^{\ep}(\CT) \ra MP^{\ep}(\CT)$.
\end{rmk}
For any polarization $\CP$ of $(\CT, k \om)$, an element $b \in Sp(\CZ)$ induces a bundle morphism $b: GL(\CP) \ra GL(b \CP)$ between the bundles of linear frames of $\CP$ and of $b \CP$.
The bundle automorphism $\tilde{b} : MP(\CT) \ra MP(\CT)$ determined by $\tilde{b} \in Mp(\CZ)$ induces a morphism of metalinear frame bundles $\tilde{b} : ML(\CP) \ra ML(b \CP)$.
Hence there is an induced morphism between the corresponding half-forms bundles, $\tilde{b} : (\text{Det} \, {\mathcal P}^{*} )^{\frac{1}{2}} \ra (\text{Det} \, (b {\mathcal P})^{*} )^{\frac{1}{2}} $. 
In conclusion, any $\tilde{b}$ in $Mp(\CZ)$ gives rise to a map of line bundles
\begin{equation*}
\begin{CD}
\CL \otimes (\text{Det} \, {\mathcal P}^{*} )^{\frac{1}{2}} @>\tilde{b}>> \CL \otimes (\text{Det} \, (b {\mathcal P})^{*} )^{\frac{1}{2}} \\
@VVV     @VVV \\
\CT @>b>> \CT
\end{CD}
\end{equation*}
The induced map on the space of sections sends a section $s \otimes \mu$ of $\CL \otimes (\text{Det} \, {\mathcal P}^{*} )^{\frac{1}{2}}$ to the section $(\tilde{b} \cdot s) \otimes (\tilde{b} \cdot \mu)$ of $\CL \otimes (\text{Det} \, (b {\mathcal P})^{*} )^{\frac{1}{2}}$ defined at any $x \in \CT$ by
\begin{equation}
(\tilde{b} \cdot s)(x) \, = \, (b \cdot s)(x) \, = \, b \, s(b^{-1}(x))
\end{equation}
and
\begin{equation} \label{e:Mpd}
(\tilde{b} \cdot \mu)^{\sharp} \, = \, \mu^{\sharp} \cdot \tilde{b}^{-1} \, ,
\end{equation}
where $\mu^{\sharp} : ML(\CP) \ra \BC$ is the $ML(g,\BR)$-equivariant function associated to $\mu$. 
Therefore we have a unitary map of Hilbert spaces $\tilde{b} : \tilde{\CH}_{\CP} \ra \tilde{\CH}_{b \CP}$.
Its composition with the operator $\tilde{F}_{\CP, b \CP} : \CH_{b \CP} \ra \CH_{\CP}$ arising from the BKS pairing gives a unitary operator $\tilde{U}_{\CP}(\tilde{b}) = \tilde{F}_{\CP, b \CP}  \circ \tilde{b}$ on $\tilde{\CH}_{\CP}$.
Theorem (\ref{t:transit}) together with the fact that $\tilde{F}_{b \CP_2, b \CP_1} \circ \tilde{b} = \tilde{b} \circ \tilde{F}_{\CP_2 \CP_1}$, for any two polarizations $\CP_1$ and $\CP_2$, imply that 
\begin{equation}
\tilde{U}_{\CP}(\tilde{b}) \, \tilde{U}_{\CP}(\tilde{b}') \, = \, \tilde{U}_{\CP}(\tilde{b} \tilde{b}'), \quad \text{ for any } \, \tilde{b}, \tilde{b}' \in Mp(\CZ) \, .
\end{equation}
Thus, the assignment $\tilde{b} \mapsto \tilde{U}_{\CP}(\tilde{b})$ defines a unitary representation of $Mp(\CZ)$ on $\tilde{\CH}_{\CP}$ and we can state:
\begin{thm}
The group $Mp(\CZ)$ has a unitary representation on the Hilbert space $\tilde{\CH}$ associated through quantization to $(\CT, k \om, Mp(\CT))$.
\end{thm}
We determine the explicit form of the representation.
We fix a polarization $\CP$ and a symplectic frame $(W_i \, ; W_i^{\perp})$ of $(\CT, \om)$ as in (\ref{e:Zbasis}), with $\CP = \text{span}_{\BR} \{ W_1,\dots,W_g \}$.
Then we take $(U_i = W_i \, ; V_i = \frac{1}{k} W_i^{\perp})$ as the symplectic frame for $(\CT, k \om)$ with respect to which we choose to trivialize the symplectic frame bundle and to define the metaplectic frame bundle.
Recall from (\ref{ss:met}) that the linear frame bundle $GL(\CP)$ is a subbundle of the bundle of Lagrangian frames $\CF l(\CT) = SP(\CT)/N = \CT \times Sp(2g,\BR)/N$. 
The frame field $\underline{W}_x = (W_1,\dots, W_g)_x \in GL(\CP)_x$ is identified to $(x,eN) \in \CT \times Sp(2g,\BR)/N$. 
The metalinear frame bundle $ML(\CP)$ is a subbundle of $\tilde{\CF l}(\CT) = MP(\CT)/N = \CT \times Mp(2g,\BR)/N$ covering $GL(\CP)$.
We define the section $\de$ of $ (\text{Det} \, {\mathcal P}^{*} )^{\frac{1}{2}}$ by setting $\de^{\sharp}(\underline{\tilde{W}}) = 1$, where $\underline{\tilde{W}} \in ML(\CP)$ is the metalinear frame field projecting onto $\underline{W}$ given by $\underline{\tilde{W}}_x = (x, \tilde{e}N)$.
Then, as shown in Sect.\ref{ss:met}, we have a uniquely determined standard unitary basis $\{ \si_{\bq} = s_{\bq} \otimes \de_{\bq} \}$ of the Hilbert space $\tilde{\CH}_{\CP}$.

Let $\tilde{\be} \in Mp(\CZ)$ be the element projecting onto the element $\be \in Sp(\CZ)$ given in (\ref{e:Spgen}) and such that $\tilde{\be} \cdot \underline{\tilde{W}} = \underline{\tilde{W}}$. 
Then $\tilde{\be} \cdot s_{\bq} = \be \cdot s_{\bq}$ and $\tilde{\be} \cdot \de_{\bq} = \de_{\bq}$, so that $\tilde{\be}$ is represented on $\tilde{\CH}_{\CP}$ by the unitary operator
\begin{equation*}
\tilde{U}_{\CP}(\tilde{\be})  \, \si_{\bq} \, = \, \te^{\frac{\pi \ti}{k} \,{ ^t \bq} B \bq} \, \si_{\bq} \,.
\end{equation*}
Let $\tilde{\vep} \in Mp(\CZ)$ be the element covering the identity $e \in Sp(\CZ)$ and such that $\tilde{\vep} \neq \tilde{e}$.
Then $\tilde{\vep}$ belongs also to the metalinear group $ML(g,\BZ) \subset Mp(2g,\BZ)$ and $\chi(\tilde{\vep}) = \te^{\pi \ti} = -1$.
According to (\ref{e:Mpd}) we have
\begin{equation*}
(\tilde{\vep} \cdot \de)^{\sharp}(\underline{\tilde{W}}) = \de^{\sharp}(\tilde{\vep}^{-1} \cdot \underline{\tilde{W}}) = 
\de^{\sharp}( \underline{\tilde{W}} \cdot \tilde{\vep}^{-1}) 
= \chi(\tilde{\vep})^{-1} \de^{\sharp}(\underline{\tilde{W}}) = \chi(\tilde{\vep})^{-1}  = \te^{\pi \ti} \,.
\end{equation*}
Thus $\tilde{\vep} \cdot \de = \te^{\pi \ti} \, \de$ and we have
\begin{equation*}
\tilde{U}_{\CP}(\tilde{\vep})  \, \si_{\bq} \, = \, \te^{\pi \ti} \, \si_{\bq} \,.
\end{equation*}
Let $\tilde{\al} \in Mp(\CZ)$ be an element projecting onto $\al = \left( \begin{smallmatrix} A & 0 \\ 0 & ^t A^{-1} \end{smallmatrix} \right)$ of (\ref{e:Spgen}).
Then, since $\al \in GL(g,\BZ) \subset Sp(2g,\BZ)$, we have $\tilde{\al} \in ML(g,\BZ) \subset Mp(2g,\BZ)$. 
Therefore we get
\begin{equation*}
(\tilde{\al} \cdot \de)^{\sharp}(\underline{\tilde{W}}) = \de^{\sharp}(\tilde{\al}^{-1} \cdot \underline{\tilde{W}}) = 
\de^{\sharp}( \underline{\tilde{W}} \cdot \tilde{\al}^{-1}) 
= \chi(\tilde{\al})^{-1} \de^{\sharp}(\underline{\tilde{W}}) = \chi(\tilde{\al})^{-1}
\end{equation*}
That is, $\tilde{\al} \cdot \de = \chi(\tilde{\al})^{-1} \de$ and we have
\begin{align*}
\tilde{\al} \cdot s_{\bq} \, &= \, \al \cdot s_{\bq} \, = \, s_{^t \! A^{-1} \bq} \\
\tilde{\al} \cdot \de_{\bq} \, &= \, \chi(\tilde{\al})^{-1} \, \de_{^t \! A^{-1} \bq} \notag
\end{align*}
Thus we get
\begin{equation*}
\tilde{U}_{\CP}(\tilde{\al})  \, \si_{\bq} \, = \, \chi(\tilde{\al})^{-1} \, \si_{^t \! A^{-1} \bq}
\end{equation*}
Let $\tilde{j} \in Mp(\CZ)$ be the element introduced in Sect.\ref{ss:forms} which projects onto $j = \left(\begin{smallmatrix} 0 & -I \\ I & 0 \end{smallmatrix} \right) \in Sp(\CZ)$.
Note that $j \ga = e$ where $\ga$ is the element defined in (\ref{e:Spgen}).
Let $\CP' = j \CP$ and $\underline{W}' = j \cdot \underline{W}$.
Let $\underline{\tilde{W}}' = \tilde{j} \cdot \underline{\tilde{W}}$ and define the section $\de'$ of $(\text{Det} \CP^{'*})^{\frac{1}{2}}$ by $\de^{' \sharp}(\underline{\tilde{W}}') = 1$.
Then $\tilde{j} \cdot \de_{\bq} = \de'_{\bq}$ and $\tilde{j} \cdot s_{\bq} = s'_{\bq}$.
Therefore we get
\begin{equation*}
\begin{split}
\tilde{U}_{\CP}(\tilde{j})  \, \si_{\bq} \, &= \, (\tilde{F}_{\CP, j \CP} \circ \tilde{j}) (s_{\bq} \otimes \de_{\bq}) \, = \, \tilde{F}_{\CP, j \CP} (s'_{\bq} \otimes \de'_{\bq}) \\
&= \, k^{- \frac{g}{2}} \, \te^{-\frac{\pi \ti}{4} g} \,  \sum_{\bq_1 \in (\BZ /k \BZ)^g} \, \te^{-\frac{2\pi \ti}{k} \,{ ^t \bq} \bq_1} \: \si_{\bq_1}
\end{split}
\end{equation*}
where the last equality follows from the general formulas (\ref{e:trmat}) and (\ref{e:phase}).
Let $\tilde{\ga} \in Mp(\CZ)$ be the element covering $\ga$ and satisfying $\tilde{\ga} \tilde{j} = \tilde{e}$. 
Then $\tilde{U}_{\CP}(\tilde{\ga}) \tilde{U}_{\CP}(\tilde{j}) = I$ and therefore
\begin{equation}
\tilde{U}_{\CP}(\tilde{\ga})  \, \si_{\bq} \, = \, k^{- \frac{g}{2}} \, \te^{\frac{\pi \ti}{4} g} \,  \sum_{\bq_1 \in (\BZ /k \BZ)^g} \, \te^{\frac{2\pi \ti}{k} \,{ ^t \bq} \bq_1} \: \si_{\bq_1}
\end{equation}
A comparison of the above results to those of the previous subsection shows that we obtain the same representation of $Mp(\CZ)$.

\bsk

\nin {\bf Acknowledgements}. I would like to thank Prof. Dan Freed for helpful conversations and suggestions regarding the material contained in this paper. 
 
\bsk



\end{document}